\documentclass[10pt]{article}
\usepackage{amssymb,amsmath,amsfonts,mathrsfs,upgreek,textgreek,bbold,relsize} 
\usepackage{mathabx}

\usepackage[nottoc]{tocbibind}

\usepackage{graphicx} 

\usepackage[colorlinks=true]{hyperref}
\usepackage{comment}

\renewcommand\comment[1]{{\iffalse #1 \fi}}

\newcommand\sopra[2]{\genfrac{}{}{0pt}{}{#1}{#2}}

\newtheorem{theorem}{Theorem}[section]
\newtheorem{corollary}{Corollary}[section]

\newtheorem{lemma}{Lemma}[section]
\newtheorem{proposition}{Proposition}[section]

\newtheorem{assumption}{Assumptions}[section]

\newtheorem{remark}{Remark}[section]
\newtheorem{notation}{Notation}[section]
\newtheorem{definition}{Definition}[section]
\newcommand{\diam}{{\,\rm diam\,}}

\newcommand{\Lip}{{\,\rm Lip}}
\newcommand{\io}{{\infty}}

\newcommand{\real}{ {\mathbb R}   }
\newcommand{\torus}{ {\mathbb T}    }
\newcommand{\integer}{ {\mathbb Z}   }
\renewcommand{\natural}{ {\mathbb N}   }
\newcommand{\complex}{ {\mathbb C}   }

\newcommand{\cA}{ {\mathcal A}   }

\newcommand{\cF}{ {\mathcal F}   }
\newcommand{\cM}{ {\mathcal M}   }
\newcommand{\cH}{ {\mathcal H}   }

\newcommand{\cT}{ {\mathcal T   } }

\renewcommand{\Im}{\, {\rm Im}\,}
\renewcommand{\Re}{\, {\rm Re}\,}

\newcommand{\eproof}{\qed}

\newcommand\beq[1]{ \begin{equation}\label{#1} }
\newcommand{\eeq}{ \end{equation} }
\newcommand{\beqno}{ \[ }
\newcommand{\eeqno}{ \] }
\newcommand\beqa[1]{ \begin{eqnarray} \label{#1}}
\newcommand{\eeqa}{ \end{eqnarray} }
\newcommand{\beqano}{ \begin{eqnarray*} }
\newcommand{\eeqano}{ \end{eqnarray*} }

\newcommand\dfn[1]{ \begin{definition}\label{#1} }
\newcommand\edfn{ \end{definition} }
\newcommand\ass[1]{ \begin{assumption}\label{#1} }
\newcommand\eass{ \end{assumption} }

\newcommand\notat[1]{ \begin{notation} \label{#1} 
 }
\newcommand\enotat{\end{notation}}

\newcommand\rem{\begin{remark} 
\rm 
}
\newcommand\erem{\end{remark} 
}

\newcommand{\proof}{\noindent{\bf Proof\ }}
\newcommand\equ[1]{{\rm (\ref{#1})}}
\newcommand{\nl}{{\smallskip\noindent}}
\newcommand{\giu}{{\medskip\noindent}}
\newcommand{\Giu}{{\bigskip\noindent}}

\newcommand{\qed}{\hskip.5truecm
\vrule width 1.7truemm height 3.5truemm depth 0.truemm
\par\Giu}
\newcommand{\qedeq}{\hskip.5truecm
\vrule width 1.7truemm height 3.5truemm depth 0.truemm}

\newcommand\bolla{{\footnotesize $\bullet\,$}}

\newcommand\casitwo[4]{ \left\{  \begin{array}{ll}
 {#1} & \mbox{ {\rm if} ${#2}$} \\
 {#3} & \mbox{ {\rm if} ${#4}$}
 \end{array} \right.}

\newcommand{\x}{\xi}

\newcommand{\e}{\varepsilon}

\renewcommand{\a }{\alpha }
\renewcommand{\b }{\beta }
\newcommand{\s }{\sigma }
\newcommand{\ii }{{\rm i} }
\renewcommand{\d }{\delta }

\newcommand{\g }{\gamma}
\newcommand{\f }{\varphi}

\renewcommand{\l }{\lambda }

\renewcommand{\t }{\tau }
\renewcommand{\o }{\omega }
\newcommand{\C}{\mathbb{C}}
\newcommand{\Z}{\mathbb{Z}}
\newcommand{\z }{z}

\newcommand{\normadue}{|}

\newcommand{\ttM}{\mathtt M}
\newcommand{\ttD}{\mathtt D}
\newcommand{\cd}{{c_{{}_n}}}
\newcommand{\Nf}{\mathtt N}

\renewcommand\AA{{\rm A}}
\newcommand\hAA{\hat\AA}

\newcommand\ttx{{\mathtt x}}
\newcommand\tty{{\mathtt y}}

\newcommand\ttq{{\mathtt q}}
\newcommand\ttp{{\mathtt p}}
\newcommand\ttG{{\mathtt G}}

\newcommand{\cgot}{\mathfrak c}

\newcommand{\Ggot}{{\mathfrak G}}

\def\R{\mathbb R}
\def\T{\mathbb T}

\def\dalla{\updelta_{\!\rm o}}  
\def\balla{\upgamma}    

\def\const{{\, \rm const\, }}
\def\dst{\displaystyle}
\def\bks{\, \backslash\, }
\DeclareMathOperator{\meas}{meas}
\newcommand\eqby[1]{\stackrel{\equ{#1}}{=}}

\newcommand\proiezione{\, { \pi}}

\newcommand\modulo{|}
\newcommand\ham{\mathtt{H}}
\newcommand\hamk{\ham_k}

\newcommand\hamo{\ham_{\rm o}}
\newcommand\hamsec{\overline\ham}

\newcommand{\Hpend}{{\ham}_{\flat}}
\newcommand{\bHpend}{{\bar\ham}_{\flat}}
\newcommand{\Hsharp}{{\ham}_k} 

\newcommand\Gf{\ttG_{{}_k}}
\newcommand\bGf{\bar\ttG_{{}_k}}

\newcommand\lalla{\upmu}
\newcommand\lella{\upmu_{\rm o}}
\newcommand\loge{\uplambda} 
\newcommand\logeo{\loge_{\rm o}}
\newcommand\blogemax{\bar\loge_{{}_{\rm max}}}
\newcommand\logesharp{\loge_{{}_\sharp}}
\newcommand\logeuno{\loge_{{}_1}}

\newcommand\sa{\theta} 

\newcommand\Gm{{\mathtt G}}
\newcommand\GO{{\bar{\mathtt G}}} 

\newcommand{\tettao}{\vartheta_{\rm o}}

\newcommand\chr\ro

\newcommand\chs{{\check \so}}

\newcommand{\K}{{\mathtt K}}
\newcommand{\KO}{{\mathtt K}_{\rm o}}
\newcommand{\suca}{\upepsilon}
\newcommand\morse{\upbeta}
\newcommand{\Ro}{{\mathtt R}} 
\newcommand{\ro}{{\mathtt r}}
\newcommand{\so}{{\mathtt s}}  
\newcommand{\Eb}{\mathtt E_\flat}

\newcommand{\rkam}{{\mathfrak r}}
\newcommand{\skam}{{\mathfrak s}}

\newcommand\act{I}  
\newcommand\ang{\f}  

\newcommand\Fiuno{\Phi_{\!{}_1}} 
\newcommand\cFiuno{\check\Phi_{\!{}_1}} 
\newcommand\Fidue{\Phi_{\!{}_2}}
\newcommand\cFidue{\check\Phi_{\!{}_2}}
\newcommand\Fitre{\Phi_{\!{}_3}}
\newcommand\cFitre{\check\Phi_{\!{}_3}}

\newcommand\Fiq{\Upphi^i}                       

\newcommand\cFiq{\check\Upphi^i}           

\newcommand\Bu{{\mathcal B}}
\newcommand\Buc{\mathcal B_{\!{}_{\rm near}}}
\newcommand\Bua{\mathcal B_{\!{}_{\rm far}}}

\newcommand{\acci}{a}
\newcommand{\bacci}{b}

\newcommand{\logemax}{{\loge_{{}_{\rm max}}}}

\newcommand\cin{\upnu}
\newcommand\cins{\upnu_{\!{}_k}}

\newcommand{\ts}{\textstyle}

\newcommand{\fico}{\eta} 

\newcommand{\fproj}{\pi_{\!{}_{\integer k}}}

\usepackage{appendix}

\newcommand\Appendix{\hyperref[appendicite]{Appendix}}

\newcommand\gen{{\cal G}^n}
\def\genk{{\cal G}^n_{K}}
\def\genKO{{\cal G}^n_{\KO}}
\def\genK{{\cal G}^n_{\K}}

\renewcommand{\cgot}{\mathtt c}

\newcommand\noruno[1]{  |#1|_{{}_1} }

\renewcommand\ln{\log}
\newcommand\st{:\ }

\newcommand\hol{{\mathbb B}}
\newcommand{\Bns}{{\mathbb B}^n_s}
\newcommand{\Gns}{{\mathbb G}^n_{s}}

\newcommand\pk{{\proiezione^\perp_k}}
\newcommand\pko{{\proiezione^\perp_k}}

\newcommand\DD{{\rm B}}

\newcommand\Rz{{\mathcal R}^0}	
\newcommand\Ru{{\mathcal  R}^1}	
\newcommand\Ruk{{\mathcal  R}^{1,k}}	
\newcommand\Rd{{\mathcal  R}^2}	
\newcommand\Rzt{\widetilde{\mathcal R}^0}
 
\newcommand\Rukt{\widetilde{\mathcal  R}^{1,k}}

\newcommand\Bi{B^i}

\newcommand\Bik{B^i_k}
\newcommand\Bui{{\mathcal B}^i}
\newcommand\Buik{{\mathcal B}^i_k}

\newcommand\cO{{\mathcal O}}
\newcommand\cOr{{\mathcal O}_{\!\varrho}}
\newcommand\DDD{\mathscr D}

\newcommand\bfEs{{{\mathtt E}_{{}_{\scalebox{0.7}{\rm *}}}^{i}}}           

\newcommand\bfI{{\rm{I}}}

\renewcommand\diamond{{{}\!\varstar}}

\newcommand\htk{\hat {\mathcal Q}_k}
\newcommand\hzk{\hat h_{\!{}_k}}

\newcommand\gu{{\bar{\rm a}}}

\newcommand\mmu{{\bar{\rm m}}}
\newcommand\mtt{{\mathtt m}}

\newcommand\lippo{{\mathtt b}}
\newcommand\vulva{\upzeta}           
\newcommand\dt{\tilde\updelta}			
\newcommand\ddi\updelta			
\newcommand\dds{\updelta_\sharp}  
\newcommand\ddb{\bar\updelta}		

\newcommand\bn{{\bar n}}

\newcommand\ta{{\mathtt g}}

\newcommand\rhol{{\rho_{{}_\loge}}}
\newcommand\prhol{{\rho_{{}_{\!\loge}}'}}
\newcommand\sil{{\sigma_{\!{}_\loge}}}
\newcommand\rstar{\mathtt r_{\!{}_{\scalebox{0.7}{\rm *}}}}
\newcommand\rhs{{\rho_{{}_\diamond}}}
\newcommand\shs{{\s_{{}_\diamond}}}

\newcommand\giuno{{\mathtt g}_{{}_1}}  
\newcommand\gitre{{\mathtt g}_{{}_3}} 
\newcommand\gstar{{\mathtt g}_{{\,}_\diamond}}
\newcommand\Istar{j_{\!{}_{\gstar}}} 
\newcommand\Istarin{j_{\!{}_{-\gstar}}}

\newcommand\chk{\chi_{{}_k}}

\newcommand\xish{{\xi_{{{}_\sharp}}}}

\newcommand\vallo{\mathcal I}

\newcommand\cdr{{c}_{{}_{\star}}}
\newcommand\ttcnr{{\mathtt C}_{\rm o}}

\newcommand\cteo{{\mathfrak c}}

\newcommand\cteos{\mathfrak c_\star}

\newcommand\cgotu{{\mathfrak c_1\,}}
\newcommand\cteot{{\mathfrak c_{{{}_0}}}}   
\newcommand\cgoth{\hat\cteo}

\newcommand\hcteot{{\hat{\mathfrak c}_{{{}_0}}}} 
\newcommand\ckam{{\rm C_*}}

\newcommand\bfco{{\bf c_{{{}_0}}}}

\newcommand\bfcast{{\bf c}_{{}_\diamond}}

\newcommand\bfc{{\bf c}}
\newcommand\hbfc{{\hat \bfc}}
\newcommand\bfcu{{\bf c}_{{{}_2}}}
\newcommand\bfcd{{\bf c}_{{{}_1}}}
\newcommand\bfct{{\bf c}_{{{}_3}}}
\newcommand\bfcq{{\bf c}_{{{}_4}}}
\newcommand\bfcc{{\bf c}_{{{}_6}}}
\newcommand\bfcs{{\bf c}_{{{}_5}}}
\newcommand\bfcst{{\bf c}_{{{}_7}}}
\newcommand\bfcot{{\bf c}_{{{}_8}}}

\newcommand\itcast{c_{{\!{}_\coasterisk}}}

\newcommand\itcu{c_{{{}_1}}}
\newcommand\itcd{c_{{{}_2}}} 
\newcommand\itct{c_{{{}_3}}} 

\newcommand\ttcs{{\mathtt c}_s}       
\newcommand\ttcq{\bar{\mathtt C}_{{{}_1}}}
\newcommand\ttCu{{\mathtt C}_{{{}_1}}}

\newcommand\bcteo{\bar{\mathfrak c}\,}

\DeclareMathOperator{\id}{id}

\newcommand\hik{h^i_k}
\newcommand\fik{f^i_k}

\newcommand\Gdag{\Ggot_{\!{}_\dag}}
\newcommand\Ggoto{\Ggot_{\!{}_0}}

\makeatletter
\newcommand{\pushright}[1]{\ifmeasuring@#1\else\omit\hfill$\displaystyle#1$\fi\ignorespaces}
\newcommand{\pushleft}[1]{\ifmeasuring@#1\else\omit$\displaystyle#1$\hfill\fi\ignorespaces}
\makeatother

\title{\bf 
Singular KAM Theory
}

\begin{document}

\author{ 
\footnotesize L. Biasco  \& L. Chierchia
\\ \footnotesize Dipartimento di Matematica e Fisica
\\ \footnotesize Universit\`a degli Studi Roma Tre
\\ \footnotesize Via della Vasca Navale, 84 - 00146 Roma, Italy
\\ {\footnotesize luca.biasco@uniroma3.it, luigi.chierchia@uniroma3.it}
\\ 
}

\maketitle

\hfill{\footnotesize \it Dedicated to the memory of our friend and colleague Walter L. Craig}


\begin{abstract}\noindent
The question of the total measure of invariant tori in analytic, nearly--integrable Hamiltonian systems is considered.
In 1985, Arnol'd, Kozlov and Neishtadt, in 
the Encyclopaedia of Mathematical Sciences \cite{AKN1}, and in subsequent editions, conjectured that in  $n=2$ degrees of freedom the measure of the non torus set of general analytic nearly--integrable systems away from critical points is exponentially small with the size $\e$ of the perturbation, and that for $n\ge 3$ the measure is, in general,  of order $\e$ (rather than $\sqrt\e$ as predicted by classical KAM Theory).
\\
In the case of generic natural Hamiltonian systems, we prove lower bounds on the measure of primary and secondary invariant tori, which are in agreement, up to a logarithmic correction, with the above conjectures.
\\
The proof is based on a new {\sl singular} KAM theory, particularly designed to study   analytic properties in neighborhoods of the secular separatrices generated by the perturbation at simple resonances.

\nl

\end{abstract}

\section*{Introduction}

Classical KAM Theory\footnote{\cite{Ko54}, \cite{A63}, \cite{Mo62}, \cite{A63b}, \cite{Mo68};  for a divulgative account, see  \cite{KAMstory}.} 
deals with the persistence of Lagrangian invariant tori of integrable Hamiltonian systems under the effect of  small perturbations. 
In the early 1980's it was clarified that an analytic integrable system, which is Kolmogorov non--degenerate (i.e., such that the action--to--frequency map is a local diffeomorphism), preserves, under a perturbation of size $\e>0$, all its Diophantine Lagrangian invariant tori in a bounded domain  {\sl up to a set of measure proportional to\footnote{\cite{Laz}, \cite{Nei}, \cite{P82}, \cite{Sva}} $\sqrt{\e}$}.
 In fact, this estimate cannot be improved, since trivial examples -- such as a classical pendulum with  Hamiltonian $\frac{p^2}2+\e \cos q$ -- show that, in bounded domains,  {\sl the measure of 
the complement of persistent primary tori} (i.e., tori which are a deformation of integrable ones) {\sl is  exactly proportional to the square root of the perturbing function} --  the rest of the phase space being filled, in the case of the pendulum,  by  {\sl secondary tori} (curves) enclosed by the pendulum separatrix. In fact,  positive measure sets of secondary Lagrangian tori  (i.e., tori, which are not a smooth deformation of integrable ones) appear  in general nearly--integrable systems, for example, near elliptic equilibria
(\cite{MNT}).

\nl
The natural question is therefore: {\sl What is  the measure of all Lagrangian tori in general nearly--integrable analytic Hamiltonian systems?}

\nl
In 1985 Arnold, Kozlov and Neishtadt, motivated by the exponentially small splitting of separatrices in general systems with two degrees of freedom,  conjectured that\footnote{Compare \cite[p.~189]{AKN1} and \cite[Remark 6.17, p. 285]{AKN}.} 

\nl
``It is natural to expect that in a generic (analytic) system with two degrees of freedom and with frequencies that do not vanish simultaneously the total measure of the `non-torus' set corresponding to all the resonances is exponentially small.''

\nl
In \cite{AKN},  again Arnold, Kozlov and Neishtadt, arguing on the basis of a simple rescaling argument in neighbourhoods of double resonances\footnote{\label{qqq}From p. 285 of \cite{AKN}: 
``Indeed, the $O(\sqrt{\e})$--neighbourhoods of two resonant surfaces intersect in a domain of measure $\sim \e$. In this domain, after the partial averaging taking into account the resonances under consideration, normalizing the deviations of the ``actions'' from the resonant values by the quantity $\sqrt\e$, normalizing time, and discarding the terms of higher order, we obtain a Hamiltonian of the form $1/2(Ap, p) + V(q_1, q_2)$, which does not involve a small parameter. Generally speaking, for this Hamiltonian there is a set of measure $\sim 1$ that does not contain points of invariant tori. Returning to the original variables we obtain a ``non--torus'' set of measure $\sim \e$.''}, 
conjectured that

\nl
``It is natural to expect that in a generic system with three or more degrees of freedom the measure of the `non--torus' set has order $\e$.''

\nl
In this paper,  we develop a `singular KAM theory' for generic analytic {\sl nearly--integrable natural systems}, apt to deal, in particular,  with the construction of maximal KAM tori that live exponentially close to  the separatrices appearing near simple resonances, which are singularities of  the action--angle variables of the integrable secular (averaged) systems. As a consequence,  
we can prove lower bounds on the total measure of KAM tori, which are 
in agreement with the above conjectures  up to a logarithmic correction $|\log \e|^\g$.
We announced  these results  in 2015 in \cite{BClin} (see also \cite{BClin2}), and it goes without saying that  to complete proofs took {\sl much} longer than we thought.

\nl
The reason for dealing with the special class of nearly--integrable natural systems, namely, 
Hamiltonian systems on $\real^n\times \torus^n$ (endowed with the standard symplectic form $dy\wedge dx$) with Hamiltonian given by 
\beq{ham}\ts
\ham(y,x;\e):=\frac12 |y|^2 +\e f(x)\,,\quad (y,x)\in\real^n\times\torus^n\,,\ 0<\e<1\,,\quad \big(|y|^2:=y\cdot y:=\sum_j|y_j|^2\big)\,,
\eeq
is twofold. On one side, this choice allows to avoid  technical unessential details, which would make even heavier the already highly technical methods. On the other hand, and more importantly, it allows to formulate the generiticity condition (whose definition is part of the problem) in a simple way, singling out a suitable class of {\sl generic analytic} potentials $f$'s (Definition~\ref{sicuro}), which guarantees, in particular,  a uniform behaviour of the {\sl secondary nearly--integrable structure}  at simple resonances with high modes  (compare Open Problems, (i) in \S~\ref{omnia}).

\giu
Let us informally discuss the overall picture.  

\nl
Analogously to what is done  in Nekhoroshev theory (compare, e.g.,  \cite{P93}, and \cite{AKN} for general information), 
fixed a maximal size of resonances $\K$ to be taken into account\footnote{Actually, we will need to consider {\sl two} orders of resonances
$\KO$ and $\K>\KO$; but for the purpose of this introduction we call them both $\K$. }, 
one covers the action space with three sets: a non--resonant set $\Rz$, a $\sqrt\e \K^c$--neighbourhood\footnote{In this introduction, we indicate with `$c$' various different constants, which are independent of $\e$. In general, keeping track of the quantities, on which the various constants appearing in singular KAM Theory depend, is a somewhat important matter (for example, from the constructive point of view) and we try to devote some care to it; compare, e.g., Remark~($\rm R_5$)  in \S~\ref{omnia}.} $\Ru$ of simple resonances,  
and a neighborhood  $\Rd$ of double (and higher) resonances. 
Eventually, the number of resonances $\K$ is taken as a suitable function of $\e$ tending to $+\io$ as $\e\to 0$ (e.g., $\K\sim1/\e^c$, or $\K\sim|\log \e|$).\\
The set $\Rd$, which has measure proportional to 
$\e \K^c$,  is a {\sl non perturbative set} in the sense
that the dynamics ruled by $\ham(y,x;\e)$ on $\Rd\times \torus^n$ is essentially equivalent to the dynamics of the 
parameter free Hamiltonian $\frac12 |y|^2+ f(x)$ (compare the argument given by Arnold, Kozlov and Neishtadt, reproduced in  footnote~\ref{qqq} above). Therefore, no further perturbative analysis on the set $\Rd\times \T^n$ is possible. 
\\
On the non resonant phase space $\Rz\times \T^n$, after high order averaging, classical  KAM theory yields the existence of primary maximal KAM tori up to a set of measure $O(\sqrt\e e^{-c\K})$.

\nl
The main game has then to be played on the  simple--resonance neighborhood $\Ru\times\T^n$.
\\
$\Ru$ is defined as union of  sets $\Ruk$, which are  $\sqrt\e\K^c$--neighborhoods  
of simple resonances $\{y\big|\, y\cdot k=0\}$ with  $k\in\integer^n$ and co--prime entries. On $\Ruk$ high--order averaging theory can be applied 
so as to remove, up to order $\e e^{-c\K}$, the angle dependence, apart from the resonant combination $k\cdot x$, obtaining a 
 symplectically conjugated   real analytic Hamiltonian  of the form
 \beq{hamk}\ts
\hamk(y,x) =  {\dst \frac{|y|^2}2}+\e \big( g^k_{\rm o}(y)+ 
g^k(y,k\cdot x) +
f^k (y,x)\big)\,, \qquad f^k\sim e^{-c \K}\,.
\eeq
Now, all these Hamiltonian systems labelled by the simple--resonance index $k$ ($|k|\le \K$), have a {\sl secondary} (secular) {\sl near--integrability structure}, as, disregarding the exponentially small terms $f^k$, they are Arnol'd--Liouville integrable,   depending effectively only on one resonant angle $\ttx_1=k\cdot x\in \T^1$. \\
Then, the plan is obvious: Put all these systems into their Arnol'd--Liouville action--angle variables, check twist (i.e., Kolmogorov's non--degeneracy), and apply KAM so as to obtain Lagrangian primary and secondary tori (with different topologies; compare Remark~\ref{tortore}). 
\\
However, a considerable series of problems arise in trying to carry out such a plan. Let us try to highlight the most important points.

\nl
First of all, as already mentioned, the resonance cut--off $\K$ will go to $+\io$ as $\e\to 0$ and therefore one has to deal, {\sl de facto}, with infinitely many Hamiltonian systems and unless there is some uniform way of treating them, there is no hope. The idea, here, has been suggested in \cite{BCnonlin} and refined in \cite{BCuni}: The  {\sl secular  Hamiltonians} $\hamsec_k(\tty,\ttx_1)$, i.e., the
 integrable Hamiltonians in \equ{hamk} obtained disregarding
$f^k$ and setting $\ttx_1=k\cdot x$,   
are one--degree--of--freedom Hamiltonians, with external parameters, and  with potentials $g^k(y,\ttx_1)$, which are close to the projections  $(\fproj f)(\ttx_1)$ over the Fourier modes proportional to $k$ of the potential $f(x)$; compare \equ{kproj} below. 
Now, one can show that for high Fourier modes $|k|>\Nf$ ($\Nf$ suitable but independent of $\e$), $(\fproj f)(\ttx_1)$ behaves generically as a shifted cosine\footnote{Compare item (iii) in Theorem~\ref{normalform} below.} 
$$
2|f_k|\e \cos(\ttx_1 +\sa_k)
$$
for a suitable $\sa_k\in\real$; where `generically' means that $f$ belongs to a suitable class of generic real--analytic potential, whose Fourier coefficients, for large $k$'s with co--prime entries,  behave as $e^{-|k|s}|k|^{-n}$ for a suitable $s>0$.
Incidentally, to obtain such a result, one has to use  a non--standard averaging theory, allowing for essentially {\sl no analyticity loss} in the angle variables; for more information on this point, see the Introduction in \cite{BCnonlin}.
\\
Analytic properties of the action--angle variables for the pendulum are quite well known, and this is encouraging (and it was also the basis for the optimistic 2015 announcement \cite{BClin}). \\
However, for low modes $|k|\le \Nf$,  the secular leading potentials $(\fproj f)(\ttx_1)$ are, in general, quite arbitrary functions, and one needs, therefore, a general holomorphic, quantitative theory of action--angle variable for one--degree--of freedom systems containing parameters. Such a theory is discussed in \cite{BCaa23} for a special class of real analytic Hamiltonians -- called there {\sl Generic Standard Form Hamiltonians} -- given by
\beq{stendardo}
\Hpend(p,q_1)= \big(1+ \cin(p,q_1)\big) p_1^2  +\Gm(\hat p, q_1)\,,
\eeq
where $p=(p_1,\hat p)$, $p_1$ is the momentum conjugated to the angle $q_1$ and  $\hat p=(p_2,...,p_n)$ are the `external parameters'; see Definition~\ref{morso} for specifications. In particular, the properties of the energy--to--action functions are discussed in the limit as the energy approaches the critical values (i.e., the energy levels of the hyperbolic points and the associated separatrices): It turns out that such functions have the form
\beq{yesterday}
\act_1(E_{\rm crit} \pm \suca z)=a(z) + b(z)\, z\log z
\eeq
where $\suca$ is a suitable reference energy, $E_{\rm crit}$ is a fixed  critical energy level of some equilibrium of the secular system, 
 $a$ and $b$ are analytic functions of $z$ (and, of course, everything depends on other  $(n-1)$ dumb action; compare Theorem~\ref{glicemiak} below).
This representations will  play a crucial r\^ole in studying the twist of the secular Hamiltonians at simple resonances in their Arnold-Liouville action--angle variables.
\\
Now, one can prove (\cite{BCuni}) that {\sl all secular Hamiltonians  $\hamsec_k$ can be put into  standard form as in \equ{stendardo}}, so that  the main rescaling  properties are controlled by one single parameter $\upkappa$, which is independent of $\e$ {\sl and} $k$ (compare Theorem~\ref{sivori} below). The draw back of this  uniformization is that the symplectic transformations performing the task  {\sl are not well defined in the fast angle--variables} $(\ttx_2,...,\ttx_n)$, and preserves periodicity only in the resonant angle $\ttx_1$.

\nl
This is the starting point of this paper. 

\nl
In \S~\ref{GS} we show how to overcome the homotopy problem of the uniformization of \cite{BCuni}:
Exploiting the particular group structure of the various symplectic transformations involved, we  show that, introducing  a special {\sl ad hoc} symplectic `semi--conjugacy', one can indeed obtain, for all $|k|\le \K$,   well defined symplectic action--angle maps $\upphi_k^i$, which conjugate the original Hamiltonian $\ham$ in \equ{ham} on 
$\Ruk\times \T^n$ to the nearly--integrable form
$$
\cH^i_k:=\ham\circ \upphi_k^i(\act,\ang)
=
\hik(\act)+
\e \fik(\act,\ang)\,,\qquad  \fik\sim  e^{-c \K}\,,
$$
where $i$ labels the various regions in which the phase spaces of the secular Hamiltonians $\hamsec_k$ are split by their separatrices (compare Theorem~\ref{garrincha} below). 

\nl
As in all KAM applications, the main problem is to prove (a suitable) {\sl non--degeneracy} of the frequency map   
$\act\to\o=\partial_\act \hik$.\\
 It should be clear from the context, that the original non--degeneracy of $\ham|_{\e=0}$ plays a little r\^ole here, as the action--structure depends on analytic properties of the secular potentials. Indeed, it is a fact that {\sl in the phase space of standard Hamiltonians \equ{stendardo}  there are, in general,  points where the twist vanishes}. For instance, points of vanishing twist appear {\sl always} in regions bounded by two separatrices (with different energy); but they appear also in very simple examples with only one separatrix near the elliptic equilibrium enclosed by the  separatrix,  like, e.g.,  in the case of the Hamiltonian 
 $$\ts p_1^2 + \cos q_1 - \frac18 \cos(2q_1)\,;
 $$
 compare Remark~\ref{bejazet} below. Furthermore (and more seriously), when the distance in energy from the separatrices goes to zero, the problem becomes a {\sl singular perturbation problem} with dramatic singularities.
\\
Therefore, entirely new methods have to be developed in order to prove that {\sl the measure where the twist vanishes is actually exponentially small in the whole phase space of all Hamiltonians $\cH^i_k$}. This is the main result of the paper; compare  the Twist Theorem~\ref{zelensky} in  \S~\ref{thisistheend}. 
\\
The proof of the Twist Theorem  is based on two different approaches according 
to whether one considers  regions {\sl far} from separatrices or  regions {\sl close} to separatrices. 
\\
In regions far from separatrices the analysis is significantly simpler, since it is partly perturbative. In such a case, one fist proves that the (normalized) second derivative of the action--to--energy functions are non--degenerate (i.e., at each point of their domains, some derivative is different from zero); then, uniform estimates can be worked out and, using standard tools from the theory of Diophantine approximations (\cite{pyartli}, \cite{E}), one can show that $\eta$--sub--levels  of the twist determinant have measure smaller than $\eta^c$, which easily yields the claim.
\\
The real heart of the matter is the analysis of the twist in regions close to separatrices, where no perturbative arguments can be used, nor uniform estimates hold. 
The proof, in this case, rests on the construction of a suitable differential operator with non--constant coefficients, which, exploiting in a subtle way the analytic structure \equ{yesterday},  can be shown not to vanish on a suitable regularization of the twist determinant. This is good enough to prove that the twist determinant is non--degenerate also near separatrices, and to conclude the proof of the Twist Theorem.

\nl
At this point, choosing carefully the various free parameters of the game, a suitable KAM Theorem 
(Theorem~\ref{KAM} below) yields the existence of maximal primary and secondary KAM tori, which fill the complementary phase set of $\Rd\times \T^n$ up to a very small set $\cA$. 
How small is $\cA$ -- which dynamically is very rich and where, e.g., Arnol'd diffusion can take place -- 
depends on how big is chosen the order $\K$ of resonances considered. 
For example,   if $\K$ is chosen as $|\log \e|^2$, then 
 $\meas(\cA)$ will be almost--exponentially small (i.e., smaller than any power of $\e$), 
while $\meas(\cA)$ is actually  exponentially small in $1/\e^c$, if $\K$  is taken to be an inverse power of $\e$; compare Remark ($\rm R_2$) in \S~\ref{omnia}.

\Giu

\giu
{\footnotesize
{\bf Acknowledgments} During the decade needed to complete this paper the authors benefited from conversations with
A.~Delshams, M.~Guardia, V.~Kaloshin, S.~Kuksin,  
C.~Liverani, A.~Neishtadt, L.~Niederman, 
M.~Procesi, V.~Rom Kedar, T.~Seara, A.~Sorrentino and D.~Treshev.\\
The authors acknowledge the encouragement of Caterina Biasco.
\\
This work is supported by the Italian grant PRIN 2022FPZEES.
}

{\small 
\tableofcontents
}

\section{Results, Remarks, and Open Problems}\label{omnia}

In order to state  the main results of this paper, we  recall a few  standard definitions. 

\nl 
\bolla{\sl Maximal KAM tori}:
A  set $\cT\subset {\cal M}=B\times\real^n$ is called a {\sl maximal KAM torus} for a real analytic  Hamiltonian $H:{\cal M} \to \real$ if there exist a real analytic embedding $\phi:\torus^n\to {\cal M}$ and a Diophantine frequency vector\footnote{A vector $\o\in\real^n$ is called  {\sl Diophantine} if there exist $\a>0$ and $\t\ge n-1$ such that $\dst |\o\cdot k| \ge \a/|k|_{{}_1}^\t$, for any non vanishing integer vector $k\in\integer^n$, where $\noruno{k}:=\sum |k_j|$.} $\o\in\real^n$ such that $\cT=\phi(\torus^n)$, and for each $z\in\cT$, $\Phi^t_H(z)=\phi(x+ \o t)$, where\footnote{In particular, maximal KAM tori are minimal invariant invariant sets for $\Phi^t_H$.}
$x=\phi^{-1}(z)$ and $t\to \Phi^t_H(z)$ denotes the standard Hamiltonian flow governed by $H$ starting at $z\in {\cal M}$.

\nl 
\bolla{\sl Generators of 1d maximal lattices}:
Let $ \integer^n_\varstar$ be the set of integer vectors $k\neq 0$ in $\integer^n$ such that the  first non--null  component is positive:
\beqno
 \integer^n_\varstar:=
 \big\{ k\in\integer^n:\ k\neq 0\ {\rm and} \ k_j>0\ {\rm where}\ j=\min\{i: k_i\neq 0\}\big\}\,.
 \eeqno
$\gen$  denotes the set of {\sl generators of 1d maximal lattices} in $\integer^n$, namely, the set of  vectors $k\in  \integer^n_\varstar$ such that the greater common divisor (gcd)  of their components is 1:
\beq{gen}
\gen:=\{k\in \integer^n_\varstar:\ {\rm gcd} (k_1,\ldots,k_n)=1\}\,;
\eeq
for $K\ge 1$, we set\footnote{As usual $\noruno{k}:=\sum_{j=1}^n|k_j|$.}:
\beq{genk}
\genk:=\gen \cap \{\noruno{k}\leq K \}\,.
\eeq

\nl 
\bolla{\sl 1d Fourier projectors}:
Given  a zero--average  real analytic 
periodic function 
$$
f: x\in\torus^n:=\real^n/(2\pi \Z^n)\mapsto f(x):=\sum_{\integer^n\bks\{0\}}
f_k e^{\ii k\cdot x}
$$
and fixed  a vector $k\in \integer^n\bks\{0\}$,
we denote by $\fproj f$ the (real analytic) periodic function of  
{\sl one variable} $\sa\in\torus$ given by
\beq{kproj}
\sa\in\torus\mapsto \fproj f (\sa):=\sum_{j\in\integer} f_{jk} e^{\ii j\sa}\,.
\eeq
Notice that one has the following (unique) decomposition:
\beqno
f(x)=
\sum_{k\in \gen} \fproj f (k\cdot x)\,.
\eeqno

\nl 
\bolla{\sl Resonances}:
Given $k\in\gen$,  a resonance ${\cal R}_k$ 
with respect to the free Hamiltonian $\frac12 |y|^2$ is the
 set $\{y\in \real^n: y\cdot k=0\}$. 
We call
${\cal R}_{k,\ell}$  a {\sl double resonance} if  ${\cal R}_{k,\ell}={\cal R}_k\cap {\cal R_\ell}$ with $k$ and $\ell$ in $\gen$ linearly independent; the order of a double resonance is given by $\max\{\noruno{k},\noruno{\ell}\}$.

\nl 
\bolla{\sl Morse functions with distinct critical values}:
A  $C^2$--function of one variable $\sa\to F(\sa)$ is a Morse function if its critical points are non--degenerate, i.e., $F'(\sa_0)=0\implies F''(\sa_0)\neq 0$; `distinct critical values' means that if 
$\sa_1\neq\sa_2$  are distinct critical points, then $F(\sa_1)\neq F (\sa_2)$.

\nl 
\bolla{\sl Banach spaces of real analytic periodic functions}:
For $s>0$ and  $n\in\natural=\{1,2,3...\}$,
consider  the Banach space of zero--average
real analytic periodic functions on $\torus^n$ 
with finite norm       
\beq{cripta}
 \|f\|_s:=\sup_{k\in \integer^n}  |f_k| e^{\noruno{k}s}\,,
\eeq
and denote by $\Bns$ {\sl its closed unit ball}.

\nl
Now we are ready to introduce a suitable generic class of potentials $f$ and 
state the main result of this paper concerning the typical dynamics of nearly--integrable natural systems with Hamiltonian $\ham$ as in \equ{ham}.

 \dfn{sicuro} {\bf (The  class of potentials $\Gns$)}
We  denote by $\Gns$ the subset  of functions $f\in \Bns$  such that the following two properties 
hold:
\beqa{P1}
&&
\varliminf_{\sopra{\noruno{k}\to+\io}{k\in\gen}} |f_k| e^{\noruno{k}s} \noruno{k}^n>0\,,
\\
&& 
 \forall \ k\in\gen\,,\ \fproj f\ {\rm is \ a\ Morse\ function\ with \ distinct \ critical\ values}\,.
 \label{P2}
\eeqa
\edfn

\rem
(i) For natural systems 
-- whose natural phase space is  $\R^n\times\T^n$ -- any ball in action space can be transformed (by translation and rescaling) to the unit ball; thus, rescaling time and renaming the smallness parameter one can always restrict oneself to study the Hamiltonian $\ham$  in \equ{ham} on the 
unit ball in $\real^n$ and with $\|f\|=1$. 

\nl
(ii) The requirement in \equ{P2} that the projections $\fproj f$ have different critical values is not really necessary, but it is generic and it simplifies the proof.
\erem

\subsubsection*{Main Results}

\begin{theorem}\label{prometeo} Let  $n\ge 2$, $s>0$, $\DD:=\{y\in \R^n\st  |y| <1\}$, $\g:=11 n +4$,   $f\in \Gns$ with $\|f\|_s=1$.
Then, 
there exist a constant
$\cteo>1$,
such that for all $\K$ and $\e>0$  satisfying
\beq{torino}
\K\ge \cteo\,,\qquad \quad
\e\, \K^\gamma\le 1\,,   
\eeq
the following holds.
There exist three sets 
$\Rd\subseteq \DD$, 
${\cA}\subseteq \DD\times\torus^n$, 
$\cT\subseteq \real^n\times\torus^n$
such that:

\begin{itemize}
\item[\rm (i)] $\DD\times\torus^n\subseteq (\Rd\times\torus^n)\cup {\cA}\cup \cT$; 

\item[\rm (ii)] $\Rd$ is a neighborhood 
 of double resonances of order smaller than $\K$ satisfying the measure estimate  
 $$\meas\Rd\le \cdr\, \e\, \K^\gamma\,,$$
 where  $\cdr$ is a suitable constant depending only on $n$; 

\item[\rm (iii)]  $\cA$ is exponentially small with  respect to $\K$: 
$$\meas\cA\le e^{-  \K/\cteo}\,;$$

\item[\rm (iv)]  $\cT$ is union  of maximal KAM tori for the natural Hamiltonian $\ham(y,x;\e):=\frac12 |y|^2 +\e f(x)$.
\end{itemize}

\end{theorem}
An immediate corollary of this theorem is that the measure of the `non--torus set' for $\ham$ does not exceed  $O(\e|\log\e|^\g)$:

\begin{corollary}\label{prometeobis} 
Under the  assumptions of Theorem~\ref{prometeo}, there exists $0<\e_{\rm o}<1$ such that for $\e<\e_{\rm o}$, 
all  points in $\DD\times \T^n$
 lie on a maximal KAM torus for $\ham$, except for a subset  whose measure is bounded by
 $
 \bcteo\, \e |\ln \e|^{\g}$ where 
 $\bcteo=1+(2\pi)^n \cdr \cteo^\gamma$.
 \end{corollary}
In particular, this corollary
implies  the  theorem announced in \cite[p. 426]{BClin}.\\
Indeed, from items (i), (ii) and (iii) of Theorem~\ref{prometeo}, there follows
\beq{lapalisse}
\meas\big((\DD\times\torus^n)\bks \cT\big)\le (2\pi)^n \cdr  \, \e \K^\g + e^{- \K/\cteo }\,,
\eeq
and Corollary~\ref{prometeobis} follows immediately by choosing 
$\K:=\cteo|\ln\e|$ (and $\e_{\rm o}$ small enough). 

\nl
The two degrees of freedom is special: in this case the only double resonance is the origin and one can take as $\Rd$ a disk of measure $\e^a$
with any $0<a<1$ getting a set of KAM tori of exponential density in the complementary of $\Rd\times\torus^2$.  This is the content of  next corollary (compare, also, 
 \cite{BC-DCDS}).

\begin{corollary}\label{prometeotris} Let the assumptions of Theorem~\ref{prometeo} hold and let $n=2$.
Then, there exists $0<\e_{\rm o}<1$, such that
for  $\e<\e_{\rm o}$ and   $0<a<1$, 
   all  points in the set $\{y\in \DD: |y|>\e^{a/2}\}\times\T^2$
 lie on a maximal KAM torus for $\ham$ in \equ{ham}, except for
 an exponentially
small set of measure bounded by 
$
e^{-1/(2\cteo\e^{\hat a})}
$,
with $\hat a:=(1-a)/24$.
\end{corollary}

\subsection*{Remarks and Open Problems}

First, we briefly discuss the class of potentials $\Gns$, for which the above results hold; for more information on $\Gns$ and complete proofs, see  
\cite[Sect. 2]{BCuni}.

\nl
It is very simple to give explicit examples of functions in $\Gns$, a prototype being\footnote{Recall the definitions given in \equ{gen}, \equ{cripta} and \equ{kproj}.}
\beq{puro}
f(x):=2 \sum_{k\in\gen} e^{-\noruno{k}s} \cos k\cdot x\,,
\eeq
which (as it is trivial to verify) satysfies
\beqno
\|f\|_s=1\,,\qquad 
\lim_{\sopra{\noruno{k}\to+\io}{k\in\gen}} |f_k| e^{\noruno{k}s} \noruno{k}^n=+\io\,,\qquad
\fproj f (\sa)=2 e^{-\noruno{k}s} \cos \theta\,.
\eeqno
The class of potentials $\Gns$  is quite general from various points of view. For example, the following result (proven in  \cite[Sect. 1]{BCuni}) holds:

\begin{proposition}
{\rm (i)}  The  class  $\Gns$ contains an open and dense set in $\Bns$.

\nl
{\rm (ii)} Let ${\rm F}$ denote the weighted Fourier isometry
 $ {\rm F}:f
\in  \Bns \to \big\{f_k e^{|k|_{{}_1}s}\big\}_{k\in {\integer^n_\varstar}}\in \ell^\io({\integer^n_\varstar})$.
Then 
 ${\rm F}(\Gns)$ is a set of probability 1 with respect to the standard product probability measure on ${\rm F}(\Bns)$.
\end{proposition}

\nl
We shall not use this proposition  here, however, we shall  often use a  suitable quantitative characterization of $\Gns$. \\
To state such a characterization one needs to make quantitative the notion of Morse functions with different critical values and to introduce a uniform Fourier cut--off
function, depending on the dimension $n$ and on two parameters $0<\d\le 1$ and $s>0$.

\begin{definition}\label{buda} {\bf ($\b$--Morse functions)}
Let $\b>0$.
 $F\in C^2(\T,\R)$  is called  $\b$--Morse,  if 
\beq{ladispoli}
\min_{\sa\in\T} \ \big( |F'(\sa)|+|F''(\sa)|\big) 
\geq\b \,,
\quad
\min_{i\neq j } 
|F(\sa_i)-F(\sa_j)|
\geq \b \,, 
\eeq
where $\sa_i\in\T$ are the  critical points of $F$.
\end{definition}

\begin{definition}\label{duda} {\bf (The cut--off function $\Nf$)}
Given $0<\d\le 1$ and $n,s>0$ define the following   `Fourier cut--off function':
 \beq{enne}
\Nf=\Nf(\d;s,n):=2\, \max\Big\{1\,,\,\frac1s \,   \log \frac\cd{s^n\, \d}\Big\} \,,\qquad \cd:= 2^{44}\ (2n/e)^n\,.
\eeq 
\end{definition}
Then, the following elementary result holds:

\begin{lemma}
 \label{telaviv} Let $n,s>0$. Then, 
 $f\in\Gns$ if and only if $f\in\Bns$ and  there exist $0<\d\le 1$ and $\b>0$ such that 		
\beq{P1+}
|f_k|\geq \d \noruno{k}^{-n}\, e^{-\noruno{k} s}\,,\qquad \forall \ k\in\gen\,,\ \noruno{k}\ge\Nf\,,
\eeq 
 \beq{P2+}
\fproj f\ {\rm is \ \b\!-\!Morse}\,,\qquad \,\,\, \ \ \forall \ k\in\gen\,,\ \noruno{k}\le \Nf\,.
\eeq
 \end{lemma}
The  proof of this lemma is given in \Appendix.

\nl 
{\bf Remarks}

\nl
($\rm R_1$) By item (i) in Theorem \ref{prometeo}, we see that
the phase space of a generic nearly--integrable natural system can be covered by two `small' sets, 
namely, $\Rd\times\torus^n$ and $\cA$,
and one `large' set, namely, $\cT$. Such sets exhibit quite different dynamics and satisfy the following properties.
\begin{itemize}
\item[\bolla]
$\Rd\times\torus^n$  is contained in a tubolar  neighborhood   of double resonances
${\cal R}_{k,\ell}\cap \DD$ of order not exceeding $\K$; compare \equ{due} and  \equ{defi}  below.
As mentioned in the Introduction, the set $\Rd\times\torus^n$ contains  a set of measure $\e$ where {\sl the dynamics
 is not perturbative} in the sense that, as long as  trajectories lie in this set, the dynamics is ruled by an effective Hamiltonian having no small 
 parameters.

\item[\bolla] The set $\cA$, which has measure $\sim e^{- \K/\cteo}$,  is  dynamically very interesting. For example, it is where the asymptotic manifolds of lower dimensional tori break up (`exponentially small splitting of separatrices') giving rise, e.g.,  to local horse shoe dynamics and, most likely,  to Arnol'd diffusion\footnote{This statement has not yet been proven in any generality in the analytic class; for references to Arnol'd diffusion see, e.g.,  \cite{A}, \cite{CG}, \cite{BBB},
\cite{T}, \cite{DH}, \cite{Z}, 
\cite{Ma}, \cite{BKZ}, \cite{DLS}, \cite{CdL22}, \cite{GLS}, \cite{KZ}.
}.

\item[\bolla] In the complement of the two above small sets, namely in $\cT\cap (\DD\times \T^n)$, 
{\sl all trajectories lie on maximal KAM tori} for $\ham$ and the 
dynamics is quasi--periodic (with possibly exponentially small Diophantine constant).

\end{itemize}

\nl
($\rm R_2$) The sets $\Rd$ and $\cA$ depend, in particular, upon $\e$ and $\K$ and choosing $\K$ as a suitable function of $\e$ so as to obtain  larger neighborhoods of double resonances, leads to different coverings of the phase space with improved  measure estimates on $\cA$. For example:

\giu
\bolla{\sl Almost--exponential density outside a region of measure $O(\e|\log \e|^{2\g})$:}\\ {  Letting  $\K=\log^2 \e$ in Theorem~\ref{prometeo} 
with $\e$ small enough (so that \equ{torino} is met), the   sets $\Rd$ and $\cA$ 
satisfy the  estimates
\beqno
\meas(\Rd\times \torus^n)\le \cdr\, \e\,  |\log \e|^{2\g}\,,\qquad \quad\meas \cA\le \e^{ (\log|\e|)/\cteo}\,.
\eeqno
In other terms, {\sl outside a neighborhood of $O(\e\,  |\log \e|^{2\g})$ of double resonances, the `non--torus set'  is almost exponentially small} (i.e., smaller than any power of 
$\e$).}

\giu
If we allow for a neighborhood of double resonances of size $\e^a$ with $a<1$, we get a pure exponential density of KAM tori outside $\Rd\times \torus^n$:

\nl
\bolla{\sl Exponential density outside outside a region of measure $O(\e^a)$:}\\ {   Let   $0<a<1$ and choose $\K=1/\e^{\bar a}$ with ${\bar a}:=(1-a)/\g$  in Theorem~\ref{prometeo} and let $\e$ be small enough. Then,  the  sets $\Rd$ and $\cA$
satisfy the estimates
\beq{ex}
\meas(\Rd\times \torus^n)\le \cdr\, \e^a\,,\qquad \quad\meas \cA\le e^{-1 /\cteo\e^{\bar a}}\,.
\eeq
In other terms, {\sl outside a neighborhood of $O(\e^a)$ of double resonances, the `non--torus set'  is  exponentially small in
$1/\e$}.}

\nl
($\rm R_3$) As anticipated in the Introduction, 
Corollary~\ref{prometeotris} and Corollary~\ref{prometeobis} prove -- or, more precisely, are in agreement with -- 
the conjectures made by  Arnol'd, Kozlov and Neishtadt mentioned in the Introduction.
Notice, however,  that the argument sketched by Arnol'd, Kozlov and Neishtadt in \cite{AKN} to support their conjecture for $n\ge 3$ (reported in footnote~\ref{qqq} above)   only suggests  a {\sl lower} bound on the measure of the non--torus set, while, here,  we provide a rigorous {\sl upper} bound on it.
\\
We also mention that  that Corollary~\ref{prometeotris} is a particular case (with slightly better constants) of  \equ{ex}, but, since we are in two action--dimensions, it is possible to take $\Rd$ simply as a small ball around the origin (while for $n\ge 3$  it is a  more complicate set). 

\nl
($\rm R_4$) The `Kolmogorov's set' $\cT$ if formed by primary and secondary tori: Such secondary tori are {\sl not} deformation of integrable tori and, in particular, they are never graphs over $\T^n$. \\
We remark also that the  set $\cT$ is {\sl not} contained in $\DD\times \T^n$, and indeed many of the invariant tori in $\cT$  (corresponding to a set of measure $\sim \sqrt\e$) have oscillations outside $\DD\times \T^n$; this fact is unavoidable, as near the boundary tori do oscillate by a quantity of order $\sqrt\e$. 

\nl
($\rm R_5$) In Theorem~\ref{prometeo} there appear two constants $\cteo$  and $\cdr$ 
(the other constants appearing in the corollaries of Theorem~\ref{prometeo} are simply related to such two constants).
The constant $\cdr$ depends only on the action--dimension $n$
(compare Lemma~\ref{coverto} below).
More relevant  for measure estimates in phase space is the constant $\cteo$.\\ 
 The constant $\cteo$ 
can  be  calculated in terms of a few analytic properties of the potential $f$. In fact, $\cteo$ depends on 
six parameters:  $n\ge 2$; $s>0$; a positive number $\d$ quantifying  property
\equ{P1} (compare \equ{P1+}); a positive number $\b$ quantifying   property
\equ{P2} (compare \equ{P2+});  and two more positive parameters $\upxi$ and $\mathtt m$ (introduced in Definition~\ref{dracula}), which, in turn, 
are suitable non--degeneracy parameters associated to the (normalized) second derivative of the action--to--energy maps  associated to the integrable 1--degree--of--freedom Hamiltonians $p^2+  \fproj f(q)$ with $\noruno{k}\le \Nf$ where $\Nf$ is as in Definition~\ref{duda}.
For more details on how the constant $\cteo$ depends on the various parameters, see \equ{verace} and \equ{finally}  below.

\nl
($\rm R_6$) The actual effective hypothesis of Theorem~\ref{prometeo} is the condition  $\K\ge \cteo$ in \equ{torino}, which is 
used almost in  all the proofs given in this paper.\\
The second condition, $\e\K^\g\le 1$, strictly speaking, is not really necessary, and it is used (for simplicity) only in Proposition~\ref{lessing} below. 
However, as it obvious even from the statement of the Theorem~\ref{prometeo} (compare with the measure estimate in item~(ii)), if $\e\K^\g$ is not small, 
some  statements might be empty.

\nl
($\rm R_7$) One of the main issues in singular KAM theory in analytic class is the identification of a suitable  generic class of analytic  potentials.
We stress that the choice of the class $\Gns$ is tailored on the simple structure of natural Hamiltonian systems.

\subsubsection*{Open problems}

(i) It might  not be difficult to prove that Theorem~\ref{prometeo} can be generalized to natural systems with Hamiltonians
 of the form $h(y)+\e f(y,x)$ with $h$ Kolmogorov non--degenerate and $f$ verifying \equ{P1+} and \equ{P2+} {\sl uniformly in $y\in \DD$}. However, in such a case, the class of perturbations $f$ would not be generic, since,
 in general, one expects that the Fourier coefficients $f_k(y)$ may vanish at some points $y\in\DD$. Selecting an analytic generic class of perturbations, to which Theorem~\ref{prometeo} extends, is a non--trivial issue.

 \nl
(ii)
The results in this paper hold for generic potentials $f\in \Gns$, and, do not cover special  cases such as, e.g.,  the case of $f$ trigonometric polynomial, or other cases with special symmetries, as they arise, e.g., in Celestial Mechanics.

\nl
(iii)   In view of our techniques, the logarithm appearing in Corollary~\ref{prometeobis} appears to be unavoidable, and one may wonder if it is possible to get rid of it.

\nl
(iv)${}^*$
The argument 
sketched by Arnol'd, Kozlov and Neishtadt
for the lower bound on the measure 
of the non--torus set rests on the claim
that a general  real analytic  Hamiltonian
system with no small parameters has
a positive measure set free of invariant
Lagrangian tori, however this is not been proved\footnote{For  
related results in
 smooth
category, see \cite{LMS}.}.

\nl
(v)${}^*$
{\sl Generic Arnol'd diffusion in analytic class:} It is natural to expect that, for $n\ge 3$ and for generic potentials  $f\in\Gns$,  almost every non--empty energy  level  of $\ham=\frac12|y|^2+\e f(x)$  is orbit--connected, i.e.,  arbitrary
neighborhoods of two points on such levels  intersect an orbit of $\phi_\ham^t$   no matter how small   $\e$ is.

\section{Prerequisites}

In this section we recall a few prerequisites, which are needed to discuss the `secondary'  nearly--integrable structure that appears near simple resonances.

\nl
We begin by recalling the averaging theory for  nearly--integrable real analytic Hamiltonian systems as discussed in \cite{BCnonlin} and \cite{BCuni}, especially designed for neighborhoods of simple resonances. 

\nl
On one hand, apart from a finite (although arbitrarily large) number  of simple resonances of order less than $\Nf$, the secular (averaged) Hamiltonians have a uniform normal form with a potential close to a shifted cosine (\S~\ref{averroe}). On the other hand, the  secular Hamiltonians at simple resonances of order less or equal than $\Nf$ admit a simple normal form called `Generic Standard Form' (\S~\ref{gufo}). Such Generic Standard Form Hamiltonians can, then,  be put into action--angle variables that can be analytically and {\sl uniformly} controlled thanks to the theory developed in \cite{BCaa23} (\S~\ref{criseide}).

\subsection{Averaging}\label{averroe}
{\bf Non--resonant and simply resonant sets}\ \\
First, we need to introduce suitable non/simply--resonant sets, which  depend upon some quantitative parameters measuring Fourier cut--offs and small divisors. In particular, we define the `non--perturbative' set $\Rd$ -- which is a neighborhood of double resonances up to order $\K$ -- and prove item (ii) of Theorem~\ref{prometeo}.

\nl
For $k\neq 0$, denote by $\proiezione_k$ and $\proiezione_k^\perp$, 
the standard   orthogonal projections
\beq{barocco}\ts
\proiezione_k y:=(y\cdot \frac{k}{|k|}) \, \frac{k}{|k|}\,,
\qquad
\pko y:=y-\proiezione_k y\,,
\eeq
and  let $\KO$, $\K$ and $\a$ be positive numbers such that
 \beq{athlone}
\KO\geq 2\,,\qquad 
\K\geq 6\KO
\,,\qquad
\a:= \sqrt\e \K^\nu\,,\ \qquad
\nu:= \ts\frac92n+2\,.
\eeq
Recall the definition of $\gen_K$ in \equ{genk} and  define  the following real subsets of $\DD=\{y\in\R^n\st |y|<1\}$:
\beqa{neva} 
&&\Rz:=\{y\in \DD: \ |y\cdot k|>\ts \frac\a2\,,\ \forall  k\in \genKO \}\,, \ 
\\
\label{sonnosonnoBIS}
&&
\left\{\begin{array}{l}
\Ruk:=
\big\{y\in \DD:
 |y\cdot k|<\ts\a;\,  |\pko y\cdot \ell|> \frac{3 \a \K}{|k|}, \forall
\ell\in {\mathcal G^n_{\K}}\bks\Z k\big\}\,,\quad
(k\in\genKO);
\\
\Ru:=\bigcup_{k\in  \genKO} 	\Ruk\,;
\end{array}\right.\\
\label{sonnosonnosonno}
&&
\Rd:=\DD\setminus (\Rz\cup\Ru)\,.
\eeqa
The first key remark  is that {\sl  the measure of $\Rd$ is proportional to $\e$},  as shown in the following 
lemma\footnote{Compare, also, Lemma 2.5 in   \cite{BCuni}.}, which  proves item (ii) of Theorem~\ref{prometeo},

\begin{lemma}\label{coverto}\  
There exists a constant $\cdr=\cdr(n)>1$ such that 
\beq{teheran44}
\meas \Rd \le \cdr\,  \a^2\  \K^{2n} =  \cdr\, \e\ \K^\g\,,\qquad \g:=11 n +4\,.
\eeq
\end{lemma}
\proof 
First observe that from the definitions of $\Rz$, $\Ruk$ and $\Rd$ in \equ{neva}, \equ{sonnosonnoBIS} and \equ{sonnosonnosonno}, it follows immediately that 
\beq{due}
\Rd \subseteq\bigcup_{k\in \genKO} 
\bigcup_{ \sopra{\ell\in \mathcal G^n_{\K}}{\ell\notin  \Z k}}  \Rd_{k,\ell}\,,
\eeq
with
\beq{defi}
\Rd_{k\ell}:= \big\{y\in \DD:
 |y\cdot k|<\ts\a;\,  |\pko y\cdot \ell|\le \frac{3 \a \K}{|k|}\big\}\,,\qquad (k\in\genKO\,,\ \ell\in \genK\bks\Z k)\,.
\eeq
Let us, then,  estimate the measure of $\Rd_{k,\ell}$ in \equ{defi}.
Denote by  $v\in\R^n$  the projection of $y$ onto  
the plane generated by $k$ and $\ell$
(recall that, by hypothesis, $k$ and $\ell$ are not parallel);
then,   
\begin{equation}\label{soldatino}
|v\cdot k|=|y\cdot k|<\a\,, \qquad |\proiezione_k^\perp v \cdot \ell|
=|\proiezione_k^\perp y \cdot \ell|
 \le 
3\a\K /|k|\,.
\end{equation}
Set
\beq{bacca}\ts
h:=\pk \ell= \ell -\frac{\ell\cdot k}{|k|^2} k\,.
\eeq
Then, $v$ decomposes in a unique way as
$v=a k+ b h$
for suitable $a,b\in\R$.
By \eqref{soldatino},
\beq{goja}
|a|<\frac{\a}{|k|^2}\,,\qquad
|\pk v\cdot\ell|
=|bh\cdot \ell| \le 3\a\K /|k|\,,
\eeq
and, since $ |\ell|^2 |k|^2-(\ell\cdot k)^2$ is a positive integer (recall, that $k$ and $\ell$ are integer vectors not parallel), 
$$
|h\cdot \ell|
\eqby{bacca}
\frac{ |\ell|^2 |k|^2-(\ell\cdot k)^2  }{|k|^2}
\ge \frac1{|k|^2}\,.
$$
Hence, 
\beq{velazquez}
|b|\le 3 \a \K |k|  \,.
\eeq
Then,  write $y\in \Rd_{k,\ell}$ as $y=v+v^\perp$ with 
$v^\perp$ in the orthogonal
complement of the plane generated by $k$ and $\ell$. Since $|v^\perp |\le |y|< 1$  and $v$ lies in the plane spanned by $k$ and $\ell$ inside a rectangle of sizes of length $2\a/|k|^2$ and $6 \a \K |k|$ (compare \equ{goja} and \equ{velazquez}),
we find
\beqno\ts
\meas(\Rd_{k,\ell})\le \frac{2\a}{|k|^2}\, (6 \a \K |k|)\ 2^{n-2}=3\cdot 2^n \, \a^2 \frac{\K}{|k|}\,,\quad  \forall 
\left\{\begin{array}{l}
k\in\genKO\,,\\
 \ell\in  \mathcal G^n_{\K}\bks \Z k\,.
\end{array}\right.
\eeqno
Thus, since $\dst\sum_{k\in\genKO}|k|^{-1}\le c\, \KO^{n-1}$ for a suitable $c=c(n)$,  \equ{teheran44}  follows immediately  
taking $\cdr=18 c$. 
\qed

\rem  (i) By the second relation in \equ{torino} it follows  that  $\a<1/\K^n$.
\\
(ii)
$\Rz$ is a  non resonant  set up to order $\KO$; $\Ruk$  is a  simply resonant set around ${\cal R}_k$ but far away from any ${\cal R}_\ell$ with  $\ell\in \genK$ ($\ell\neq k$); $\Rd$  is contained in a neighborhood of  double resonances of 
order $\K$ (compare relation \equ{due} below). According to the terminology in \cite{P93}, $\Rz$ is $(\a/2,\KO)$ completely non--resonant,  while, for each $k\in\genKO$, the set $\Ruk$ is $(2 \a \K/|k|)$--non resonant modulo $\integer k$ up to order~$\K$.
\\
(iii)
From the definition of $\Rd$ in \equ{sonnosonnosonno} it follows trivially that $\{{\cal R}^i\}$ is a covering of $\DD$.
\\
(iv) Having  two different Fourier cut--offs $\KO$ and $\K$ is necessary in order to  obtain high order `cosine--like' normal forms as described  in point (iii) of the averaging Theorem~\ref{normalform} below; compare also \cite{BCnonlin}.

\erem
{\bf Notations} \ \\
Given $m\ge 1$,  $D\subseteq \real^m$, and $r>0$, let us denote $D_r$  the complex neighborhood of $D$ given by
$$D_r := \bigcup_{z\in D} \{y\in\complex^m\st |y-z|<r\}\,.$$
For $s>0$,  let $\torus^m_s$ denote  the complex neighborhood of width $2s$ of $\torus^m$ given by
\beqno
\torus^n_s:=\{x=(x_1,...,x_m)\in\complex^n:\ \ |\Im x_j|<s\}/(2\pi \integer^m)
\,.
\eeqno
We shall also use the notation $\Re(V_r)$ to denote the {\sl real} $r$--neighbourhood of $V\subseteq \real^n$, namely,
\beq{dire}
\Re(V _r) := V_r\cap \real^n= \bigcup_{z\in V}\ \{y\in\R^n\st |y-z|<r\}\,.
\eeq
Given  $D\subseteq \real^n$ and a function $f$ defined, respectively, on $D_r$, $\T^m_s$, $D_r\times \torus^m_s$, we denote its sup norm, respectively,  by
\beqno
|f|_{D,r}:= \sup_{y\in D_r}|f(y)|\,,\quad
|f|_{s}:= \sup_{x\in\torus^m_s}|f(x)|\,,\quad
|f|_{D,r,s}:= \sup_{(y,x)\in D_r\times \torus^m_s}|f(x,y)|\,.
\eeqno
{\bf Analytcity parameters} \ 
To formulate properly next (normal form) theorem, we  need to introduce a few `analyticity parameters'
related to the analyticity width $s$ and the numbers $\a$, $\KO$ and $\K$ in \equ{athlone}. We define:
\beqa{dublino}
&&\ts r_{\rm o}:=\frac{\a}{16 \KO}\,, \quad r_{\rm o}':= \frac{r_{\rm o}}2\,,\quad\
\ts s_{\rm o}:=s\big(1-\frac1{\KO}\big)\,, \quad s_{\rm o}':=s_{\rm o}\big(1-\frac1{\KO}\big)\,,
\nonumber
\\
&&\ts  s_{\varstar}:=s\big(1-\frac1{\K}\big)\,,\quad \ts   s_{\varstar}':=s_{\varstar}\big(1-\frac1{\K}\big)\,,\quad
s'_k:=|k|_{{}_1}s_{\varstar}' \,,\\
\nonumber
&&\ts r_k:=\frac{\a}{ |k|}\,,\quad 
r_k':=\frac{r_k}2
\,,\nonumber\\
&&\ts \tilde r_k:=\frac{{r_k}}{ \itcu |k|}\,, \qquad\tilde s_k:=
\frac{s}{\itcu |k|^{n-1}}\,,
\ \  {\rm where}\quad \ts \itcu:=5n(n-1)^{\frac{n-1}2}\,.\nonumber
\eeqa
We also need the following  consequence of Bezout's Lemma, which will allow to define the effective `resonant angle' near  simple resonances:

\begin{lemma}\label{murino}
For any  $k\in \gen$  there exists   a matrix
$\hAA\in\Z^{(n-1)\times n}$
such that\footnote{$|M|_{{}_\infty}$, with $M$ matrix (or vector), denotes the maximum norm $\max_{ij}|M_{ij}|$ (or $\max_i |M_i|$).
}
\beqa{atlantide}
&& \AA:=\binom{k}{\hAA}
=\binom{k_1\cdots k_n}{\hAA}\in\ {\rm SL}(n,\Z)\,,
\nonumber
\\
&&
|\hAA|_{{}_\infty}\leq |k|_{{}_\infty}\,,\ \ 
|\AA|_{{}_\infty}=|k|_{{}_\infty}\,,\ \ 
|\AA^{-1}|_{{}_\infty}\leq 
(n-1)^{\frac{n-1}2} |k|_{{}_\infty}^{n-1}\,.
\eeqa
\end{lemma}
\proof
From B\'ezout's lemma it follows easily that\footnote{See Lemma~A.1 in \cite{BCnonlin} for a proof.}: 

\nl
{\sl 
Given $k\in\integer^n$, $k\neq 0$ there exists a matrix $\AA=(\AA_{ij})_{1\le i,j\le n}$ with integer entries such that $A_{nj}=k_j$ $\forall$ $1\le j\le n$, $\det\AA={\rm gcd}(k_1,...,k_1)$, and 
$|\AA|_{{}_\infty}=|k|_{{}_\infty}$.}

\nl
Since $k\in \gen$,  it is ${\rm gcd}(k_1,...,k_1)=1$ and, therefore,  $\det\AA=1$.\\
 The  first two relations in \equ{atlantide} are  consequence of the above statement.
\\
Observing that
for  any $m\times m$ matrix ${\rm M}$, one   has 
$|\det {\rm M}|\leq m^{m/2} |{\rm M}|_{{}_\infty}^m$,  the bound on $|\AA^{-1}|_{{}_\infty}$ follows from D'Alembert expansion of determinants. \qed

\nl
The following normal form result -- proven in\footnote{See Theorem 2.1 and  the Covering Lemma 2.3 in \cite{BCuni}.
} \cite{BCuni} --  holds:

\begin{theorem}[Normal Form Theorem]\label{normalform}
Fix  $n\ge 2$, $s>0$. 
Let $\ham$ be as in \equ{ham} 
with  $f\in\Bns$ satisfying   \eqref{P1+} with $\Nf$ as in \equ{enne}; let 
  \equ{athlone}$\div$\equ{sonnosonnoBIS} and \equ{dublino} hold.
For $k\in\genKO$, let $\AA$  be  the matrix 
in Lemma~\ref{murino} and define the following real sets: 
\beq{codino}
\Rzt:= \Re (\Rz_{r'_{\rm o}/2})\,,\quad \Rukt:= \Re(\Ruk_{r'_k/2})\,, 
\quad
\DDD^k:= \AA^{-T}\Rukt\,,\qquad   (k\in\genKO)\,.
\eeq
Then, there exists a constant $\bfco=\bfco(n,s,\d)\ge \Nf$ such that if
$\KO\ge \bfco$,   there exist
 real analytic symplectic maps
\beq{trota2}
\Psi_{\rm o}: \Rz_{r_{\rm o}'}\times \T^n_{s_{\rm o}'} \to 
\Rz_{r_{\rm o}}\times \T^n_{s_{\rm o}} 
\,,
\qquad
\Psi^k:
\DDD^k_{\tilde r_k}\times \T^n_{\tilde s_k} 
\to 
\Ruk_{r_k} \times \T^n_{s_\varstar}
\eeq
having the following properties.

\nl
{\rm (i)} 
In the  symplectic variables $(y,x) \in \Rz_{r_{\rm o}}\times \T^n_{s_{\rm o}}$,   $\ham$  takes the form:
$$
\hamo(y,x) := \big(\ham\circ\Psi_{\rm o}\big)(y,x)
=\frac{|y|^2}2+\e\big( g^{\rm o}(y) +
 f^{\rm o}(y,x)\big)\,,\qquad \langle f^{\rm o}\rangle=0\,,$$
with  $g^{\rm o}$ and $f^{\rm o}$ satisfying
\beq{552}
| g^{\rm o}|_{r_{\rm o}'}
\leq
\tettao:=
\frac{1}{\K^{6n+1}}\,, \qquad  
\modulo f^{\rm o} \modulo_{r_{\rm o}',s_{\rm o}'} 
\leq  e^{-\KO s/3}\,.
\eeq
{\rm (ii)}  Let $k\in\genKO$. In the  symplectic variables $(\tty,\ttx)=\big(\tty,(\ttx_1,\hat\ttx)\big)\in \DDD^k_{\tilde r_k}\times \T^n_{\tilde s_k}$,   $\ham$  takes the form:
\beq{dopomedia}
\cH_k(\tty,\ttx):=\ham\circ \Psi^k(\tty,\ttx)=\hamsec_k(\tty,\ttx_1)+ \e
\bar f^k(\tty,\ttx) \,,\quad (\tty,\ttx) \in 
\DDD^k_{\tilde r_k}\times \T^n_{\tilde s_k}\,,
 \eeq
 where
 \beq{hamseccippa}
\hamsec_k(\tty,\ttx_1):= \frac12 |\AA^T\tty|^2+\e \mathtt g^k_{\rm o}(\tty)+
\e \mathtt g^k(\tty,\ttx_1)
\eeq
is real analytic   in $\tty\in  \DDD^k_{\tilde r_k}$ and $\ttx_1\in \torus_{s'_k}$.
In particular
 $\mathtt g^k(y,\cdot)\in\hol_{s'_k}^1$
for every $y\in \DDD^k_{\tilde r_k}$. Furthermore, the following estimates hold:
\begin{equation}\label{cristinacippa}
|\mathtt g^k_{\rm o}|_{\tilde r_k}
\leq \tettao\,,\qquad
\modulo  \mathtt g^k-\fproj f\modulo_{{\tilde r_k},s'_k} 
\leq  \tettao\,,\qquad 
\modulo \bar f^{k} \modulo_{{\tilde r_k},{\tilde s_k}} \le
 e^{- \K s/3}\,.
\end{equation}
{\rm (iii)} If  $k\in\genKO$ satisfies $\noruno{k}\geq \Nf$, then there exists  $\sa_k\in[0,2\pi)$ such that 
\beqno
\cH_k
=
 \frac12 |\AA^T\tty|^2+\e \mathtt g^k_{\rm o}(\tty)+
2|f_k|\e\ 
\big[\cos(\ttx_1 +\sa_k)+
F^k_{\! \varstar}(\ttx_1)+
\mathtt g^k_{\! \varstar}(\tty,\ttx_1)+
\mathtt f^k_{\! \varstar} (\tty,\ttx)
 \big]\,,
\eeqno
where   
$$F^k_{\! \varstar}(\sa):=\frac{1}{2|f_k|}\sum_{|j|\geq 2}f_{jk}e^{\ii j \sa}\in\hol_1^1\,,\qquad \modulo F^k_{\! \varstar} \modulo_1\leq 2^{-40}\,.$$ 
Moreover,
 $\mathtt g^k_{\! \varstar}(y,\cdot )\in\hol_1^1$
(for every $y\in \DDD^k_{\tilde r_k}$), $\fproj\mathtt f^k_{\! \varstar}=0$, and  one has 
\beqno\ts
\modulo \mathtt g^k_{\! \varstar}\modulo_{\tilde r_k,1}\le 
\frac{1}{\K^{5n}}\,,\quad\qquad
\modulo \mathtt f^k_{\! \varstar} \modulo_{\tilde r_k,\tilde s_k} 
\leq
e^{-\K s/7}\,.
\eeqno
{\rm (iv)} Finally, the following `coverings' holds:
\beqa{surge}
&&\Psi_{\! \rm o}\big(\Rzt\times \T^n\big)\supseteq
 \Rz\times \T^n\,,\qquad
\Psi^k
\big(\DDD^k\times \T^n\big)\supseteq
\Ruk \times \T^n \,.
\eeqa
\end{theorem}

\rem\label{munro}
(i) Beware that, while $\Psi_{\rm o}$ is a map close to the identity, $\Psi^k$ is not, as it is  the composition of a linear transformation\footnote{Namely,  the symplectic transformation, the generating function of which is given by  $\tty\cdot \AA x$, and
which maps  the resonant combination $k\cdot x$ to the `resonant'  angle $\ttx_1$.}
with a near--to--identity map.

\nl
(ii) The larger covering in \equ{codino} is introduced so that \equ{surge} holds: Such a property will be essential in covering also boundary regions by KAM tori without leaving out (as it happens in standard KAM theory) regions of size of order $\sqrt\e$, a fact that, for our purposes, would be clearly not acceptable.

\nl
(iii)
Point (iii) in Theorem~\ref{normalform}  shows that the secular Hamiltonian $\hamsec_k$ (obtained disregarding the exponentially small perturbation 
$\mathtt f^k_{\! \varstar}$) has a potential, which is $O(1/\K^{5n})$--perturbation of the  `cosine--like function' 
$$\cos(\ttx_1 +\sa_k)+F^k_{\! \varstar}(\ttx_1)\,,\qquad {\rm where}\quad   \modulo F^k_{\! \varstar} \modulo_1\leq 2^{-40}\,.$$
This means that for $\noruno{k}\ge \Nf$, {\sl the secular Hamiltonians at simple resonances all look the same}, allowing, in particular, for a uniform analysis in terms  of action--angle variables (compare \S~\ref{criseide} below).

\nl
Notice also that the perturbation $\mathtt f^k_{\! \varstar}$, which is bounded by $e^{-\K s/7}$ has a factor $|f_k|$ in front of it and that such a factor, in turn, may be 
exponentially small (since $|f_k|\sim e^{-\noruno{k}s}$ for large $\noruno{k}$).

\nl
(iv) For later use we  observe 
that\footnote{If $s\ge 1$ then $\Nf\ge 2\ge 2/s$, while if $s<1$ then the logarithm in
\equ{enne} is larger than one, so that $\Nf\ge 2/s$ also in this case.}
\beq{bollettino}
\KO\ge \Nf\ge2 \ttcs\,, \qquad {\rm where}\qquad \, \ttcs:=\max\{1,1/s\}\,.
\eeq

\erem

\subsection{Generic Standard Form  at simple resonances}\label{gufo}

It turns out that the secular Hamiltonians $\hamsec_k$ in \equ{dopomedia}--\equ{hamseccippa} in the Normal Form Theorem~\ref{normalform} for $k\in\genKO$, have a common uniform analytic structure: They 
can be put into a standard form, which has  uniform (in $k\in\genKO$) analytic characteristics. The precise formulation of this fact is the main theorem in \cite{BCuni}, whose statement needs some preparation. 

\begin{definition} \label{morso}
{\bf (1D Hamiltonians in standard form)} 
Let $\hat D \subseteq \R^{n-1}$ be a bounded domain,  $\Ro>0$ and $D:=  (-\Ro ,\Ro ) \times\hat  D
$. We say that a real analytic Hamiltonian $\Hpend$ is in Generic Standard Form (in short, `standard form') with  respect to the symplectic variables $(p_1,q_1)\in (-\Ro ,\Ro )\times\torus$ and  `external actions' 
$\hat p=(p_2,...,p_n)\in \hat D$, if $\Hpend$ has the form 
 \beq{pasqua}
  \Hpend(p,q_1)=
\big(1+ \cin(p,q_1)\big) p_1^2  
  +\Gm(\hat p, q_1)
  \,,
\eeq
where $p=(p_1,\hat p)=(p_1,p_2,...,p_n)$, and the following specifications hold.

\begin{itemize}

\item[\bolla]  $\cin$ and $ \Gm$ are real analytic functions defined on, respectively, $D_\ro\times\T_\so$ and $\hat D_\ro\times \T_\so$ for some $0<\ro\leq\Ro$ and $\so>0$;

\item[\bolla]  
$\Gm$ has zero--average and there exists a  zero--average
 function $\GO$ (the `reference potential')  depending only on $q_1$ such that, for some $\morse>0$, 
 $\GO$ is $\morse$--Morse\footnote{Recall Definition~\ref{buda}.};
 
 \item[\bolla] 
 the following estimates hold:
 \beq{cimabue}
 \left\{\begin{array}{l} \dst \sup_{\torus^1_\so}|\GO|\le \suca\,,\\
 \dst \sup_{\hat D_\ro\times \torus^1_\so}|\Gm-\GO| \leq
\suca
 \lalla
\,,\quad{\rm for\ some}\quad 0<\suca\le \ro^2/2^{16} 
\,,\ \ 0\le \lalla<1\,,
\\
\dst \sup_{D_\ro\times \torus^1_\so}|\cin| \leq
\lalla\,.
\end{array}\right.
\eeq
\end{itemize}
\end{definition}
We shall call $(\hat  D,\Ro,\ro,\so,\morse,\suca,\lalla)$ {\sl the `analyticity characteristics' of $\Hpend$
with respect to the reference  potential $\GO$}.

\rem
(i) A Hamiltonian in standard form $\Hpend$ has the analytic features of its reference natural Hamiltonian 
$$\bHpend:=p_1^2+\GO(q_1)\,.
$$ 
In particular, for $\lalla$ small with respect to $1/\upkappa$, $\Hpend$ has the same finite (because of analyticity) number of equilibria (which lie on the $q_1$ axis) of $\GO$ and in the same relative order,  which is also preserved  by the corresponding critical energies; compare  Lemma~\ref{enricone} below. 

\nl
(ii)
 If $\Hpend$ is in standard form, then $\morse$ and $\suca$ satisfy the relation\footnote{By \equ{cimabue}, $\morse \le |\bar \Gm(\sa_i)- \bar \Gm(\sa_i)|\le 2 \max_\torus|\bar \Gm|\le 2\suca$.}
$\suca/\morse\ge 1/2$.
Furthermore,  one can always
fix   a number $\upkappa\geq 4$  so that:
  \begin{equation}\label{alce}
1/\upkappa\leq \so\leq 1\,,\qquad
1\leq
\Ro/\ro\leq \upkappa\,,\qquad
1/2\leq
\suca/\morse
\leq \upkappa \,.
\end{equation}
Such a  parameter $\upkappa$ rules the main scaling properties of these Hamiltonians.

\nl
(iii) Hamiltonians in standard form are particularly suited for the analytic theory of action--angle variables (in neighborhoods of separatrices) as developed  in \cite{BCaa23}, where the notion of Generic Standard Form has been introduced. Such action--angle variables will be reviewed in \S~\ref{criseide} below.

\nl
(iv) The smallness of the `adimensional ratio' $\suca/\ro^2$ in \equ{cimabue} is needed in the analytic theory of action-angle variables for Hamiltonians in standard form developed in \cite{ BCaa23}, however the factor $1/2^{16}$ is rather arbitrary and  not optimal.
\erem
{\bf Notation} 
If $w$ is a vector with $n$ or $2n$ components, $\hat w=(w)^{\widehat{}}$ denotes the last $(n-1)$ components; 
if $w$ is vector with $2n$ components, $\check w=(w)^{{\!\!\widecheck{\phantom{a}}}}$ denotes the first $n+1$ components. 
Explicitly:
\beq{checheche}
w=(y,x)=\big((y_1,...,y_n),(x_1,...x_n)\big)\quad \Longrightarrow\quad
\left\{
\begin{array}{l}
\hat w=(w)^{\widehat{}}=(x_2,...,x_n)=\hat x\,,\\
\hat y=\,(y)^{\widehat{}}=(y_2,...,y_n)\,,\\
\check w=(w)^{{\!\!\widecheck{\phantom{a}}}}=(y,x_1)\,,\\
w=(\check w,\hat w)\,.
\end{array}
\right.\eeq

\nl
Next, we introduce a special simple group of symplectic transformations, which will appear in Theorem~\ref{sivori} below.

\dfn{dadaumpa0}
Given a domain $\hat {\rm D}\subseteq \R^{n-1}$,
we denote by
 $\Gdag$   the 
abelian group of 
symplectic diffeomorphisms $\Psi_{\! \ta}$ 
of $(\R\times\hat {\rm D})\times \R^n$
given by
\beq{psiab}
(p,q)\in(\R\times\hat {\rm D})\times \R^n\stackrel{\Psi_{\! \ta}}\mapsto (P,Q)=
(p_1+\ta(\hat p),\hat p,q_1,\hat q-q_1\partial_{\hat p}
\ta(\hat p))\in(\R\times\hat {\rm D})\times \R^n\,,
\eeq
with $\ta:\hat {\rm D}\to \R$ smooth.
 \edfn

\rem\label{quicascalasino}
The group properties of $\Gdag$ are trivial: 
\beq{cotugno}
{\rm id}_{\Gdag}=\Psi_{\! 0}\,,\qquad   
\Psi_{\! \ta}^{-1}=\Psi_{\! -\ta }\,,\qquad
\Psi_{\! \ta}\circ\Psi_{\! \ta'}=\Psi_{\! \ta+\ta'}\,.
\eeq
Notice that, unless $\partial_{\hat p} \ta\in\integer^{n-1}$,
maps $\Psi_{\! \ta}\in \Gdag$  {\sl do not induce  well defined maps}\footnote{In general, given $A\in {\rm SL}(n,\Z)$ and  a $2\pi$--multi--periodic function $f:\real^n\to\real^n$,  we identify the $\real^n$--map $x\in\real^n\to f(x)=A x+ g(x)\in\real^n$ with the $\torus^n$--map given by 
$\theta\in\torus^n\to F(\theta)=\pi_{{}_{\torus^n}}\big(Ax+f(x)\big)\in\torus^n$ where $\theta=x+2\pi\integer^n$ and $x\to \pi_{{}_{\torus^n}}(x)=x+2\pi\integer^n$ is the projection of $\real^n$ onto $\torus^n$.}   
 $$q\in \torus^n\mapsto (q_1,\hat q-q_1\partial_{\hat p}
\ta(\hat p))\in\torus^n\,,$$ 
a fact that will create a problem in applying the theory of this and next section to the normalized Hamiltonians $\cH_k$ of Theorem~\ref{normalform}; compare Remark~\ref{portogallo}--(ii) below.
\erem

\nl
Let us now spell out all the {\sl assumptions  and definitions, which,  from now, will be part of the 
hypotheses
 of all statements regarding the natural system with  Hamiltonian $\ham$ as in \equ{ham}}.

\begin{assumption}\label{assunta}
Fix  $n\ge 2$, $s>0$, and let $\ham$ be as in \equ{ham} 
with  $f\in\Gns$ (Definition~\ref{sicuro}) satisfying   \eqref{P1+} and \equ{P2+} for some $0<\d\le 1$ and $\b>0$  with\footnote{By Lemma~\ref{telaviv} such $\d$ and $\b$ always exist.} $\Nf$ as in \equ{enne}.
\end{assumption}

\begin{definition}\label{assunto} Given $\ham$ as in Assumption~\ref{assunta}, we define the following sets and parameters.
\\
\bolla  Let $\KO$, $\K$ and $\a$ be as in   \equ{athlone};
let ${\mathcal R}^i$'s be the domains defined in \equ{neva}$\div$\equ{sonnosonnoBIS}; let the definitions in \equ{dublino} hold (`analiticity parameters').
\\
\bolla For $k\in\genKO$, let $\AA$  be  the matrix 
in Lemma~\ref{murino}. Let $\Rzt$, $\Rukt$ and $\DDD^k$ be the real domains defined in~\equ{codino}.
\\
\bolla Define the following parameters\footnote{Here and in what follows we shall  not always indicate explicitly the dependence upon $k$.
Recall the definitions of $\itcu$, $\hAA$ and $\ttcs$ in, respectively, \equ{dublino},   Lemma~\ref{murino} and \equ{bollettino}.}:
\beq{cerbiatta}
\begin{array}l
\Ro={\a}/{|k|^2}={\sqrt\e \K^\nu}/{|k|^2}\,,\quad \itcd=4\,  n^{\frac32} \itcu\,,
\quad \ro= {\Ro}/{\itcd}\,,
\quad\e_k=\frac{2\e}{|k|^2}\,,
\phantom{\dst\int}
\\
\hat D =
\big\{ \hat\act\in\R^{n-1}: \ 
|\proiezione_k^\perp \hAA^T \hat\act|<1\,,
\ \dst
\min_{\sopra{\ell\in \genK}{\ell \notin \Z k}}
\big| \big(\proiezione_k^\perp \hAA^T \hat\act\big)\cdot \ell\big|
\geq  {\textstyle \frac{3\a\K}{|k|}}   \big\}\,,\quad D=(-\Ro,\Ro)\times \hat D\,,
\\
 \morse=\casitwo{\e_k \b,}{\noruno{k}<\Nf}{\e_k |f_k|,}{\noruno{k}\ge \Nf}\,,\,\qquad
 \chk=\casitwo
{1 \,,}{\noruno{k}<\Nf}
{ |f_k|\,,}{\noruno{k}\ge \Nf}\,,
\quad \suca=\ts \ttcs \e_k\,  \chk\,,     
\\
\so=\casitwo{\min\{\frac{s}2,1\}\,,}{\noruno{k}<\Nf}{1\,,}{\noruno{k}\ge \Nf}\,,\quad
\ts
\chs=\casitwo{s'_k\,,}{\noruno{k}<\Nf\,,}{1\,,}{\noruno{k}\ge \Nf}\,, 
\quad \dst\ 
\lalla=\frac{1}{\K^{5n}}\,. 
\end{array}
\eeq
\end{definition}
\rem
(i)  Since $|f_k|\le 1$  one has: 
\beq{roger} 
|\chk|\le1\,.
\eeq
Furthermore, by the definitions in \equ{cerbiatta} and \equ{athlone}, by \equ{roger} and \equ{bollettino}, one has
\beq{orso}
\sqrt{\suca}<\ttcs\Ro/\K^{\nu-1}<\Ro/\K^{\frac92 n}\,.
\eeq

\nl
(ii) Since $(1-\frac1\K)^{-2}<2$, by  definition of  $s'_k$ in \equ{dublino}, one has
\beq{fangorn}
\so\le 2 \chs\,.
\eeq
\erem

\nl
We can, now, state the main result of\footnote{Compare Theorem 3.1  in  \cite{BCuni}.} 
 \cite{BCuni}: 
 
\begin{theorem}[Generic Standard Form at simple resonances]\label{sivori} \ \\
Let  Assumptions~\ref{assunta} and Definitions~\ref{assunto} hold, let $\bfco$ be the constant  defined in Theorem~\ref{normalform}, and 
assume that $\KO\ge \max\{\itcd,\bfco\}$. Then,  for all $k\in\genKO$,  the following holds.

\nl
{\rm (i)}
There exists a real analytic
 symplectic transformation 
\beq{diamond}
\Phi_\diamond:(\ttp,\ttq)\ \in D\times \R^n \to
(\tty,\ttx)=\Phi_\diamond(\ttp,\ttq)\in \R^{2n}\,,
\eeq
such that: $\Phi_\diamond$ fixes $\hat\ttp$ and\footnote{I.e., in \equ{diamond} it is $\tty=\hat\ttp,\ttx_1=\ttq_1$.} $\ttq_1$; 
for every $\hat\ttp\in\hat D$ the map  $(\ttp_1,\ttq_1)\mapsto (\tty_1,\ttx_1)$ is symplectic; 
the $(n+1)$--dimensional map\footnote{Recall the notation in \equ{checheche}.}  $\check\Phi_\diamond$ depends only on the first $n+1$ coordinates $(\ttp,\ttq_1)$, is $2\pi$--periodic in $\ttq_1$ 
and, if $\DDD^k= \AA^{-T}\Ruk$ and $\hamsec_k$ are as in Theorem~\ref{normalform},  one has\footnote{
$r_k$, $\tilde r_k$ and $s'_k$ are defined in \equ{dublino}.} 
\beqa{tikitaka}
 && \check\Phi_\diamond: 
 D_{\chr}\times \torus_\chs\ \, \to
\DDD^k_{\tilde r_k}\times \T_\chs\,,\nonumber
\\
&&
\hamsec_k\circ \check\Phi_\diamond(\ttp,\ttq)=:
\textstyle{\frac{|k|^2}{2}}( \Hsharp(\ttp,\ttq_1)+
\hzk(\hat\ttp))\,, 
\\
&&\sup_{\hat\ttp\in \hat D_{2\ro}}
\big|\ts\hzk(\hat\ttp)-
\htk(\hat\ttp)
\big|
\leq
\ts\frac{12}{|k|^2} \e\lalla\,,\qquad \quad\  \htk(\hat\ttp):={\ts \frac{1}{|k|^2}} 
\normadue \proiezione^\perp_k \hAA^T \hat\ttp\normadue^2\,.
\nonumber
\eeqa
{\rm (ii)} 
$\Hsharp$ in \equ{tikitaka} is  in Generic Universal Form according to Definition~\ref{morso}:
\beq{tikitaka2}
\Hsharp(\ttp,\ttq_1)=\big(1+\cins(\ttp,\ttq_1)\big)\,  \ttp_1^2 +  \Gf(\hat\ttp,\ttq_1)\,,
\eeq
having reference potential 
 \beq{paranoia}
\GO=\bGf:= \e_k\, \fproj f\,,
\eeq  
analyticity characteristics  given in \equ{cerbiatta}, and 
$\upkappa$
verifying \equ{alce}  with\footnote{$\ttcs$ is defined in \equ{bollettino}.}
\beq{kappa}\ts
\upkappa=\upkappa(n,s,\b):=\max\big\{\itcd\,, 4\ttcs\,, \ttcs/\b
\big\}\,.
\eeq
{\rm (iii)} The map $\Phi_\diamond$ is obtained as composition of three symplectic maps: 
\beq{toto}\Phi_\diamond=\Fiuno\circ\Fidue\circ\Fitre\,,
\eeq
 where\footnote{Recall Definition~\ref{dadaumpa}.}: 
\begin{itemize}
\item[\bolla] $\Fiuno:=\Psi_{\!\giuno}\in\Gdag$ with 
$\giuno(\hat \ttp):=-\ts \frac1{|k|^2}{(\hAA k)\cdot \hat \ttp}$; 
 \item[\bolla]$\Fidue(\ttp,\ttq)=(\ttp_1+\upeta_{{}_2},\hat \ttp, \ttq_1,\hat \ttq+\upchi_{{}_2})$ for suitable real analytic functions $\upeta_{{}_2}=\upeta_{{}_2}(\hat \ttp,\ttq_1)$ and 
$\upchi_{{}_2}=\upchi_{{}_2}(\hat \ttp,\ttq_1)$ 
 satisfying
\beq{tess2}\ts
|\upeta_{{}_2}|_{4 \chr,\chs}< \frac{\e_k\chk}{\chr}\lalla\,,\qquad |\upchi_{{}_2}|_{2 \chr,\chs}< \frac{4\e_k\chk}{\chr^2}\,\lalla
\,;
\eeq
\item[\bolla]
$\Fitre:=\Psi_{\!\gitre}\in\Gdag$  for a suitable real analytic function $\gitre(\hat\ttp)$ satisfying
\beq{tess}\ts
|\gitre|_{4 \chr}< \frac{\e_k\chk}{\chr}\lalla\,.
\eeq 
\end{itemize}
\end{theorem}
\rem\label{portogallo}
{\rm (i)} The main point of the above theorem is item (ii), which shows that the `simply--resonant Hamiltonians' $\Hsharp$ in \equ{tikitaka} are in {\sl uniform} Generic Standard Form. The word  `uniform' refers to the fact that the 
 parameter $\upkappa$ (defined in \equ{kappa} and satisfying \equ{alce}) --  which rules the scaling properties of the normalized Hamiltonians $\Hsharp$  -- {\sl  does not depend upon $k$}, allowing, e.g.,  for a uniform (in $k\in \genKO$) treatment of action--angle variables (compare next Section \ref{criseide}).

\nl
(ii) There is, however, a drawback in the construction of  the above normal forms, namely,  that the maps $\Fiuno$ and $\Fitre$ appearing in the definition of $\Phi_\diamond$ (item (iii) in the above theorem), {\sl do  not induce  well defined maps on $\torus^n$}; compare Remark~\ref{quicascalasino}. Therefore, 
a non trivial homotopy issue  will have to be faced in considering the global secondary nearly--integrable structure  of the system near simple resonances. On the other hand, the map $\Fidue$ is well defined also on $\T^n$.
This matter will be discussed in details in  Section~\ref{GS}. 
\erem
The following remark explains the individual purpose of the three symplectic transformations $\Phi_{\!{}_j}$ whose composition  forms $\Phi_\diamond$. 

\rem
(i)
The  map $\Fiuno$ in the definition of $\Phi_\diamond$   is a linear map that has the purpose of block--diagonalize the quadratic part  $|\AA^T \tty|^2$ appearing in \equ{hamseccippa}, so as {\sl to obtain a kinetic part which is  the sum of a quadratic part in $p_1$ and a quadratic 
 $(n-1)$--dimensional part in   $\hat p$}. Indeed, rewriting $\Fiuno$ as
\beq{Fiu}
(\tty,\ttx)= \Fiuno(p,q):=\big({\rm U} p, {\rm U}^{-T} q\big)
\,,\qquad{\rm where}\qquad
{\rm U}:= 
\left(\begin{matrix}
1  &  - \frac1{|k|^2}(\hAA k)^T \cr 0 & \quad{\id}_{{}_{n-1}} \cr
\end{matrix}\right)\,,
\eeq
and
observing  that 
\beqno
\AA^T \tty=\AA^T{\rm U} p=p_1 k + \proiezione_k^\perp \hAA^T \hat p\,,
\eeqno
one sees that  
\beq{venere}
|\AA^T{\rm U} p|^2=|k|^2 p_1^2 + |\proiezione_k^\perp \hAA^T \hat p|^2 = 
|k|^2(p_1^2+  \htk(\hat p))
\,,
\eeq
$\htk$ being the positive definite quadratic form in $\hat \ttp=\hat p$ defined in \equ{tikitaka}.

\nl
Furthermore,
$y= (\AA^T{\rm U}) p$ if and only if $y\cdot k=p_1 |k|^2$ and $\proiezione_k^\perp y=\proiezione_k^\perp
\AA^T \hat p$, which, recalling the definition of $\ts\Rukt$ in \equ{sonnosonnoBIS}, shows that
\beq{sonnosonno2}
\AA^T{\rm U} D= \ts\Rukt \qquad \implies \quad \meas D= \meas \Rukt\,.
\eeq
Notice also that, from \equ{cerbiatta}, the definitions of $\Fiuno$ and ${\rm U}$, and the definition of 
$\DDD^k$ in Theorem~\ref{normalform}--(i), 
 it follows that 
\beq{chirone}
\cFiuno(D\times\T)={\rm U}D \times\T=\DDD^k\times\T\,.
\eeq
 Incidentally, observe that from \equ{atlantide} it follows that
 the norms of ${\rm U}$  and its inverse  satisfy the bounds\footnote{As usual, for a matrix $M$ we denote by
$\dst |M|=\sup_{u\neq 0} |Mu|/|u|$  the standard  operator norm.} 
\beq{UU}
|{\rm U}|, |{\rm U}^{-1}|\le n\sqrt{n}\,.
\eeq  

\nl
(ii)
The second map $\Fidue$ is a near--to--identity symplectic (globally well--defined) transformation, which is introduced so as to 
transform $\Hsharp$  into a Hamiltonian with  a potential  {\sl  independent of $p_1$}.

\nl
(iii)
 $\Fitre$ is  a near--to--identity symplectic map, which sets {\sl all   critical points on the line $p_1=0$}.
\erem

\subsection{Action--angle variables for 1D standard Hamiltonians}\label{criseide}

In this subsection  we review the general theory of action--angle variables for Hamiltonian systems in standard form as developed in \cite{BCaa23}, where complete proofs may be found.  

\nl
This subsection is independent from the previous ones; in particular the analytic characteristics $\hat  D$, $\Ro$, $\ro$, etc., are  
arbitrary (and do not refer to the definitions given in \equ{cerbiatta} in the specific case of the secular Hamiltonians $\hamsec_k$).

\subsubsection*{Topology of the phase space of 1D Hamiltonians in standard form}
We begin by describing  the topological structure of the $\hat p$--dependent phase space of a givern Hamiltonian $(p_1,q_1)\mapsto \Hpend(p_1,\hat p,q_1)$ in generic standard form according to Definition~\ref{morso}.

\nl
For a fixed $\hat p \in \hat D$, we take as phase space of $\Hpend$
the  subset of $\real\times \T$   given by
\beq{emmone}
\cM=\cM(\hat p):=
\{
(p_1,q_1)\in \real\times\T\, \big|  \ \ 
\Hpend(p_1,\hat p,q_1)<
\Eb\}\,,\quad \Eb:=
\Ro^2+\Ro\ro
 \,,
\eeq
where $\Ro$ and $\ro$ are as in Definition~\ref{morso}.
Although such sets depend on the parameter  $\hat p\in \hat  D$,   for $\mu$ small enough, they are close to a box:
\begin{lemma}\label{container}
Let $\Hpend$ be as in Definition~\ref{morso} and  $\cM$ be as in \equ{emmone}, and 
assume that\footnote{Recall the definition of $\upkappa$ in \eqref{alce}.}
\beq{minore}
\lalla\leq 1/(4\upkappa)^2\,.
\eeq
Then, for all $\hat p\in \hat D$, one has
\beq{gonatak}\ts
\big(-\Ro -\frac\ro3,\Ro+\frac\ro3\big)\times\T
\subseteq\cM(\hat p)\subseteq     
\big(-\Ro -\frac\ro2,\Ro+\frac\ro2 \big)\times\T\,.
\eeq
\end{lemma}
The simole proof is given in \Appendix.

\nl
Since the reference potential $\GO$ is a $\morse$--Morse function, it  has $2N$  critical points, for some $N\in\natural$,  with different critical values. Let 
$\bar \sa_0\in [0,2\pi)$ be  the unique point of absolute maximum of the reference potential $\GO$ of $\Hpend$. 
Then, the relative strict non--degenerate maximum and minimum
 points of $\GO$, $\bar\sa_i\in[\bar \sa_0,\bar \sa_0+2\pi]$, ($0\leq i\leq 2N$)
follow in alternating order,
$\bar\sa_0< \bar\sa_1<\bar\sa_2<\ldots <\bar\sa_{2N}:=\bar\sa_0+2\pi$,
 in particular,  {\sl $\bar\sa_i$ are
 relative maxima/minima points 
 for $i$ even/odd}. 
The corresponding  distinct critical energies will be denoted  by 
\beq{cesarini}
\bar E_i:=\GO(\bar\sa_i)\,, \quad
\bar E_{2N}=\bar E_0\ {\rm \ being\ the\ unique\  global\ maximum\ of\ } \GO\,.
\eeq
By the Implicit Function Theorem,
 for $\lalla$ small enough with respect to $\upkappa$,
 one can continue  the $2N$ critical points $\bar\sa_i$ of $\GO$    obtaining $2N$ 
critical points
$\sa_i=\sa_i(\hat{p})$
of $\Gm(\hat{p},\cdot)$ for $\hat{p}\in\hat D$.
The corresponding distinct  critical energies become
\beq{chinottocaldo}
E_i=E_i(\hat{p}):=\Gm(\hat{p}, \sa_i(\hat{p}))\,.
\eeq
Furthermore, for $\lalla$ small , 
 the functions
$\sa_i(\hat{p})$ and 
$E_i(\hat{p})$ preserve the same order
of $\bar\sa_i$ and $\bar E_i$. 
Indeed,  from Definition~\ref{buda} and the Implicit Function Theorem, the following   result  proven in \cite{BCaa23} holds\footnote{See Lemma 3.1 in \cite{BCaa23}. }:
\begin{lemma}\label{enricone}
Let $\Hpend$ be as in Definition~\ref{morso} and 
assume that\footnote{Notice that condition \equ{mino} is stronger than \equ{minore}.}
\beq{mino}
\lalla\leq 1/(2\upkappa)^6\,.
\eeq
Then,
the  functions   $\sa_i(\hat{p})$
and $E_i(\hat{p})$ defined above are  real analytic in $\hat{p}\in \hat D_\ro$
and
\begin{equation}\label{october}\ts
\sup_{\hat{p}\in \hat D_\ro}|\sa_i(\hat{p})-\bar\sa_i| \leq
\frac{2\suca\lalla}{\morse \so}\,,\qquad 
\sup_{\hat{p}\in \hat D_\ro} |E_i(\hat{p})-\bar E_i|
 \leq 3\upkappa^3 \suca\lalla\,.
\end{equation}
Furthermore, the relative order of  $\sa_i(\hat{p})$ and $E_i(\hat{p})$ is, for every $\hat{p}\in \hat D_\ro$, the same as that of, respectively, $\bar\sa_i$ and $\bar E_i$.
\end{lemma}
Therefore, under the assumption \equ{mino}, we see that the phase space
 $\cM$ is disconnected by the separatrices\footnote{I.e., the stable manifolds (curves) of the hyperbolic points $(0,\sa_{2j})$.} into  exactly $2N+1$ 
open connected components $\cM^i=\cM^i(\hat {p})$, for $0\le i\le 2N$, which can be labelled so that:

\begin{itemize}
\item[\bolla]the {\sl odd regions} $\cM^{2j-1}$ (for $1\le j\le N$) contain the elliptic points $(0,\theta_{2j-1})$ and have as boundary parts of separatrices; topologically, such regions are discs;

\item[\bolla]
the {\sl outer even regions} $\cM^{0}$ and $\cM^{2N}$ are 
homotopically {\sl non  trivial} annuli bounded by the most external separatrices  and one of the two  curves $\Hpend^{-1}(\Eb)$;

\item[\bolla]
when  $N>1$, the  {\sl inner even regions} $\cM^{2j}$ (for $1\le j\le N-1$) are  homotopically trivial annuli\footnote{I.e., annuli in the cylinder $\real\times\T$ which are contractible.} whose boundary is given by two pieces of separatrices (with different energies).
\end{itemize}

\centerline{\includegraphics[width=7cm,height=6.5cm]{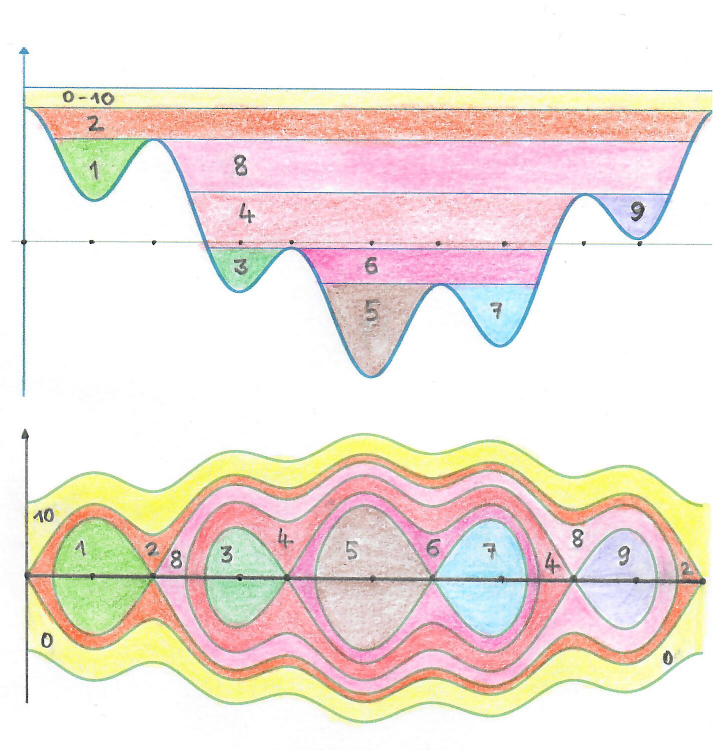}}

\nopagebreak

\centerline{\footnotesize \it  The Morse potential $\GO(q_1)=\sin q_1 + \frac12\cos(5q_1)$ with 10 critical points   (top)}
\centerline{\footnotesize \it and the phase portrait of $\Hpend:=p_1^2+\GO(q_1)$ (bottom)}
\centerline{\footnotesize \it Labels of corresponding regions are as in Definition \ref{piggy}}
\vglue0.3truecm

\nl
More formally, we can define the $2N+1$ regions $\cM^i$ in terms of suitable energy intervals $(E^{(i)}_-,E^{(i)}_+)$ as follows. 

\nl
Let $E_i$ be the critical energies defined in \equ{chinottocaldo}, and let $\Eb$ the reference energy defined in \equ{emmone}.

\dfn{piggy} {\rm (i) (Outer regions)} For $i=0,2N$, let $E^{(0)}_-=E^{(2N)}_-:= E_0$, and 
$E^{(0)}_+=E^{(2N)}_+:= \Eb$. Then, the `lower outer region' $\cM^{(0)}$ is the connected component of $\Hpend^{-1}\big((E^{(0)}_-,E^{(0)}_+)\big)$ contained in $\{p_1<0\}$, while 
the `upper outer region' $\cM^{(2N)}$ is the connected component of $\Hpend^{-1}\big((E^{(2N)}_-,E^{(2N)}_+)\big)$ contained in $\{p_1>0\}$.
\\
{\rm (ii) (Inner region, $N=1$)}  When $N=1$, $\cM^{(1)}$ is just the region enclosed by the unique separatrix
 $\Hpend^{-1}(E_0)$; the orbits in $\cM^{(1)}$ have energies ranging in the critical interval $[E^{(1)}_-,E^{(1)}_+):=[E_1,E_0)$.
\\
{\rm (ii) (Inner regions, $N>1$)} Define $E^{(i)}_-:=E_i$.\\
For $i$ odd, let $E^{(i)}_+:=\min\{ E_{i-1}, E_{i+1}\}$ and define $\cM^{(i)}$ as the connected component of 
$\Hpend^{-1}\big([E^{(i)}_-,E^{(i)}_+)\big)$ containing the elliptic equilibrium $(0,\sa_i)$.
\\
Finally, for $0<i=2j<2N$ even, define
 $$
j_-:=\max\{ \ell<j\big|\   E_{2\ell}>
 E_{2j} \}\,,\  
j_+:=\min\{ \ell>j\big|\ E_{2\ell}>
 E_{2j} \}\,,\ 
 E^{(i)}_+:=\min\{ E_{2j_-}, 
E_{2j_+}\}\,;
 $$
and define $\cM^{(i)}$ as the connected component of  $\Hpend^{-1}\big((E^{(i)}_-,E^{(i)}_+)\big)$ whose boundary contains the hyperbolic point $(0,\sa_i)$.
\edfn
Notice that the phase space $\cM$ is the union of the regions $\cM^{(i)}$ and the singular
 zero--measure set $S=S(\hat p)$
 formed by  the $N$ 
separatrices:
\beq{enkidu} 
\cM=\cM(\hat p)=\bigcup_{i=0}^{2N} \cM^i  \ \cup\ S =
\bigcup_{i=0}^{2N} \cM^i(\hat {p}) \ \cup\ S(\hat p)\,.
\eeq
Below we shall also consider the following $(n+1)$--dimensional domains:
\begin{eqnarray}\label{didone}
\check\cM &:=&
\{ (p,q_1)  \ \ {\rm s.t.}\ \ 
\hat p\in\hat D,\ (p_1,q_1)\in\cM(\hat p)
\}\,,\nonumber
\\
\check\cM^i &:=&
\{ (p,q_1)  \ \ {\rm s.t.}\ \ 
\hat p\in\hat D,\ (p_1,q_1)\in\cM^i(\hat p)
\}
\,.
\end{eqnarray}
Notice that $\bigcup_{0\leq i\leq 2N} \check\cM^i$ covers $\check\cM$ up to a   set of measure zero.

\subsubsection*{Arnol'd--Liouville's  action/energy functions}
Let  $E\in [E^{i}_-(\hat{p}),E^{i}_+(\hat{p})]$  and let $ \g^i$  be the (possibly, piece--wise) smooth closed curve in the 
clusure of $\cM^i(\hat{p})$
given by
\beqno
 \g^i= \g^i(E, \hat {p}):=
\{ ({p}_1,{q}_1)\in\overline{\cM^i(\hat{p})}\ \ {\rm s.t.}\ \ 
 \Hpend({p}_1,\hat{p},{q}_1)=E
\}\,,
\eeqno
 oriented clockwise\footnote{For the non contractible  curves ($i=0,2N$) the orientation is `to the right'
on $\cM^{2N}$, `to the left' on 
$\cM^0$.}; for $2\leq j\leq N$  consider also the trivial curves $\g^i_j=\{(p_j,s): s\in \torus\}$.  

\nl
Then,   
the  classical {\sl Arnol'd--Liouville's  action functions} are given by
\beqno
\begin{array}l
\dst \act_1^{(i)}(E)=\act_1^{(i)}(E,\hat p):= \frac1{2\pi}  
\oint_{ \g^i} p_1d q_1\,, \\
\dst  \act_j= \frac1{2\pi}  
\oint_{ \g^i_j}  p_jd q_j
=\frac{ p_j}{2\pi}\int_{\torus} d q_j =  p_j\,,\qquad \forall\ 2
\le j \le N\,.
\end{array}
\eeqno
The action  function
 $E\to\act_1^{i}(E,\hat\act)$ is  strictly monotone and its inverse is, by definition, the {\sl energy function
$\act_1\to {\mathtt E}^{i}(\act_1,\hat\act)$}. 
We also define   
$\bar\act_1^{i}:=\act_1^{i}|_{\lalla=0}$
and its inverse function\footnote{Note that when $\lalla=0$, $\Hpend$ becomes simply $\bHpend=p_1^2+\GO(q_1)$.}
$\bar{\mathtt E}^{i}:=
{\mathtt E}^{i}|_{\lalla=0}$.

\nl
We can now describe the fine analytic properties of the  action/energy functions.

\subsubsection*{Critical holomorphic behaviour and action estimates}

\nl
The first  result describes the exact  behaviour of the action functions as the energy approaches the critical energy of separatrices
and  contains  estimates on the derivatives  of the action functions that
will play a central r\^ole in the discussion on the twist Hessian matrix in \S~\ref{thisistheend}.  
The following theorem has been proven in  \cite[Theorem 3.1]{BCaa23}.

\begin{theorem}\label{glicemiak} Let $\Hpend$ be a   Hamiltonian in standard form as in Definition~\ref{morso},
let $\upkappa\ge 4$ be such that \equ{alce} holds and let $2N$ be the number of critical points of the reference potential $\GO$.
Then, there exists   a suitable  constant
$\bfc=\bfc(n,\upkappa)\geq  2^8\upkappa^3$ 
such that, if \footnote{Note that \equ{caviale2}
implies the hypothesis of Lemma \ref{enricone}. Thus, in particular also $\Hpend$ has $2N$ critical points.}
 \begin{equation}\label{caviale2}
\lalla\leq 1/\bfc^2\le 1/(2^{16}\upkappa^6)\,,
\end{equation}   
then, for all $0\le i\le 2N$ and $\hat \act\in\hat D$, the action functions $E\in (E^{i}_-(\hat \act),E^{i}_+(\hat \act))\mapsto \act_1^{i}(E,\hat \act)$ verify the following properties.

\nl
{\rm (i)} {\bf(Universal behaviour at critical energies)} 
There exist  functions $\phi^{i}_-(\z,\hat \act),$ 
$\psi^{i}_-(\z,\hat \act)$ for  $0\le i\le 2N$, 
and,  functions $\phi^{i}_+(\z,\hat \act),$ 
$\psi^{i}_+(\z,\hat \act)$, for $0<i<2N$,
 which are 
real analytic in a complex neighborhood of the set  $\{\z=0\}\times \hat D$ and satisfy
\begin{equation}\label{LEGOk}
\act_1^{i}\big(E_\mp^{i}(\hat \act)\pm \suca \z, \,\hat \act\big)=
\phi^{i}_\mp(\z,\hat \act) +\psi^{i}_\mp(\z,\hat \act)\ \z \log \z\, 
 \,,\quad \forall \ 0<\z<{ 1/\bfc}\,,\,
 \hat \act\in \hat D\,.
\end{equation}
the functions $\phi^{i}_\pm(\z,\hat \act)$,
$\psi^{i}_\pm(\z,\hat \act)$  are real analytic on 
 $\{\z\in\complex: |\z|< 1/\bfc\}\times  \hat D_{\ro}$, where 
satisfy:
\beq{pappagallok}
\begin{array}{l}
\dst
\sup_{|\z|<1/\bfc, \,\hat \act\in \hat D_{\ro}}\big(
|\phi^{i}_\pm|+|\psi^{i}_\pm|\big)
\leq \bfc\sqrt\suca 
\,,
\\
\dst
\sup_{|\z|<1/\bfc, \,\hat \act\in \hat D_{\ro/2}}
\big(|\partial_{\hat \act}\phi^{i}_\pm|
+|\partial_{\hat \act}\psi^{i}_\pm|\big)
\leq \bfc
\lella
\,,\qquad \lella:=
{\ts \frac{\sqrt\suca}{\ro}\lalla} \stackrel{\eqref{cimabue}}\leq 
2^{-8}\lalla
\,.
\end{array}
\eeq
Moreover,
 \begin{equation}\label{trippa}
 | \phi^{i}_\pm -  \bar\phi^{i}_\pm|\,,\ 
| \psi^{i}_\pm -  \bar\psi^{i}_\pm|\ \leq\ 
\bfc\sqrt\suca \lalla\,,
\end{equation}
where $\bar\phi^{i}_\pm:=\phi^{i}_\pm|_{_{\lalla=0}}$
and 
$\bar\psi^{i}_\pm:=\psi^{i}_\pm|_{_{\lalla=0}}$.

\nl
{\rm (ii)} {\bf (Limiting critical values)}
The following bounds at the limiting critical energy values hold:
\beq{ciofecak}
 \begin{array}{l}
|\psi^{i}_+(0,\hat \act)|\geq \sqrt\suca/ \bfc \,,
 \quad \phantom{r} 0< i< 2N\,,\quad  \forall \ \hat \act\in \hat D_{\ro}\,,\\
 |\psi^{2j}_-(0,\hat \act)|\geq \sqrt\suca/\bfc\,,
 \quad \ 0\leq j\leq N\,, \quad\   \forall  \ \hat \act\in \hat D_{\ro}\,,
 \\
 \psi^{i}_+(0,\hat \act)>0 \,,
 \quad \phantom{r} 0< i< 2N\,,\quad  \forall \ \hat \act\in \hat D\,,\\
 \psi^{2j}_-(0,\hat \act)<0\,,
 \quad \ 0\leq j\leq N\,, \quad\   \forall  \ \hat \act\in \hat D\,,
 \end{array}
\eeq
while, in the case of relative minimal critical energies, one has, $\forall$ $\hat \act\in \hat D$, $0<\z<{ 1/\bfc}$,
\begin{equation}\label{lamponek}
\phi^{2j-1}_-  (0,\hat \act)=0\,,
\qquad     
\psi^{2j-1}_-(\z,\hat \act)=0\,,
\qquad  \forall\ 1\le j\le N\,.
\end{equation}
\nl
{\rm (iii)}
{\bf (Estimates on  derivatives of  actions on real domains)}
The derivatives of the actions with respect to energy verify, on real domains, the following estimates:
\begin{equation}\label{vana} 
\inf_{( E^{i}_-, E^{i}_+)}
\partial_E  \act_1^{i} \ge \frac{1}{\bfc\sqrt\suca}\,,\qquad \forall\ \hat \act \in \hat D\,,\ \forall\ 0<i<2 N
\,;
\end{equation}

\beqno
\min\big\{\partial_E  \act_1^{2N}\,,
\
\partial_E  \act_1^{0}\big\}
\ge \frac{1}{\bfc\sqrt{E+\suca} }
\,,
\ \ 
\forall\,
E> E_{2N}\,,\  \forall\ \hat \act\in\hat D\,.
\eeqno
{\rm (iv)}
{\bf  (Estimates on  derivatives of  actions on complex domains and perturbative bounds)}
For $\loge>0$ satisfying
 \beq{bassoraTH}
  \bfc\lalla
  \leq
  \loge
  \leq 
  1/\bfc
  \,,
\eeq
define the following complex energy--domains:
\beq{autunno2}
{\mathcal E}^{i}_\loge :=
\left\{
\begin{array}{ll}
 \{z\in\complex:\bar E^i_- -{\ts \suca/\bfc}\ \,<\Re  z<\bar E^i_+ - \loge\suca\,,\   |\Im z|<{\ts \suca/\bfc} \}
\,,
&i\ {\rm odd}\,,
\\
 \{z\in\complex:\bar E^i_- + \loge\suca <\Re  z<\bar E^i_+ - \loge\suca\,,\  |\Im z|<{\ts \suca/\bfc}\}
\,,
&i\ {\rm even}\,, i\neq 0,2N
\,,
\\
 \{z\in\complex:\bar E^i_- + \loge\suca <\Re  z
 <\bar E^i_+
 \,,\  |\Im z|<{\ts \suca/\bfc} \}
\,,
& i=0,2N
\,.
\end{array}\right.
\eeq
Then, for $0\leq i\leq 2N$, the functions 
$
 \act_1^{i}$
and
$
 \bar \act_1^{i}$ are holomorphic 
on the domains
 $\mathcal E^{i}_\loge
\times\hat D_{\ro}$,  and satisfy the following estimates:
 \begin{equation}
\sup_{\mathcal E^{i}_\loge
\times\hat D_{\ro/2}}
|\partial_{\hat \act} \act_1^{i}|
\leq \bfc^2\, 
\lella
\,,
\ \ 
\sup_{\mathcal E^{i}_\loge}
\big|
\partial_E \bar \act_1^{i}
\big|
\leq \bfc^2\, 
\frac{|\log\loge|}{\sqrt\suca }
\,,\ \ 
\sup_{\mathcal E^{i}_\loge
\times\hat D_{\ro}}
\big|
\partial_E \act_1^{i}
-
\partial_E \bar \act_1^{i}
\big|
\leq 
\frac{\bfc^2\lalla}{\loge \sqrt\suca} \,
\label{rosettaTH}\,.
\end{equation}
\end{theorem}

\rem\label{caciottella}
(i)
Eq. \eqref{lamponek}  confirms the known analyticity at minima of actions as  function
of energy.

\nl
(ii)  A formula similar to \eqref{LEGOk} is given in  \cite{BFS} (compare Eq. (5.8) of Theorem 5.2 there).
\erem

\nl
We finally report a remarkable property of standard Hamiltonians $\Hpend$, whose reference potential 
$\GO$ is close enough to a cosine. In such a case, in fact,  one has uniform concavity of the second derivative of the energy function:

\begin{proposition}\label{kalevala}
 Assume that, for some $\sa_0\in\real$,   $\GO$ satisfies
 \beq{cosinelike}
{\modulo}\GO(\sa)-\cos (\sa+\sa_0) {\modulo}_1
 := \sup_{\torus_1}{\modulo}\GO(\sa)-\cos (\sa+\sa_0) {\modulo}
 \leq 2^{-40}\,.
\eeq
Then $N=1$ and 
\beqno
\partial^2_{\act_1} \bar{\mathtt E}^{1}           
(\bar \act_1^{1}(E))
\, \leq \,
-\frac{1}{27}\,,
\qquad
\forall E\in(\bar E_1,\bar E_2)
\,.
\eeqno
\end{proposition}
Also this result is proven in  \cite{BCaa23}; compare 
 Proposition 5.12 there.

\subsubsection*{Arnol'd--Liouville's action--angle variables in $n$ d.o.f.}
Let us now  discuss the  Arnol'd--Liouville's action--angle variables for the Hamiltonian $\Hpend$
viewed as a $n$ degrees of freedom Hamiltonian on the $2n$--dimensional phase space $\check\cM^i\times\T^{n-1}$.

\nl
For every fixed $\hat p\equiv\hat\act \in\hat D$,
the  map $(p_1,q_1)\to \act_1^{(i)}\big(\Hpend(p_1,\hat\act,q_1),\hat\act\big)$
can be symplectically completed with the angular term\footnote{Such completion is unique if one fixes, e.g., $\f^{(i)}_1(p_1,0;\hat\act)=0$.}
$(p_1,q_1)\to \f^{(i)}_1(p_1,q_1;\hat\act)=\f^{(i)}_1(p_1,\hat\act,q_1)$.\\
Defining the normal domains\footnote{Recall Definition~\ref{piggy}. For $i$ odd, 
$\act_1^{(i)}(E^i_-(\hat\act),\hat\act)=0$,
which is the action of the elliptic point.}
 \begin{equation}\label{aristakk}
\Bui :=\big\{ \act=(\act_1,\hat \act)\ |\ \hat \act\in\hat D,\ \  \act_1^{(i)}(E^i_-(\hat\act),\hat\act)<\act_1
<\act_1^{(i)}(E^i_+(\hat \act),\hat \act)
\big\}\,, 
\end{equation}
we see that, by construction,  the map\footnote{Recall the definition of $\check\cM^i $ in \equ{didone}.}
$$
(p,q_1)\in\check\cM^i \to (\act,\f_1)=
\big(\act_1^{(i)}(\Hpend(p,q_1),\hat\act),\hat\act,\f^{(i)}_1(p,q_1)\big)\in 
\Bui\times\T
$$
is surjective and invertible;  
let us denote by
\beqno
\cFiq:(\act,\ang_1)
\in 
\Bui\times\T
\ \to \ 
( p, q_1)\in\check\cM^i\,,\qquad (\hat p=\hat \act)\,,
\eeqno
its inverse map. 
Note that such
`Arnol'd-Liouville  suspended'
 transformation
$\cFiq$ integrates  $\Hpend$, i.e.,
\begin{equation}\label{red2}
\Hpend\circ\cFiq (\act,\ang_1)
={\mathtt E}^{(i)} (\act)\,,\qquad dp_1\wedge dq_1|_{\hat\act=\const}= d\act_1\wedge d\f_1\,.  
\end{equation}
By the standard Arnol'd--Liouville construction of the angle variables, one sees easily that the complete symplectic action--angle map 
$\Fiq:(\act,\ang)\mapsto (p,q)$ has the form
\beq{bruford2} 
\Fiq(\act,\ang)=\casitwo{(\upeta^i,\hat \act, \uppsi^i,\hat \ang+\upchi^i)\,,}
{0<i<2N\,,}
{(\upeta^i,\hat \act, \ang_1+ \uppsi^i,\hat \ang+\upchi^i)\,,}
{i=0,2N\,,}
\eeq
where $\upeta^i ,\upchi^i ,\uppsi^i $ are function of $(\act,\ang_1)$ only and are 
$2\pi$--periodic in $\ang_1,$
and, in the case $i=0,2N,$
$\sup|\partial_{\ang_1}\uppsi^i |<1$.

\nl
By construction, 
$\Fiq: \Bui\times\T^n\ \stackrel{\rm onto}{\longrightarrow}\ \check\cM^i\times\T^{n-1}$ {\sl is a global symplectomorphism}, and 
by \equ{red2}, one  has 
\begin{equation}\label{red3}
(\Hpend\circ \Fiq)(\act,\ang)=
(\Hpend\circ\cFiq) (\act,\ang_1)
={\mathtt E}^{(i)} (\act)\,,\qquad \forall\ 0\le i\le 2N\,.
\end{equation}
\nl
Next, we introduce suitable decreasing subdomains 
$\Bui (\loge)$ of $\Bui$ 
depending on a non negative parameter $\loge$ so that $\Bui(0)=\Bui$  and 
such that 
the map $\Fiq$ 
has, for positive $\loge$,  a holomorphic extension
on a suitable complex neighborhood
of $\Bui(\loge)\times\T^n$.

\nl
Define
\begin{equation}\label{sinistro}
\logemax=\logemax(\hat \act):=\big(E_+(\hat \act)-E_-(\hat \act)\big)/\suca\,,\qquad
\blogemax:=
\big(\bar E_+-\bar E_-\big)/\suca
\,.
\end{equation} 
Notice  that, by \equ{alce}, Definitions~\ref{buda}, \ref{morso}, and \equ{cimabue}  one has 
\beq{latooscuro}
1/\upkappa\leq 
\morse/\suca
\leq
\blogemax
\leq
2\,;
\eeq
notice also that, by  \eqref{october},  we have\footnote{Recall that $\lalla\leq 1/\bfc^2$ and  $\bfc \geq  2^8\upkappa^3$ (compare Theorem~\ref{glicemiak}).} 
\begin{equation}\label{destro}
|\logemax-\blogemax|\leq 6\upkappa^3\lalla\,,\qquad \logemax\geq 1/2\upkappa\,.
\end{equation}
Then, for   $0\le\loge\leq \logemax$
define\footnote{Recall the definition of $\Eb$ in \equ{emmone}.}:
\beqa{arista}   
\acci^i_\loge(\hat\act)
&:=&
\act_1^i(E_-^i(\hat \act)+\loge\suca,\hat \act)
\,,
\quad\ \ 
\forall\, 0\le i\le 2N\,,
\nonumber
\\
\bacci^i_\loge(\hat\act)
&:=&
\left\{
\begin{array}{ll}
\act_1^i(E_+^i(\hat \act)-\loge\suca,\hat \act)\,,& \forall\  0<i<2N\\
\act_1^i\big(\Eb,\hat \act\big)\,, &\, i=0,2N\,.  
\end{array}
\right.
\\
\acci^i(\hat\act)&:=&\acci^i_0(\hat\act)\,,\  \bacci^i(\hat\act):=\bacci^i_0(\hat\act)\,, \ 
 \ \forall\  0\le i\le 2N\,,
\nonumber\\
\Bui(\loge)&:=&\{
\act=(\act_1,\hat\act): \ \hat\act\in\hat D\,, \ 
\acci^i_\loge(\hat\act)<\act_1< \bacci^i_\loge(\hat\act)\}\,, \quad 0\le \loge\le\logemax\,.
\nonumber
\eeqa

\rem (i) By the above definitions one has that
\beq{accidenti}
\acci^{2j-1}(\hat \act):=\acci^{2j-1}_0(\hat \act)=\act_1^{2j-1}(E_-^{2j-1}(\hat \act),\hat \act)
\equiv 0\,,
\eeq
reflecting  the analyticity at the 
 elliptic points; compare Remark~\ref{caciottella}--(i) above. 
 \\
(ii)  
By \equ{aristakk} and \equ{arista} one sees that
$\Bui=\Bui(0) = \bigcup_{0<\l<\logemax}\Bui(\loge)$.  
\erem
The holomorphic properties of the Arnol'd--Liouville symplectic maps are described in  following theorem, proven in
\cite[Theorem 4.1]{BCaa23}. Recall the definition of the constant $\bfc$ in Theorem~\ref{glicemiak}.

\begin{theorem}\label{barbabarba}    
Under the  hypotheses of Theorem~\ref{glicemiak}
there exists a constant  
$\hbfc=\hbfc(n,\upkappa)\ge 4 \bfc^2$ 
depending only on
$n$ and $\upkappa$ such that,
taking
\beq{caviale3}
\lalla\leq 1/\hbfc\,,
\eeq
 the symplectic transformation $\cFiq$
 extends, for any  
$0\leq i\leq 2N$ and  
 $0< \loge\leq 1/\hbfc$,
to a real analytic map 
\begin{equation}\label{pediatra2bis}
\Fiq :\big(\Bui (\loge)\big)_{\!\rhol}\times\T^n_\sil\ \to\ 
D_{\ro}\times\T^n_{\so/4}\,,
\qquad \forall\, 0< \loge\leq 1/\hbfc\,,
\end{equation}
 where
\begin{equation}\label{blueeyes}\ts
\rhol	:=\frac{\sqrt\suca} {\hbfc}\, \loge |\log \loge| 
\,,\qquad \sil:= \frac1{\hbfc|\log \loge|}\,.
\end{equation}
Now, let $0< \loge\leq 1/\hbfc$, then
the function $\mathtt E^i$ admits a  holomorphic
extension
on $\big(\Bui (\loge)\big)_{\!\rhol}$, where, setting $\hat\loge:=\loge|\ln\loge|^3$, one has 
\beq{ofena}\ts
\big|\partial_{\act_1} \mathtt E^i \big|
\leq 
\hbfc\sqrt{\suca+|\mathtt E^i|}\,, 
  \quad 
\big|\partial^2_{\act_1} \mathtt E^i \big|
\leq 
\frac{\hbfc}{\hat\loge}\,,
\quad 
\big|\partial^2_{\act_1 \hat \act} \mathtt E^i \big|
\leq
\hbfc\frac{\lella}{\hat\loge}
\,,
 \quad 
\big|\partial^2_{\hat \act} \mathtt E^i \big|
\leq \hbfc
\big(
\frac{\sqrt\suca}{\ro}\act_1^{i}
+
 \frac{\lella}{\hat\loge}
\big)
\lella\,;
\eeq
furthermore, defining 
\beq{piatto}
D^\flat:=(-\Ro -\ro/3,\Ro+\ro/3)\times\hat D\,,\qquad \check\cM^i (\loge):= \cFiq (\Bui (\loge)\times\T)\,,
\eeq
one has
\beq{inthecourtk}
\meas\big((D^\flat\times\T)\ \setminus\ 
\bigcup_{0\leq i\leq 2N}\check\cM^i (\loge)\big)
\leq \hbfc\, \sqrt\suca
\meas(\hat D)\ \loge |\log \loge|\,.
\eeq
 \end{theorem}

\rem
Observe that, by \equ{cerbiatta}, \equ{roger}, \equ{athlone}, \equ{caviale2}, \equ{dublino}, \equ{tikitaka} and \equ{kappa}, it is\footnote{Recall that $\e<1$; see \equ{ham}.}
\beq{fango}\ts
1/\upkappa< \chs/4
\,,\qquad\quad
\frac{\e_k\chk}{\chr}\,\lalla<\ro/6
\,,\qquad\quad
 \frac{4\e_k\chk}{\chr^2}\,\lalla< \frac{\chs}{2^{20} \upkappa^3}< \chs/2^{20}\,.
\eeq
Thus, since\footnote{Recall the hypotheses of Theorem~\ref{barbabarba}.} $\loge\le 1/\hbfc$, by \equ{caviale2},
 $\sil$ in \equ{blueeyes} satisfies
\beq{cindy}
\sil<\chs/2^{20}\,.
\eeq
\erem

\section{Secondary nearly--integrable structure at simple resonances}\label{GS}
 
Now we go back to the original system in the simply--resonant zones governed by the Hamiltonians $\cH_k(\tty,\ttx)$ in \equ{dopomedia} and 
discuss their  global nearly--integrable structure with exponential small perturbations  (compare Theorem~\ref{garrincha} below). 

\nl
As mentioned above (see item (ii) in Remark~\ref{portogallo}), the problem here is that the symplectic transformations of  Theorem~\ref{sivori}, which put the  simply--resonant Hamiltonians $\Hsharp$ in \equ{tikitaka} in
 standard form, are, in general, {\sl not well defined in the fast angles} 
$\hat \ttq=(\ttq_2,...,\ttq_n)$, making the construction of {\sl global} action--angle variables for the full Hamiltonians $\cH_k(\tty,\ttx)$ in \equ{dopomedia} not straightforward. 
\\
To overcome such homotopy problems, we shall exploit the particular group structure of the various symplectic transformations involved, and show that, introducing  a special {\sl ad hoc} conjugacy, one can indeed  obtain globally well defined symplectic maps; see, in particular, \equ{elettra} below.

\subsubsection*{Special sets of symplectic transformations}

Besides the group $\Gdag$ introduced in Definition~\ref{dadaumpa0} above, we shall introduce   two new special classes of symplectic transformations,  
which will be used in the proof of Theorem~\ref{garrincha}.
Recall the notation in \equ{checheche}.

\dfn{dadaumpa} 
{\rm (a)}
Given a domain $\hat {\rm D}\subseteq \R^{n-1}$, 
$\Ggot$ denotes the formal\footnote{See Remark~\ref{alice}--(iii) below.} group
 of symplectic transformations of the form 
\beqno
(p,q)\in{\rm D}\times \torus^n\stackrel{\Phi}\mapsto (P,Q)=
(\upeta,\hat p, q_1+\uppsi,\hat q+\upchi)\in\real^n\times\T^n\,,
\eeqno
 where: ${\rm D}\subseteq \real^n$ is a normal smooth domain\footnote{I.e., $D=\{(p_1,\hat p):\ \a(\hat p)<p_1<\hat \b(\hat p)\,,\ \hat p \in \hat{\rm D}\}$ where $\a$ and $\b$ are smooth function on $\hat {\rm D}$.} over
 $\hat {\rm D}$, the functions $\upeta,\uppsi,\upchi$ depend on $(p,q_1)$, are
 $2\pi$--periodic in $q_1$ and the $(n+1)$--dimensional the map
\beqno
(p,q_1)\mapsto \check\Phi(p,q_1)=(\upeta,\hat p, q_1+\uppsi)
\eeqno
is injective. 

\nl
{\rm (b)} Given a domain $\hat {\rm D}\subseteq \R^{n-1}$, 
$\Ggoto$ denotes the set
 of smooth symplectic transformations of the form   
\beqno
(p,q)\in{\rm D}\times \torus^n\stackrel{\Phi}\mapsto (P,Q)=
(\upeta,\hat p, \uppsi,\hat q+\upchi)\in\real^{n+1}\times\torus^{n-1}\,,
\eeqno
 where  ${\rm D}\subseteq \real^n$ is a normal smooth domain over
 $\hat {\rm D}$;   the functions  $\upeta,\uppsi,\upchi$ depend only on $(p,q_1)$ and  are
 $2\pi$--periodic in $q_1$.
 \edfn
Let us collect a few  observations and discuss the main properties of such classes, but, first of all, 
notice that all the above maps leave fixed the variable $\hat p\in \hat {\rm D}\subseteq \R^n$ and the set $\rm \hat D$. Thus, in the following discussion,  the domain $\hat {\rm D}$ is fixed once and for all. 

\rem\label{alice} 
(i) The Arnol'd--Liouville map $\Fiq$ in the  outer cases  \equ{bruford2} ($i=0,2N$) belongs to $\Ggot$ (since 
$\sup |\partial_{q_1}\uppsi|<1$), 
while $\Fiq$ in the  inner case  \equ{bruford2} ($0<i<2N$) belongs to $\Ggoto$. \\
Notice also that $\Fidue$ in Theorem~\ref{sivori}--(iii) is a near--to--the--identity symplectic map belonging to $\Ggot$.

\nl
(ii) In the definition of $\Ggot$ and $\Ggoto$,  the functions $\upeta$ and $\uppsi$ are scalar functions, while $\upchi$ has $(n-1)$ components. Notice  that, since $\Phi$ is assumed to be symplectic,  these maps are such that 
\beqano
d\upeta\wedge dq_1+ d\upeta\wedge d\uppsi+ d  \hat p\wedge d \upchi&=&dp_1\wedge dq_1\,, \phantom{AAAAAAAe} (\Phi\in\Ggot)\,,
\\ 
d\upeta\wedge d\uppsi+ d  \hat p\wedge d \upchi&=&dp_1\wedge dq_1\,, \qquad\qquad \quad\ \  (\Phi\in\Ggoto).
\eeqano
(iii)  All maps in the group $\Gdag$ in Definition~\ref{dadaumpa0}  have a common domain of definition, i.e.,  
$(\R\times\hat {\rm D})\times \R^n$. 
On the other hand, every map  $\Psi\in \Ggot$  has its own domain of definition ${\rm D}$. Thus, the composition $\Psi_1\circ\Psi_2$ of two maps in $\Ggot$
$$\Psi_1:{\rm D}_1\times\torus^n\to\real^n\times\torus^n\,,\qquad \ 
\Psi_2:{\rm D}_2\times\torus^n\to\real^n\times\torus^n$$
is well defined {\sl only when the compatibility condition $\Psi_2\big({\rm D}_2\times\torus^n\big)\subseteq {\rm D}_1\times\torus^n$ is satisfied}. 
This is the reason why the cautionary word `formal' appears in the definition of $\Ggot$. However, as already noticed, all maps in $\Ggot$ verify
$\pi_{\hat p}({\rm D})=\hat{\rm D}$, which is fixed {\sl a priori}.

\nl
(iv)   If  $\Phi\in\Ggot$, by definition $\check\Phi$ is injective, so that 
also $\Phi$ itself is injective.
Furthermore, for any fixed $p$,  the map $q_1\to Q_1=q_1+\uppsi$ is a continuous  injective map on the circle $\torus^1$, 
hence it is surjective, and, therefore, it is a smooth (orientation preserving) circle diffeomorphism. Thus,    {\sl $q\to Q=(q_1+\uppsi,\hat q+\upchi)$ is a global diffeomorphism of $\torus^n$}, and {\sl $\Phi:{\rm D}\times\torus^n\to \Phi({\rm D}\times\torus^n)\subseteq\R^n\times\T^n$ is a global symplectomorphism}.
\\
Notice also that 
{\sl if $\Phi,\Phi'\in\Ggot$ and the composition $\Phi\circ\Phi'$ is well defined, then $\Phi\circ\Phi'\in \Ggot$}.

\nl
(v) The definition of the first $(n+1)$ component of any member of the above families depends only on the first $(n+1)$ variables $(p,q_1)$.
Therefore, any finite compositions of maps 
$\Psi_i\in \Gdag\cup\Ggot\cup\Ggoto$, $1\le i\le m$, {\sl whenever the composition is well defined}, satisfies
\beq{yomo2}  
\Psi_i\in \Gdag\cup\Ggot\cup\Ggoto
\quad
\Longrightarrow
\quad
(\Psi_1\circ\cdots \circ \Psi_m)^{{\!\!\widecheck{\phantom{a}}}} =(\check\Psi_1\circ\cdots \circ \check\Psi_m)\,.
\eeq
(vi) Finally, one readily verifies that the following property holds:
\beq{buddah}
\Phi\in \Ggoto\quad {\rm and}\quad \Psi\in\Gdag\cup\Ggot
\qquad
\Longrightarrow
\qquad
\Psi\circ \Phi\in\Ggoto\,.
\eeq
\erem

\subsubsection*{Action--angles variables for the secular standard Hamiltonians $\Hsharp$ at simple resonances}
For each $k\in\genKO$, we  may apply the  theory of \S~\ref{criseide} to the secular Hamiltonians described in Theorem~\ref{sivori}  in standard form $\Hpend=\Hsharp$; see
 \equ{tikitaka2}, \equ{paranoia}, and \equ{kappa}.

\nl
By \equ{red3}, we get that, for every $k\in\genKO$ and $0\le i\le 2N_k$,  the Arnol'd--Liouville map 
\beq{freddy2}
\Fiq\ :\ \Buik\times\T^n\ \stackrel{\rm onto}{\longrightarrow}\ \check\cM^i_k\times\T^{n-1}
\eeq
integrates $\Hsharp$, i.e.:
\begin{equation}\label{red3bis}
(\Hsharp\circ \Fiq)(\act,\ang)=
(\Hsharp\circ\cFiq) (\act,\ang_1)
={\mathtt E}^{(i)}_k (\act)\,,\qquad \forall\ 0\le i\le 2N_k\,,
\end{equation}
where $\Buik$,  $\check\cM^i_k$ and ${\mathtt E}^{(i)}_k$
correspond to  $\Bui$,  $\check\cM^i$ and ${\mathtt E}^{(i)}$ in \S~\ref{criseide} in the case\footnote{Compare, in particular,
\equ{enkidu} and \equ{didone}  for the definitions of $\cM^i_k$ and  $\check\cM^i_k$; 
\equ{aristakk} for the definition of $\Buik$;   \equ{arista} and \equ{sinistro} for the definition of $\Buik(\loge)$;
  the definition of  $\check\cM^i_k(\loge)$ is given \equ{piatto}.} $\Hpend=\Hsharp$. 
  
\nl
Beware that, even if sometimes, for ease of notation,  we do not report the dependence upon the resonance label $k\in\genKO$,  
 we are treating different Hamiltonians in the neighbourhoods of simple resonances labelled by $k\in\genKO$.

\nl
Finally, we shall use the following notations:
Given a  function  $\ta: \hat {\rm D}\to \real$, we shall denote by $j_{\!\ta}$  
 the translation
\beq{tordissimo}
j_{\!\ta}(p):=(p_1+\ta(\hat p),\hat p)\,. 
\eeq
Notice that, by the definition of $\Psi_{\! \ta}$ in \equ{psiab}, one has 
\beq{trance}
 \check\Psi_{\! \ta}(p,q_1)= (j_{\!\ta}(p),q_1)\,.
\eeq

\subsubsection*{Global action--angle variables at simple resonances}
We are now ready to state and prove
the first step of the proof of Theorem~\ref{prometeo}, which consists in showing how to construct  
symplectic action--angle maps which put
a  generic nearly--integrable natural systems, near simple resonances, {\sl for all $k\in \genKO$}, into uniform analytic nearly--integrable form  {\sl with exponentially small perturbations}:

\nl
Let Assumptions~\ref{assunta} and Definitions~\ref{assunto} hold; let   $\bfco$ be as  in Theorem~\ref{normalform}, and $\hbfc$  as in Theorem~\ref{barbabarba} 
with $\upkappa$ as in \equ{kappa}. Let $\giuno$ and $\gitre$ be as  in (ii) of Theorem~\ref{sivori}, and define 
\beq{tordobis}
\Bik:=
\left\{ \begin{array}{ll} \Buik\,, &
 {\rm if} \ \ 0<i<2N_k\,, \\  
\Istarin \big(\Buik\big)\,, & {\rm if} \ i=0,2N_k\,,\ \end{array}\right. 
\quad  \gstar:=-(\giuno+\gitre)\,.
\eeq
Then, the following result holds.

\begin{theorem}\label{garrincha} 
{\bf (Secondary nearly--integrable structure at simple resonances)}\\
There exists
$\bfcast=\bfcast(n,s,\b,\d)\ge\max\{\itcd, \bfco,\hbfc\}$
such that 
if $\KO\ge\bfcast$, 
then for any $k\in\genKO$, $0\le i\le 2N_k$,
there exist 
 real analytic symplectomorphisms\footnote{Recall the notation \equ{dire}.}
\beq{tordor}
\upphi_k^i:  \Bik \times \T^n \to \Re(\Ruk_{r_k}) \times \T^n\,,
\eeq 
such that, if ${\mathtt E}^{i}_k={\mathtt E}^{i}_k(\act)$  is  the  integrable  Hamiltonian $\Hsharp$   of Theorem~\ref{sivori}
in its  Arnol'd--Liouville action
variables,   $\widetilde{\mathtt E}^{i}_k:= {\mathtt E}^{i}_k\circ \Istar$, and 
$\hzk$ is as in Theorem~\ref{sivori}, then\footnote{Recall \equ{sonnosonnoBIS} and \equ{trota2}.}
\beq{kant}
\begin{array}{l}
\cH^i_k:=\ham\circ \upphi_k^i(\act,\ang)
=
\hik(\act)+
\e \fik(\act,\ang)\,,\quad \mbox{with:}\\
\hik:={\ts \frac{|k|^2}{2}} \mathtt{h}^i_k\,,\ 
\mathtt{h}^i_k:=
\casitwo{
{\mathtt E}^{i}_k +\hzk
\,,}{0<i<2N_k\,,}
{
\widetilde{\mathtt E}^{i}_k  +\hzk
}{i=0,2N_k\,.}
\end{array}
\eeq
Furthermore,  for   $0<\loge\le 1/\bfcast$  define:
\beqa{tales}
&& \ts \rhs	:={ \frac{\sqrt{\suca}} {\bfcast\KO^{n}}}\  \loge |\log \loge|  
\,,\quad { \shs := \frac{1}{\bfcast\KO^{n}|\log \loge|}}\,,
\nonumber
\\ 
&& \Bik(\loge):=
\left\{ \begin{array}{ll} \Buik(\loge)\,, &
 {\rm if} \ \ 0<i<2N_k\,, \\  
\Istarin \big(\Buik(\loge)\big)\,, & {\rm if} \ \ i=0,2N_k\,,
 \end{array}\right. \quad \forall\ 0\le \loge< 1/\bfcast\,.
\eeqa
Then,
$\upphi_k^i$ admits a  holomorphic extension
\beq{tordo}
\upphi_k^i:  (\Bik(\loge))_{\rhs}\times \T^n_{\shs} \to \Ruk_{r_k} \times \T^n_{s_\varstar}\,
\eeq 
and  the perturbation $\fik$ in \equ{kant} satisfies the exponential estimate
\beq{schiaffino}
\sup_{(\Bik(\loge))_{\rhs}\times \T^n_{\shs} }|\fik|
\leq e^{-\K s/3}\,.
\eeq    
\end{theorem}

\rem
(i) 
Notice  that, since $\lalla=1/\K^{5n}$ (see \equ{cerbiatta}), and since\footnote{The constant $\bfc$ is defined in Theorem~\ref{glicemiak}.} 
$$\K>\KO\ge \bfcast> \bfc\,,$$ 
condition \equ{caviale3} -- which is stronger than condition \equ{caviale2} -- is implied by the assumption $\KO\ge\bfcast$. 
\\
Observe also that from the definitions of the constants in  Theorem~\ref{garrincha}, Theorem~\ref{glicemiak} and from \equ{destro} it follows that  
\beq{sgallettata}
\bfcast\ge\bfc\geq  2^8\upkappa^3\ge 2^{14}\,,\qquad \logemax\ge 2^{12}/\bfcast\,.
\eeq
Finally, we remark that, recalling the definitions of  $\rhol$ and $\sil$  in \equ{blueeyes}, since $\bfcast\ge \hbfc$, one has 
\beq{sillogismo}
 \rhs<\rhol\,,\qquad  \shs<\sil\,.
\eeq
(ii) In the proof of the theorem the maps $\upphi_k^i$ are explicitly given; compare \equ{agamennone} and \equ{menelao} below. 
\erem

\nl
The following simple lemma will be one of the key points of the proof of Theorem~\ref{garrincha}. 
Recall Definition~\ref{dadaumpa}.

\begin{lemma}\label{gattopazzo}
 Let $\Phi:(p,q)\in{\rm D}\times \torus^n\mapsto (\upeta,\hat p, q_1+\uppsi,\hat q+\upchi)\in 
\R^n\times\T^n$ be  in $\Ggot$,  $\Psi_{\!\ta}\in \Gdag$, and denote by $\t_{\!\ta} \Phi$ the map 
\beq{calcante}
\tau_{\!\ta} \Phi:=\t_{\!\ta} \Phi(p,q)
:=\big(
\upeta_\ta+\ta,\hat p,q_1+\uppsi_\ta, \hat  q+\upchi_\ta-\uppsi_\ta \partial_{\hat p} \ta
\big)\,,
\eeq
where  for a function $u:D\times\torus\to\R^m$, $u_\ta$ denotes te map
\beq{calcante2}
u_\ta:=u\circ\check \Psi_{\!-g}: j_{\!\ta}(D)\times\torus\to \R^m\,.
\eeq
Then,  $\tau_{\!\ta} \Phi$ belongs to $\Ggot$ and it  is a symplectomorphism satisfying
 \beq{atreo}
\tau_{\!\ta} \Phi: j_{\!\ta}({\rm D})\times\T^n\stackrel{\rm onto}\longrightarrow \big(\check\Psi_{\!\ta}\circ\Phi({\rm D}\times \T^n)\big)\times \T^{n-1}\,,
\eeq
and
\beq{atreo2}
(\tau_{\!\ta} \Phi)^{{\!\!\widecheck{\phantom{a}}}}=(\upeta_\ta+ \ta,\hat p,q_1+\uppsi_\ta)=\check\Psi_{\!\ta}\circ \check\Phi\circ \check\Psi_{\!-\ta}\,.
\eeq
\end{lemma}
\proof
First observe that  since $\upeta_\ta,\uppsi_\ta,\upchi_\ta$
are $2\pi$-periodic in $q_1$, the map
$$q\in\torus^n\mapsto \pi_{{}_Q} \tau_{\!\ta} \Phi(p,q)=
\big(q_1+\uppsi_\ta, \hat  q+\upchi_\ta-\uppsi_\ta \partial_{\hat p} \ta\big)\in\torus^n
$$
is a well defined  $\T^n$--map and \equ{atreo2} follows immediately by direct computation. 
Thus, $(\tau_{\!\ta} \Phi)^{{\!\!\widecheck{\phantom{a}}}}$ is injective being the composition of three injective maps, and, therefore, 
the whole map $\tau_{\!\ta} \Phi$ is injective, and  \equ{atreo} follows. To check symplecticity, just  note that, locally, on the universal cover 
$\R^{2n}$, $\tau_{\!\ta} \Phi$ coincides (as  it is immediate to check) with the composition $\Psi_{\!\ta}\circ\Phi\circ \Psi_{\!-\ta}$ of three symplectic maps. Hence $\tau_{\!\ta} \Phi$  is symplectic and the claim follows. \qed

\proof {\bf of Theorem~\ref{garrincha}} 
We start by defining   the maps $\upphi_k^i$.\\
Consider, first, the inner case $0<i<2N_k$. Recall Definition~\ref{dadaumpa}.
By Theorem~\ref{sivori}--(iii), $\Phi_\diamond$ is the composition of maps in 
 $\Gdag$ and $\Ggot$ while,  for $0<i<2N_k$, $\Fiq\in \Ggoto$ (Remark~\ref{alice}--(i)).
 Hence, by   \equ{buddah},  it follows that $\Phi_\diamond\circ \Fiq\in \Ggoto$ 
 and we may define\footnote{$\Psi^k$ appears in Theorem \ref{normalform}. Recall that, when $0<i<2N_k$,  $\Bik:=\Buik$.}
 \beq{agamennone}
\Phi^i_\diamond:=\Phi_\diamond\circ\Fiq\,,\quad 
\upphi_k^i:=\Psi^k\circ\Phi^i_\diamond:\Bik\times\torus^n\to \R^n\times\T^n\,,\quad\quad  (0<i<2N_k)\,,
\eeq
provided the composition is well defined.  To check that this is the case, we observe that  by  
 \equ{yomo2}, 
 \equ{pediatra2bis},   
 \equ{fangorn},
  \equ{tikitaka}, \equ{sillogismo}  for $0< \loge\le 1/\hbfc$,
we get
\begin{equation}\label{ossiriand}
\check\Phi^i_\diamond=
(\Phi_\diamond\circ
\Fiq){\!\check{\phantom{a}}}
=\check\Phi_\diamond\circ\cFiq
:\, 
(\Buik (\loge))_{\!\rhol}\times\T_\sil\ \to\ 
\DDD^k_{\tilde r_k}\times \T_\chs\,,\ \quad  (0<i<2N_k)\,,
\end{equation}
thus the composition is well defined and \equ{agamennone} is well posed.
\\
Let us now consider the outer case $i=0,2N_k$. In this case  $\Fiq\in \Ggot$  (Remark~\ref{alice}--(i)).
Recalling the definition in \equ{calcante}--\equ{calcante2}, by Lemma~\ref{gattopazzo}, we may define
\beq{elettra}
\Phi^i_{\!{}_{23}}:= \Fidue\circ \tau_{\!\gitre}\!\Fiq\,,\qquad{\rm and}\qquad \Phi_\diamond^i:=\tau_{\!\giuno} \! \Phi^i_{\!{}_{23}}\,,
\qquad\ \ \ (i=0,2N_k)\,.
\eeq
Recalling that $\Phi_2\in \Ggot$, by Lemma~\ref{gattopazzo} and
Remark~\ref{alice}--(iv),  $\Phi^i_{\!{}_{23}}\in \Ggot$ and, again by Lemma~\ref{gattopazzo}, $\Phi_\diamond^i\in \Ggot$, provided the compositions are well defined. To check that  this is the case, as above, it is enough to control the complex domains of the first $(n+1)$ components.
By \equ{atreo2} (used twice), \equ{cotugno}, \equ{toto},  and \equ{yomo2}, one finds\footnote{Recall that $\gstar=-(\giuno+\gitre)$; compare \equ{tordobis}.} 
\beq{tantalo}
\check \Phi_\diamond^i =\check\Phi_\diamond\circ \cFiq\circ\check\Psi_{\!\gstar}\,,\qquad \qquad \qquad (i=0,2N_k)\,.
\eeq
Then, by \equ{antiloco},  we get,  
\begin{equation}\label{cioccolata}
j_{\!\gstar}\big((\Bik(\loge))_{\prhol}\big)
\subseteq
\big(j_{\!\gstar} \big(\Bik(\loge)\big)\big)_{\rhol}
\stackrel{\equ{tales}}=\ts
\big(\Buik(\loge)\big)_\rhol\,,\quad 
{\rm where}\ \prhol:=\frac{\rhol}{n+2}\,,\qquad 
(i=0,2N_k)\,.
\end{equation}
Observing that 
$\check\Psi_{\!\gstar}(p,q_1)=\big(j_{\!\gstar} (p),q_1\big)$,
by \equ{tantalo}, \equ{cioccolata},  \equ{pediatra2bis}, 
\equ{fangorn} and
\equ{tikitaka},  
we get, for $0< \loge\le 1/\hbfc$,
\begin{equation}\label{ossiriand2}
\check\Phi^i_\diamond:\, (\Bik (\loge))_{\prhol}\times\T_\sil\ \to\ 
\DDD^k_{\tilde r_k}\times \T_\chs\,,\qquad (i=0, 2N_k)\,.
\end{equation}
Thus,   the composition is well defined and \equ{elettra} is well posed. So, we may define:
\beq{menelao}
\upphi_k^i:=\Psi^k\circ\Phi_\diamond^i:\Bik\times\torus^n\to \R^n\times\T^n\,,
\qquad   \Phi_\diamond^i\ {\rm as\ in}\  \equ{elettra}\,,
\quad\quad (i=0,2N_k)\,.
\eeq
We can, now,  prove \equ{kant}.
Recall the definition of $\bar f^k$ in Theorem~\ref{normalform}  and  
define
\beq{perturbazione}
 \fik:=f\circ \upphi_k^i\stackrel{\equ{dopomedia}}=\bar f^k\circ\Phi_\diamond^i\,,\qquad\quad \qquad(0\le i \le 2N_k)\,.
\eeq
Then, by definition of $\upphi_k^i$ in \equ{agamennone} and \equ{menelao}, we have, for $0\le i\le 2N_k$, 
\beq{peleo}
\cH^i_k:=\ham\circ \upphi_k^i(\act,\ang):=
\ham\circ \Psi^k\circ\Phi_\diamond^i 
\stackrel{(\ref{dopomedia},\ref{perturbazione})}=
\hamsec_k\circ \Phi_\diamond^i+ \e \fik\,.
\eeq
Since $\hamsec_k$ in \equ{hamseccippa} depends only on the first $(n+1)$ variables, by \equ{ossiriand} and \equ{tantalo}, we find
\beq{aloe}
\hamsec_k\circ\Phi_\diamond^i=\hamsec_k\circ\check\Phi_\diamond^i=
\casitwo{\hamsec_k\circ\check\Phi_\diamond\circ \cFiq\,,}{0<i<2N_k}{\hamsec_k\circ\check\Phi_\diamond\circ \cFiq\circ\check\Psi_{\!\gstar}\,,}{i=0,2N_k\,,}
\eeq
and, by \equ{tikitaka} and \equ{red3bis},
\beq{aloe2}
\hamsec_k\circ\check\Phi_\diamond\circ \cFiq=\textstyle{\frac{|k|^2}{2}}( {\mathtt E}^{(i)}_k+
\hzk)\,.
\eeq
Thus, \equ{kant} follows from \equ{peleo}, \equ{aloe}, \equ{aloe2} and \equ{trance}.

\nl
Next, we show that $\upphi_k^i$ has, for $0<\loge\le 1/\bfcast$,  a holomorphic extension satisfying \equ{tordo}.
To do this we have to consider the last $n-1$
components of $\Phi^i_\diamond$, namely\footnote{Recall the notation in \equ{checheche}.} $\pi_{\!{}_{\hat\ang}}\Phi^i_\diamond
=\hat\Phi^i_\diamond$.
By definition of $\Phi_\diamond^i$ in \equ{agamennone} and \equ{elettra}  it follows that
\beqa{cicciona} 
&&
\hat \Phi^i_{\diamond}(\act,\ang)=
\left\{ \begin{array}{ll} 
\hat \ang +\upchi^i_{\,\diamond}(\act,\ang_1)\,,
 &
\  {\rm if} \ 0< i<2N_k\,, \\ 
\hat \ang +\upchi^i_{\,\diamond}\big(\Istar(\act),\f_1\big)\,,& \  {\rm if} \  i=0,2N_k\,,
 \end{array}\right.
\eeqa
with\footnote{Recall the form of $\Fiq$ in \equ{bruford2};   $\gstar$ is defined in \equ{tordobis}; $\upchi_{{}_2}$ is as in Theorem~\ref{sivori}--(iii).}
\beq{ciccione}
\ts
\upchi^i_{\,\diamond}
:=
\upchi^i 
+\upchi_{{}_2}^\flat
+\uppsi^i\,
\partial_{\hat \act} \gstar\,,\qquad \upchi_{{}_2}^\flat(\act,\ang_1):=
\left\{ \begin{array}{ll} 
\upchi_{{}_2}(\hat\act,\uppsi^i)\,,
 &
\  {\rm if} \ 0< i<2N_k\,, \\ 
\upchi_{{}_2}(\hat\act,\ang_1+\uppsi^i)
\,,& \  {\rm if} \  i=0,2N_k\,.
 \end{array}\right.
\eeq
Now, we claim that 
\beq{fame}
|\uppsi^i|_{\rhol,\sil} <\frac34 \chs\,, \qquad \forall\ 0\le i \le 2N_k\,.
\eeq
Indeed, if $0<i<2N_k$, \equ{fame} follows directly from \equ{pediatra2bis} and \equ{fangorn}; 
in the case  $i=0,2N_k$, \equ{fame} follows again from \equ{pediatra2bis} and \equ{fangorn} observing that
$$
|\uppsi^i|_{\rhol,\sil}=|(\ang_1 +\uppsi^i)  -\ang_1|_{\rhol,\sil}
\le \frac{\so}4+\sil<\frac34 \chs\,.
$$
Next, since $\prhol=\rhol/(n+2)$, 
by \equ{cioccolata}, \equ{ciccione}, \equ{pediatra2bis}, \equ{fangorn},  \equ{tess2}, \equ{fango}, \equ{cindy}, \equ{fame}, \equ{antiloco}, we find,
for every $0\le i\le 2N_k$, and for every  $2\le \ell\le n$, 
\beqa{ciccioni}\nonumber
|\Im \hat\Phi^i_{\diamond \ell}|_{\prhol,\sil}
& \le& |\Im (\ang_\ell+  \upchi^i_{\,\diamond \ell})|_{\rhol,\sil}
\\
& \le&
|\Im(\ang_\ell+  \upchi^i_{\ell})|_{\rhol,\sil}+|\upchi_{{}_2}^\flat|_{\rhol,\sil}+
|\uppsi^i|_{\rhol,\sil}\,
|\partial_{\hat \act} \gstar|_\rhol\,\nonumber\\
&\le& \frac{\chs}2 +\frac{\chs}{2^{20}}+\frac34 (n+1)\chs<2n \chs\,.
\eeqa
Thus, by \equ{ossiriand}, \equ{ossiriand2} and \equ{ciccioni}, we get
\beqno
\Phi^i_\diamond:\, (\Bik (\loge))_{\prhol}\times\T_\sil^n\ \to\ 
\DDD^k_{\tilde r_k}\times \T_{2n\chs}^n\,,\qquad (0\le i\le  2N_k)\,.
\eeqno
We need, now, an elementary result on real analytic functions,  whose  proof is given in \Appendix:
\begin{lemma}\label{schiena}
Let $g: D_r\times \T_s^n\to \C$ be a  real analytic function 
 satisfying $|\Im g|\leq \xi$. 
Then, for every $0<\zeta \leq 1/2$, one has
$$
\sup_{D_{{}_{\zeta r}}\times \T_{{}_{\zeta s}}^n}|\Im g|\leq 8 \zeta \xi\,.
$$
\end{lemma}
Now, define 
\beq{gomma}
\zeta:= \frac1{16 n\, \itcu \ttcs\, \KO^n}\,.
\eeq
Then, since $|k|\le \KO$, by \equ{tikitaka}, \equ{cerbiatta},  we find
\beqno
8\zeta (2n \chs)< 16n\, \zeta\, \KO \max\{1,s\}\stackrel{\equ{gomma}}=\frac{\max\{1,s\}}{\itcu\ \ttcs\, \KO^{n-1}}
= \frac{s}{\itcu \KO^{n-1}}\stackrel{\equ{dublino}}=\tilde s_k\,.
\eeqno
Thus, by Lemma~\ref{schiena} (applied with  $g=\hat\Phi^i_{\diamond \ell}$  for $2\le \ell\le n$, $\zeta$ as in \equ{gomma} and $\xi=2n \chs$),
it follows that
\beqno
\Phi^i_\diamond:\, (\Bik (\loge))_{\rhs}\times\T_\shs^n\ \to\ 
\DDD^k_{\tilde r_k}\times \T_{\tilde s_k}^n\,,\qquad (0\le i\le  2N_k)\,,
\eeqno
with $\rhs$ and $\shs$ as in \equ{tales},  provided  
\beqno
\bfcast:= \max\{\itcd\,,\bfco\,,\, \hbfc \, \itcu\,\ttcs \, 16 n\, (n+2)\}\,. 
\eeqno
In conclusion, \equ{tordo} follows  by the definition of $\upphi_k^i$ in \equ{agamennone}, \equ{menelao} and by
\equ{trota2}. 

\nl
Finally, estimate \equ{schiaffino} follows at once from \equ{perturbazione}, \equ{tordo} and \equ{cristinacippa}. The proof is complete.
\qed

\nl
The following measure estimate will play a crucial r\^ole in the proof of Theorem~\ref{prometeo}.  

\begin{proposition}
For every $0\le \loge< 1/\bfcast$,
the following  measure estimate holds\footnote{The sets $\Rukt$ are defined in \equ{codino}.}:
\begin{equation}\label{inthecourt2}
\meas\big(
\big(\Ruk\times\T^n\big)\ \setminus\ 
\bigcup_{0\leq i\leq 2N_k}
\upphi_k^i
\big(\Bik(\loge)\times \T^n\big)
\big)
\le
\bfcast
\meas\big(\Rukt\times\T^n\big)\,
\loge |\log \loge| \,.
\end{equation}
\end{proposition}
\proof
Since $\check\Phi_\diamond^i$ depends only on the first $(n+1)$ variables, 
by \equ{cicciona}, \equ{agamennone}, \equ{tantalo} and the definitions of $\Bik(\loge)$ in \equ{tales} and $\check{\cal M}^i_k(\loge)$ in \equ{freddy2}, one has
\beqa{mauritania}
\Phi_\diamond^i \big(\Bik(\loge)\times\T^n \big)&=&\check\Phi_\diamond^i \big(\Bik(\loge)\times\T \big) \times \T^{n-1}
=\big(\check\Phi_\diamond\circ \cFiq(\Bu_k^i(\loge)\times\T)\big) \times \T^{n-1}\nonumber\\
&\stackrel{\equ{inthecourtk}}=& 
\big(\check\Phi_\diamond\circ \check{\cal M}^i_k(\loge)\big) \times \T^{n-1}\,.
\eeqa
Analogously, one has
\beq{mali}
\Phi_\diamond^{-1}(\DDD^k\times \T^n)=\check\Phi_\diamond^{-1}(\DDD^k\times \T)\times\T^{n-1}\,.
\eeq
Observe also that, by \equ{tess}, \equ{tess2} and (the second estimate in) \equ{fango} it follows 
that\footnote{Observe that $\cFitre^{-1}(p,q)=(p_1-\gitre(\hat p),\hat p,q_1)$ and 
$\cFidue^{-1}(p,q)=(p_1-\upeta_{{}_2}(\hat p, q_1),\hat p, q_1)$.
Recall the definition of $D^\flat$ in \equ{piatto}.}
\beq{chad}
\cFitre^{-1}\circ \cFidue^{-1} (D\times \T)\subseteq \big((-\Ro -\ro/3,\Ro+\ro/3)\times\hat D\big)\times\T=(D^\flat\times\T)\,.
\eeq
Then\footnote{The unions are over $0\le i\le 2N_k$.},  recalling  Theorem~\ref{normalform},  using the fact that $(\Psi^k)^{-1}$ and $\Phi_\diamond^{-1}$ are diffeomorphysms preserving  Liouville measure, we find
\beqano
&&\!\!\!\!\!\!\!\!\!\!\!\!\!\!\!\!\!\!\!\!\!\!\!\!
\meas(
\Ruk\times\T^n\ \setminus\ 
\mathsmaller{\bigcup}
\upphi_k^i
(\Bik(\loge)\times \T^n))\\
&\stackrel{(\ref{agamennone},\ref{menelao})}=&
\meas(
(\Psi^k)^{-1}(\Ruk\times\T^n)\ \setminus\ 
\mathsmaller{\bigcup}
\Phi_\diamond^i
(\Bik(\loge)\times \T^n))
\\
&\stackrel{\equ{surge}}\le&
\meas(
(\DDD^k\times\T^n)\ \setminus\ 
\mathsmaller{\bigcup}
\Phi_\diamond^i
(\Bik(\loge)\times \T^n)
)\\
& =&
\meas(
\Phi_\diamond^{-1}(\DDD^k\times\T^n)\ \setminus\ 
\mathsmaller{\bigcup}
\Phi_\diamond^{-1}\Phi_\diamond^i
(\Bik(\loge)\times \T^n))
\\
&\stackrel{(\ref{mauritania},\ref{mali})}=&(2\pi)^{n-1}
\meas(
\check\Phi_\diamond^{-1}(\DDD^k\times\T)\ \setminus\ 
\mathsmaller{\bigcup}
\check{\cal M}^i_k(\loge))
\\
&\stackrel{(\ref{toto},\ref{chirone})}=&(2\pi)^{n-1}\meas( \cFitre^{-1}\circ \cFidue^{-1} (D\times \T) \setminus\ 
\mathsmaller{\bigcup}
\check{\cal M}^i_k(\loge))\\
&\stackrel{\equ{chad}}\le&
(2\pi)^{n-1}\meas(D^\flat\times \T \setminus\ 
\mathsmaller{\bigcup}
\check{\cal M}^i_k(\loge))\\
&\stackrel{\equ{inthecourtk}}\le&
(2\pi)^{n-1}\hbfc\, \sqrt\suca
\meas(\hat D)\ \loge |\log \loge|\\
&\stackrel{\equ{orso}}<&
(2\pi)^{n-1} \hbfc\,  \,\Ro\, 
\meas(\hat D)\ \loge |\log \loge|\\
&\stackrel{\equ{sonnosonno2}}=&
\frac{\hbfc}{2\pi}\, 
\meas(\Rukt\times \T^n)\ \loge |\log \loge|\,,
\eeqano
which yields \equ{inthecourt2} since $\bfcast\ge \hbfc$. \qed

\rem
The measure estimate \equ{inthecourt2} holds in view of the covering property \equ{surge}, which takes care of the deformations near the boundaries.\\
The logarithmic correction is unavoidable and is related to the Lyapunov exponents  
of the hyperbolic equilibria issuing the separatrices of the secondary integrable systems at simple resonances.
\erem

\nl
The final result of this section deals with the size of the
domains $\Bik$, which depends on $k$ and actually  {\sl grows} with $k$. It is therefore important to control such a growth.

\begin{proposition}\label{lessing}
Assume that\footnote{Notice that, since $\g=2(\nu+n)$,  the hypothesis $\a=\e\K^\nu<1$ is implied by  the second condition in \equ{torino}.} $\a<1$. Then, there exists a constant   $\itcast=\itcast(n)>1$ such that 
\beq{diametro}
\diam \Bik\le \itcast |k|^{n-1} \,,\qquad \meas \Bik \le \itcast\,.
\eeq
\end{proposition}

\proof
For the purpose of this proof, we denote by `$c$' suitable (possibly different) constants greater than one and depending only on $n$. \\
Since $\a<1$, 
by the definition of $\Buik$ in \equ{aristakk},  by \eqref{alce} and the definition of $\Ro$ in \equ{cerbiatta},
we have, for every $0\le i\le 2N_k$,
\beqno  
\ts
 \diam \Buik
\leq c \big(\Ro +\diam\hat D\big)
\le c \big(\frac{\a}{|k|^2}+\diam \hat D \big)<c (1+\diam \hat D)\,.
\eeqno
By \eqref{venere}, \eqref{atlantide} and \equ{UU} it follows that 
$$\ts
|\proiezione_k^\perp \hAA^T \hat \act|=
|\AA^T{\rm U}
\binom{0}{\hat \act}|
\geq \frac{|\hat \act|}{\|\AA^{-1}\|\,\|{\rm U}^{-1}\|}
\geq \frac{|\hat \act|}{ c |k|^{n-1}}\,,
$$
and since by \equ{cerbiatta}  $\hat D\subseteq\{ \hat \act\in\R^{n-1}: \ |\proiezione_k^\perp \hAA^T \hat \act|<1\}$,
it follows that $\diam \Buik\le c |k|^{n-1}$, proving the first relation in \equ{diametro} in the case $0<i<2N_k$. \\
In the case $i=0,2N_k$, we need to estimate the Lipschitz constant of\footnote{ $\gstar$ is defined \equ{tordobis}.} $\ta_i$. The map $\giuno$ is linear and its gradient is given by $\hAA k/|k|^2$, thus, by \equ{atlantide} one gets 
\beqno\ts
|\partial_{\!{}_{\hat\act}}\giuno|=\big|\frac{\hAA k}{|k|^2}\big|
\le n\,.
\eeqno
By \equ{tess}, recalling that  $|\chk|\le 1$, the definitions in \equ{cerbiatta}, and
by Cauchy estimates\footnote{Compare, e.g., \cite{CC}.},  one sees that
\beq{nestore}\ts
|\gitre|_{4 \chr}<  \frac{2\e}{|k|^2}\frac\lalla\chr
  <\frac{\itcd \bfco}2\, \frac{\sqrt\e}{\K^{14n+5}}\,,
  \qquad
   |\partial_{\!{}_{\hat\act}}\gitre|_{3\chr}\le 
  \frac{2 \e \lalla}{\chr|k|^2 \ro^2}\le c\, \frac{\bfco}{\K^{14n+3}}<\frac14\,, 
\eeq
by taking $\KO$ big enough (recall that $\K\ge 6\KO$). Hence\footnote{$\Lip_{B}(g)$ denotes  absolute value of the Lipschitz constant of a function $g$ over a domain $B$.}, 
\beq{antiloco}
|\partial_{\!{}_{\hat\act}}\gstar|_{{}_{3\chr}}\le n+1\,,\qquad\quad  \Lip_{D_{{}_{3\chr}}}(\Istar)
\le n+2\,,
\eeq
and, choosing $\itcast$ suitably,  the first relation in  \equ{diametro} follows also in this case.\\
Let us check the second relation in \equ{diametro}. 
Since $\upphi_k^i$ in \equ{tordor} is symplectic, we have
\beqano
\meas \Bik &=&\frac1{(2\pi)^n}\meas (\Bik\times \T^n)=\frac1{(2\pi)^n}\meas \big(\upphi_k^i(\Bik\times \T^n)\big)\\
&\stackrel{\equ{tordor}}\le& \meas\big(\Re(\Ruk_{r_k})\big)\,.
\eeqano
Now, since $\Ruk\subseteq \DD$ and $r_k\le \a<1$, choosing $\itcast$ suitably, also the second relation in \equ{diametro} follows, and claim (i) has been proved.
\qed

\section{Twist at simple resonances}\label{thisistheend}

In this section -- which is the  heart of the paper -- we discuss  the main issue in singular KAM theory, namely, the twist of the integrable (rescaled) secular  Hamiltonians $\mathtt{h}^i_k$ 
in \equ{kant}
 near simple resonances and, in particular, in neighborhoods of secular separatrices, where the action become singular.

\nl
In general, it has to be expected that  
{\sl there are  points where the twist  of the  secular  Hamiltonians $\mathtt{h}^i_k$ 
 vanishes};  compare   Remark~\ref{bejazet} below. Furthermore, and more importantly, 
when
approaching separatrices, the evaluation of the twist becomes a {\sl singular perturbation problem}, where no standard tools can be applied and a new strategy is needed. \\
Our approach -- which  exploits in an essential way the fine analytic structure of the action functions described in Theorem~\ref{glicemiak} -- roughly speaking,  consists in constructing a suitable   differential operator with non--constant coefficients, which does not vanish on (a suitable regularization of) the Kolmogorov's twist determinant.  This will be enough to prove that {\sl the Liouville measure  of the set where the twist is smaller than a positive quantity $\fico$ may be bounded, uniformly in $k$,  by a power of $\fico$}.
This is the content of the Twist Theorem~\ref{zelensky} below, from which the proof of the results described in \S~\ref{omnia}  will follow easily.

\rem\label{bejazet}{\bf (Points where the twist vanishes)}   First, let us consider a region bounded by separatrices, i.e.,   (in the above setting) the case when $i$ is even and different from $0$ and $2N$. 
From \equ{arista}, \equ{LEGOk} and \equ{ciofecak} there follows that 
$\partial_{E} \act_1^i\to+\infty$ as $E$ approaches $E^i_\pm$. Thus, since  $\partial_{E} \act_1^i>0$ (always), 
$E\to\partial^2_{E} \act_1^i$ must have 
at least one zero in $(E^i_-,E^i_+)$. Since, by the chain rule,
\beq{sonnifera}
\partial^2_{\act_1}{\mathtt E}^{i} (\act_1)=
-\left.\frac{\partial^2_{E} \act_1^i}{(\partial_{E} \act_1^i)^3}\right|_{{\mathtt E}^{i} (\act_1)}\,,
\eeq
we see that $\partial^2_{\act_1}{\mathtt E}^{i} $ must vanish at some  points in the interval $(\acci^i,\bacci^i)$ defined in \equ{arista}.

\nl
Let us next consider the case $i$ odd, i.e., regions whose closure contains an elliptic point. Let us first consider the case $\lalla=0$, and let us denote 
$\bar\acci^i=a^i|_{\mu=0}$ and $\bar\bacci^i=b^i|_{\mu=0}$.
As above, by \equ{ciofecak}, the function
 $E\to\partial^2_{E}\bar \act_1^i$
 tends to $+\infty$ when $E\to \bar E^i_+$. Thus,  by \eqref{sonnifera},
  $\partial^2_{\act_1}\bar{\mathtt E}^{i}(\act_1)$
  is negative when $\act_1$ is close to $\bar\bacci^{i}$.
  Now, $\bar{\mathtt E}^{i}(\act_1)$
 is analytic at  $\act_1=\bar\acci^{i}=0$,
 and,  evaluating the 
 Birkhoff normal form of $p_1^2+\GO(q_1)$ at order 4 close to the elliptic point $(p_1,q_1)=(0,\bar\sa_i)$, one sees that 
 $$
 \bar{\mathtt E}^{i}(\act_1):=\o_0  \act_1+
 \frac12 c  \act_1^2+ O(\act_1^3)\,,
 \quad {\rm with}\ \ 
 \o_0=\sqrt{2 d_2}\,,\ \ 
 c=\frac14\left(\frac{d_4}{d_2}-\frac{5 d_3^2}{3 d_2^2} \right)\,,
 $$
where $d_j$ are the $j$-th order derivatives of the reference 
potential $\GO$ evaluated at the minimum $\bar\sa_i$.
Thus,  $\partial^2_{\act_1}\bar{\mathtt E}^{i}(0)>0$
whenever the condition 
\beq{aurora}
\updelta:=3 d_2 d_4-5 d_3^2>0\,,\qquad d_j:=(\partial_{q_1}^j\GO)(\bar\sa_i)\,,
\eeq
is satisfied,  in which case
 $\partial^2_{\act_1}\bar{\mathtt E}^{i}$ must vanish at some point  in $\big(0,\bar\bacci^{i}\big)$. 
By \equ{tikitaka}, 
$\hzk=\htk$, so that
$
\mathtt{h}^i_k|_{\lalla=0}=
\bar{\mathtt E}^{i}_k(\act_1)+\htk(\hat\act)$,
which implies 
$$
\det\partial^2_{\act}\mathtt{h}^i_k (\act)|_{\lalla=0}
=
\partial^2_{\act_1}\bar{\mathtt E}^{i}_k (\act_1)\cdot
   \det\partial^2_{\hat \act} \htk(\hat\act)\,.
$$
Thus, by continuity, for   $\lalla$ small enough  
 it follows that 
the Hessian matrix $\partial^2_{\act}\mathtt{h}^i (\act)$
is singular at some point.
 
 \nl
Condition \equ{aurora} is easily satisfied. For example, if $\GO(\sa)=\cos \sa - \frac18 \cos(2\sa)$, one finds that $\updelta=3/2$, so that, in this very simple cases, inside the (unique) region enclosed by the main separatrices, there are points where the twist vanishes. 
However, this is not the case if the potential is close enough  to a cosine, compare Proposition~\ref{kalevala}.
\erem

\subsection*{Twist Theorem near simple resonances (statement)}

To state the Twist Theorem we need to introduce two parameters
($\upxi>0$,  $\mtt \geq 1$) 
which measure the non--degeneracy (in a suitable sense to be specified below) of the energy as function of actions in the inner regions $0<i<2N_k$. This requires some preparation.

\subsubsection*{Non--degenerate functions and theirs sub--levels}
First, let us  recall a standard quantitative definition of  non--degenerate functions. 
\dfn{troll} 
Given $\xi>0$, an open set $A\subseteq \real$ and   $f\in C^m(A,\real)$, we say that $f$   
is $\xi$--non--degenerate at order $m\ge 1$ on $A$
(or, in short, $(\xi,m)$--non--degenerate), 
if
\begin{equation}\label{troll1}
\inf_{x\in A}\max_{ 1\leq j\leq m}| 
f^{(j)}(x)|\geq \xi \,.
\end{equation}
\edfn
An important property of non--degenerate functions is that one can easily estimate the measure of their sub--levels:
\begin{lemma}\label{nessuno} 
 Let $f$ be a $(\x,m)$--non--degenerate function on a bound\-ed interval $(a,b)$
and let\footnote{$\|f\|_{C^{m+1}(a,b)}:=\max_{0\leq j\leq m+1}\sup_{(a,b)}|f^{(j)}|$.} $M:=\|f\|_{C^{m+1}(a,b)}$.
Then, there exist a constants  $c_m>1$ depending only on $m$ such that, for all $\fico>0$, one has  
$$
{\rm meas}  \{ x\in (a,b):|f(x)|\leq \fico  \} 
\leq   \frac{c_m}{\xi^{1/m}}
\big(\ts\frac{M}{\xi} (b-a)+1\big) \, \fico^{1/m}\,.
$$
\end{lemma}
The proof of this  lemma can be found, e.g.,  in  \cite[Lemma B.1]{E}; compare, also,  \cite{pyartli}.

\subsubsection*{Non--degeneracy of the rescaled reference potentials for $\noruno{k}\le \Nf$}
Consider a general Hamiltonian \equ{pasqua} in standard form, recall Definition~\ref{piggy},  recall \equ{arista}, and define also, for $0\le\loge\leq \blogemax$ (defined in \equ{sinistro}),
\beq{wagner}
\bar\acci^i:=\acci^i|_{\mu=0}\,,\quad 
\bar\bacci^i:=\bacci^i|_{\mu=0}\,,\quad
\bar\acci^i_\loge:=\acci^i_\loge|_{\mu=0}\,,\quad 
\bar\bacci^i_\loge:=\bacci^i_\loge|_{\mu=0}\,,\quad\  
\forall\, 0\le i\le 2N_k\,.
\eeq
In the following, we shall explicitly indicate the dependence upon the reference  potential 
 $\GO$ and write, e.g,  $\bar \act_{1,\GO}^{i}$, $\bar{\mathtt E}^{i}_{\GO}$, $\bar\acci^{i}_{\GO}$, $\bar\bacci^{i}_{\GO}$ for 
$\bar \act_{1}^{i}$, $\bar{\mathtt E}^{i}$, $\bar\acci^{i}$, $\bar\bacci^{i}$, respectively.

\dfn{mrwhite} Given $\Hpend$ in standard form with reference potential $\GO$, 
 we denote by
 \begin{equation}\label{nutella}
{\mathtt F}^{i}_{\GO}(x):=
(\partial^2_{I_1} \bar{\mathtt E}^{i}_{\GO})
\big(\bar\acci^{i}_{\GO}+(\bar\bacci^{i}_{\GO}-\bar\acci^{i}_{\GO})x
\big)\,,\quad\forall\ x\in (0,1)\,,\qquad (0<i<2N_k)\,,
\end{equation}
the `normalized second derivative of the energy function within separatrices'.
\edfn
These functions  satisfy a remarkable rescaling property:

\begin{lemma}\label{genzano}
If ${\mathtt F}^{i}_{\GO}$ is as in Definition~\ref{mrwhite}, then, for any $ \l>0$, one has 
${\mathtt F}^{i}_{\GO}={\mathtt F}^{i}_{\l\GO}$.
\end{lemma} 
\proof Indeed, from the definition of  actions, there follows easily that
 \begin{equation}\label{simpson}
\bar \act_{1,\l\GO}^{i}(E)=\sqrt\l \bar \act_{1,\GO}^{i}(E/\l)\,,
\qquad  
\bar{\mathtt E}^{i}_{\l\GO}(\act_1)=\l\bar{\mathtt E}^{i}_{\GO}(\act_1/\sqrt\l)\,,\qquad \forall \l>0\,.
\end{equation}
Indeed, considering the case $i=2N_k$ (the other cases being similar),
one has
 $$
 \bar \act_{1,\l\GO}^{2N_k}(E)
\stackrel{\eqref{bolsena}}{=}
\frac{1}{2\pi}\int_0^{2\pi}\sqrt{E-\l\GO(x)}dx
=
\frac{\sqrt\l}{2\pi}\int_0^{2\pi}\sqrt{\frac{E}{\l}-\GO(x)}dx
=\sqrt\l \bar \act_{1,\GO}^{i}(E/\l)
\,,
 $$
which proves the first equality in \equ{simpson}, which, in turns, implies immediately the second inequality.
From \equ{simpson}, then , follows that
 \begin{equation}\label{simpson2}
 \bar\acci^{i}_{\l\GO}=\sqrt\l \bar\acci^{i}_{\GO}\,,
 \qquad
  \bar\bacci^{i}_{\l\GO}=\sqrt\l \bar\bacci^{i}_{\GO}\,,
\end{equation}
and the claim follows at once from 
 \eqref{simpson} and \eqref{simpson2}. \qed 
Let us go back to the Hamiltonians in standard form  $\Hsharp$ of Theorem~\ref{garrincha}, and let 
us prove that  the functions ${\mathtt F}^{i}_{\GO}$ -- and hence $\bar{\mathtt E}^{i}_{\l\GO}$ -- with $\GO$ as in \equ{paranoia}, are  $(\xi,m)$--non--degenerate.

\begin{lemma}\label{goblin}
For every $0< i< 2N_k$,
the function ${\mathtt F}^{i}_{\GO}$ defined in \equ{nutella} is 
 $(\xi,m)$--non--degenerate for some $\xi,m>0$.
 \end{lemma}
\proof
We consider only the case $i$ odd, the even case being similar.
 Deriving  \eqref{sonnifera} we get, for $\mu=0$,
 \begin{equation}\label{soonnifera2}
 \partial^3_{I_1} \bar{\mathtt E}^{i}(\bar \act_1^{i}(E))
 =-\frac{\partial^3_{E}\bar \act_1^{i}(E)}
 {\big(\partial_E\bar \act_1^{i}(E)\big)^4}
 +3\frac{\big(\partial^2_{E}\bar \act_1^{i}(E)\big)^2}
 {\big(\partial_E\bar \act_1^{i}(E)\big)^5}
 \,.
 \end{equation}
 By \eqref{LEGOk}--\eqref{vana} 
 (which   hold also  for $\bar \act_1^{i}$, corresponding to $\lalla=0$),
we have that the dominant term in  \eqref{soonnifera2}
 as $z:=(\bar E^{i}_+-E)/\suca\to 0^+$ has the form
 $-1/({\mathtt c}^3 z^2 \ln^4 z)$
 with ${\mathtt c}:=\psi^{i}_+(0)|_{\lalla=0}.$
Then,
$$
 \lim_{E\to (\bar E^{i}_+)^-}  
 \big|\partial^3_{I_1} \bar{\mathtt E}^{i}(\bar \act_1^{i}(E))\big|=
 \lim_{\act_1\to  (\bar\bacci^{i})^-}
  \big|\partial^3_{I_1} \bar{\mathtt E}^{i}(\act_1)\big|=
 +\infty\,.
$$
By \eqref{nutella} we obtain
 \begin{equation}\label{palato}
\lim_{x\to 1^-} |\partial_x {\mathtt F}^{i}_{\GO}(x)|=+\infty\,.
\end{equation}
Moreover $\partial_x {\mathtt F}^{i}_{\GO}(x)$ is analytic in a neighborhood of $x=0$
(recall in particular \eqref{lamponek}).
Assume now by contradiction  that \eqref{troll1} does not hold, namely that
there exists a sequence $x_m\in (0,1)$ such that 
$$
|\partial_x^j{\mathtt F}^{i}_{\GO}(x_m)|<1/m\,, \qquad \forall\, 1\leq j\leq m\,.
$$
By \eqref{palato}, up to a subsequence, $x_m$ converges to some $\bar x\in [0,1)$
such that $\partial_x^j{\mathtt F}^{i}_{\GO}(\bar x)=0$ for every $j\geq 1$. By analyticity we would have that
${\mathtt F}^{i}_{\GO}$ is constant on $[0,1)$ leading to a contradiction with \eqref{palato}.
\eproof

\nl
This lemma allows us  to introduce
{\sl uniform} non--degeneracy
parameters $\upxi>0$ and  $\mtt \geq 1$
for the function ${\mathtt F}^{i}_{\GO}$
in \eqref{nutella}
associated to 
 the 
reference potentials
$\GO\stackrel{\equ{paranoia}}= {\textstyle \frac{2\e}{|k|^2}}\, \fproj f$,
for $k\in\gen$, $\noruno{k}<\Nf$ and
$0<i<2N_k$.
 Indeed,  by Lemma~\ref{genzano},
\begin{equation}\label{prato}
{\mathtt F}^{i}_{\GO}={\mathtt F}^{i}_{\frac{2\e}{|k|^2}\, \fproj f}
={\mathtt F}^{i}_{\fproj f}
\,,
\end{equation}
 and, by  \eqref{P2+},
 every
 potential
$\fproj f$ is $\b$--Morse.
By the above Lemma \ref{goblin},  
every function in \eqref{prato}
is $(\xi,m)$--non--degenerate  
for some $\xi,m>0$.
We therefore can define  uniform $\e$--independent non--degeneracy  parameters $\upxi$, $\mathtt m$ by setting:
\dfn{dracula}
{\sl Let ${\mathtt F}^{i}_{\fproj f}$ be as in Definition~\ref{mrwhite} with rescaled reference potential $\GO=\fproj f$.
We define
$\upxi>0$ and  $\mtt \geq 1$
 to be, respectively, the largest and  smallest number such that all the functions
 ${\mathtt F}^{i}_{\fproj f}$, for $0<i<2N_k$, 
 $k\in\gen$ with $\noruno{k}\le \Nf$,  
 are $(\upxi,\mtt)$--non--degenerate (Definition~\ref{troll})}.
\edfn

\subsubsection*{The Twist Theorem}

Let Assumptions~\ref{assunta} and Definitions~\ref{assunto} hold, let  $\upkappa$ be as in \equ{kappa},  let $\upxi,\mtt$ be as in Definition~\ref{dracula}, let 
$\Bik$ be as in \equ{tordobis},  let $\mathtt{h}^i_k$ be as in \equ{kant}, and define
\beq{balla}
\dalla:= {|k|}^{-2n}\,.
\eeq
Then, the following result holds.

 \begin{theorem}\label{zelensky} There exists a constant 
 $\cteot=\cteot( n,\upkappa,\upxi,\mtt)>1$
 such that, for $\KO\ge \cteot$,   $k\in\genKO$,  $0\leq i\leq 2N_k$,  and  $0<\fico< \dalla/2^5$,   one has:
 \beq{gibboso}
\meas\big(\big\{ \, \act\in \Bik:\  \big| \det \partial^2_{\act}\mathtt{h}^i_k(\act) \big|\leq \fico\big\}\big)
\le 
\cteot  (|k|^{2n} \fico)^\lippo
\meas \Bik
 \,,\qquad \lippo:=\min\{\ts \frac1{9n^4}\,,\,\frac{1}{\mtt}\}\,.
\eeq
\end{theorem}
Theorem~\ref{zelensky} will be proven in several steps:

\nl
{\bf Step 1:}    {\sl Preliminaries}\\
{\bf (a)} Explicit  expressions for the twist  matrix in the inner case ($0<i<2N_k$) are  given; 
\\
{\bf (b)} analogous formulae are given for the outer case ($i=0,2N_k$), but, due to the presence of the translation  $\Istarin$ in \equ{tordissimo}, the measure estimate is expressed in terms of the domains
$\Buik$  rather than the domains $\Bik$ (recall that such domains  differ in the outer case; compare \equ{tordobis});
\\
{\bf (c)} uniform estimates on the   sub--matrix $\partial^2_{\hat \act} \hzk$ of order $(n-1)$, depending only on  the `trivial actions' $\hat\act$, 
 are given.

\nl
{\bf Step 2:} 
{\sl Coverings of the phase space into regions close to separatrices and far from separatrices}\\
 This is a necessary step, since 
the analysis will be  non
perturbative near separatrices, while in   regions 
away from separatrices,
the analysis will be partly
perturbative (and significantly simpler).

\nl
{\bf Step 3${}^*$:} {\sl Non--degeneracy of the twist function in neighborhoods of separatrices} \\
In such regions perturbative arguments do not hold, and, in  particular the energy function 
${\mathtt E}^{i}$ is singular at the boundary (corresponding to separatrices) and its derivatives diverge as the boundary  is approached. Furthermore,  ${\mathtt E}^{i}$ and 
$\bar{\mathtt E}^{i}={\mathtt E}^{i}|_{\lalla=0}$ have singularities in different points.
Exploiting the singularity structure
described in Theorem~\ref{glicemiak}, we will prove that a suitable regularization of the twist determinant is a
non--degenerate function allowing to control the measure of its sub--levels. 
This is the core of the proof.

\nl
{\bf Step 4:} {\sl The Twist Theorem  in  neighborhoods of  separatrices} \\
By the previous step, measure estimates in regions close to separatrices follow easily, yielding the proof of the Twist Theorem in this case.

\nl
{\bf Step 5:}  {\sl The Twist Theorem far from separatrices in the inner case}\\
 It is here (in particular, in the low mode case $\noruno{k}\le \Nf$) that the non--degeneracy condition involving the parameters $\upxi$ and $\mtt$, is needed.

\nl
{\bf Step 6:}  
{\sl Uniform twist  in outer regions far from separatrices}\\
In such regions there is uniform twist; the proof
rests on  a simple  argument based on Jensen's inequality.

\nl
{\bf Step 7:} {\sl Conclusion of the proof of  the Twist Theorem}

\subsection*{Proof of  the Twist Theorem}

Fix $k\in\genKO$,    $0\leq i\leq 2N_k$,  and   $\fico>0$. 

\nl
{\sl Throughout the proof the $(n-1)$ dimensional domain $\hat D$ (defined in \equ{cerbiatta}) will be kept fixed and often the variables $\hat \act$ will not be indicated explicitly. Also, the label
$k$ will  usually be omitted it  in the notation, as well as 
 the suffix $i$ (when this does not lead to confusion).}

\nl
{\bf Step 1} {\sl Preliminaries}\\
 {\bf (a)}  We give 
 the analytic expression of the  twist determinant  inside separatrices, i.e., for $0<i<2N$. Recall that in this case, by \eqref{kant} and  \eqref{tordobis}, one has
$\mathtt{h}^i={\mathtt E}^{i}+\hzk$ and  $\Bi=\Buik$. Then\footnote{Observe that if ${\rm S}=({\rm s}_{ij})_{i,j\le n}$ is an $(n\times n)$ matrix and $\hat {\rm S}$ denotes the  $(n-1)\times(n-1)$ sub--matrix 
$({\rm s}_{ij})_{i,j\ge 2}$, and ${\rm S}_0$ denotes the matrix obtained by ${\rm S}$ 
replacing the entry ${\rm s}_{11}$ with 0, then $\det {\rm S} = {\rm s}_{11}\cdot \det \hat {\rm S} + \det {\rm S}_0$.
}
\allowdisplaybreaks
\begin{align}\label{parata}
 \det\partial^2_{\act}\mathtt{h}^i
&=
\det\big( 
\partial^2_{\act}{\mathtt E}^i 
+
\partial^2_{\act} \hzk 
\big)
=
  \det\left(
  \begin{array}{cc}
   \partial^2_{\act_1}{\mathtt E}^i  
  &
    \quad\partial_{\hat \act}^T(\partial_{\act_1}{\mathtt E}^i  ) 
    \\ 
 \partial_{\hat \act}(\partial_{\act_1}{\mathtt E}^i  )    
    &
    \partial^2_{\hat \act}{\mathtt E}^i  +\partial^2_{\hat \act}\hzk
  \\ 
  \end{array}
  \right)
\\
&=
   (\partial^2_{\act_1}{\mathtt E}^i )
   \cdot \det ( \partial^2_{\hat \act}{\mathtt E}^i  +\partial^2_{\hat \act} \hzk )
   +
    \det\left(
  \begin{array}{cc}
    0
  &
    \quad\partial_{\hat \act}^T(\partial_{\act_1}{\mathtt E}^i  ) 
    \\ 
 \partial_{\hat \act}(\partial_{\act_1}{\mathtt E}^i  )    
    &
 \partial^2_{\hat \act}{\mathtt E}^i +\partial^2_{\hat \act} \hzk    
  \\ 
  \end{array}
  \right) \,.
  \nonumber
\end{align}
{\bf (b)} We now consider the case outside the outer  separatrices, i.e.,  $i=0,2N$.

\nl
The Hamiltonian $\mathtt{h}^i(\act)$, in this case, is given by\footnote{Recall  \eqref{kant},  \eqref{tordobis} and \eqref{arista}.}
$\mathtt{h}^i(\act)=\widetilde{\mathtt E}^i(\act)+\hzk(\hat\act)$
for $\act\in B^i=\Istarin \big(\Buik\big)$. 
Recalling \equ{tikitaka} and \equ{tess} we note that
in the evaluation of the Hessian of $\mathtt{h}$
 involves the {\sl non--small linear term} $\frac{(\hAA k)\cdot \hat \act}{|k|^2}$, a fact that  complicates analytic expressions. However, 
 such complications may be avoided, using the following trick.
\\
Let us introduce new action variables $\bfI$, defined by the relation 
$\act={\rm U}\bfI=j_{\!\giuno}(\bfI)$, where ${\rm U}$ is defined in \equ{Fiu}.
 Then, we observe that, defining
\beq{batman}
h^i(\bfI) :=\bfEs(\bfI)+\hzk(\hat\bfI)\,,\qquad   \bfEs:=
 {\mathtt E}^i_k\circ  j_{\!{}_{-\gitre}}\,,\ 
\eeq
one has that
\beq{piripicchio}
\Istar= j_{\!{}_{-\gitre}}\circ {\rm U}^{-1}\,,\qquad
\mathtt{h}^i({\rm U}\bfI)=h^i (\bfI)\,,\qquad (\forall\ \bfI\in {\rm U}^{-1} B^i)\,.
\eeq
Now,   since $\det {\rm U}=1$,
$$
\det \big[\partial^2_\bfI h^i (\bfI)\big]\stackrel{\equ{piripicchio}}{=}\det \big[\partial^2_\bfI \big(\mathtt{h}^i({\rm U}\bfI)\big)\big]
=\det \big[  {\rm U}^T  \partial^2_\act \mathtt{h}^i(\act) \  {\rm U}   \big]= \det  \big[ \partial^2_\act \mathtt{h}^i(\act) \big]\,.
$$
Thus,
\beq{questasi}
(\det \partial^2_\act \mathtt{h}^i )\circ {\rm U} 
=\det \partial^2_\bfI h^i\,.
\eeq
Recalling \eqref{batman} we then obtain
\beqa{paratafora}
\det \partial^2_{\bfI} h^i &=&
\det\big( 
\partial^2_{\bfI}\bfEs
+
\partial^2_{\bfI} \hzk 
\big)
=
  \det\left(
  \begin{array}{cc}
   \partial^2_{\bfI_1}\bfEs  
  &
    \quad\partial_{\hat \bfI}^T(\partial_{\bfI_1}\bfEs ) 
    \\ 
 \partial_{\hat \bfI}(\partial_{\bfI_1}\bfEs )    
    &
    \partial^2_{\hat \bfI}\bfEs +\partial^2_{\hat \bfI}\hzk
  \\ 
  \end{array}
  \right)
  \\
&=&
   (\partial^2_{\bfI_1}\bfEs)
   \cdot \det ( \partial^2_{\hat \bfI}\bfEs +\partial^2_{\hat \bfI} \hzk )
   +
    \det\left(
  \begin{array}{cc}
    0
  &
    \quad\partial_{\hat \bfI}^T(\partial_{\bfI_1}\bfEs ) 
    \\ 
 \partial_{\hat \bfI}(\partial_{\bfI_1}\bfEs )    
    &
 \partial^2_{\hat \bfI}\bfEs+\partial^2_{\hat \bfI} \hzk    
  \\ 
  \end{array}
  \right) 
  \nonumber
\eeqa
and, by the chain rule, 
\begin{eqnarray}
(\partial^2_{\bfI_1} \bfEs)\circ  j_{\!{}_{\gitre}}
&=&
\partial^2_{\act_1} {\mathtt E}^i 
\,,
\nonumber
\\
(\partial^2_{\bfI_1 \hat \bfI}\bfEs)\circ  j_{\!{}_{\gitre}}
&=&
\partial^2_{\act_1 \hat \act} {\mathtt E}^i -
\partial^2_{\act_1} {\mathtt E}^i\partial_{\hat\act} 
\gitre
=:{\rm {\rm \hat v}}
\,,
\label{daitarn4}
\\
(\partial_{\hat \bfI\hat \bfI}^2\bfEs)\circ  j_{\!{}_{\gitre}}
&=&
\partial^2_{\hat \act} {\mathtt E}^i +
\partial_{\act_1}^2 {\mathtt E}^i  \partial_{\hat \act}^T\gitre
\partial_{\hat \act}\gitre
-
\partial_{\act_1} {\mathtt E}^i  \partial_{\hat \act}^2\gitre
-\partial_{\hat\act}^T\partial_{\act_1} {\mathtt E}^i 
 \partial_{\hat \act}\gitre
-
\partial_{\hat \act}^T\gitre
\partial_{\hat\act}\partial_{\act_1} {\mathtt E}^i 
=:{\rm \hat M}
\,.
\nonumber
\end{eqnarray}
Recalling that by \equ{piripicchio}
 $\Istar= j_{\!{}_{-\gitre}}\circ {\rm U}^{-1}$,
 by \eqref{questasi}, \eqref{paratafora} and \eqref{daitarn4} 
we get
\begin{equation}\label{enumaelish}
\dds^i:=
(\det   \partial^2_\act \mathtt{h}^i )\circ 
\Istarin=
\partial^2_{\act_1} {\mathtt E}^i
   \cdot \det ( {\rm \hat M} +\partial^2_{\hat \act} \hzk )
   +
    \det\left(
  \begin{array}{cc}
    0
  &
    {\rm \hat v}^T
    \\ 
 {\rm \hat v}  
    &
 {\rm \hat M}+\partial^2_{\hat \act} \hzk    
  \\ 
  \end{array}
  \right)
 \,.
\end{equation}
Finally,  since the map $\Istar:B^i_k\to\Buik$ is volume preserving, it is 
\beq{pizza}
\meas B^i_k=\meas \Buik\,,\qquad i=0,2N\,,
\eeq
so that one obtains the following  

\begin{lemma}\label{depechemode} Let $i=0,2N$ and  $\dds^i$  as in \equ{enumaelish}.
Then, 
\beqno
\meas\big(\big\{ \, \act\in B^i  \ {\rm s.t.}\  
\big| \det \partial^2_{\act}\mathtt{h}^i(\act) \big|
\leq \fico
\big\}\big)
=
\meas\big(\big\{ \, \act\in \Buik  \ {\rm s.t.}\  \big| \dds^i(\act) \big|\leq \fico \big\}\big)
 \,. 
\eeqno
\end{lemma}

\nl

\nl
{\bf (c)} Here we prove the following uniform bound on the Hessian sub--matrix $\partial^2_{\hat \act} \hzk$. Recall the definition of $\dalla$ in \equ{balla}

\begin{lemma}\label{laromascaja} There exists $\itct=\itct(n)>1$ such that  
if $\K\ge \itct$  the following estimates on the sub--matrix $\partial^2_{\hat \act} \hzk$ hold:
\begin{equation}\label{lovebuzz}
\sup_{\hat D_{\ro}}|\partial^2_{\hat \act} \hzk|
\leq 
2n^5+1\,,
\qquad\quad
\inf_{\hat D_{\ro}\cap\R^{n-1}}
\det \partial^2_{\hat \act} \hzk
\geq  \dalla\,.
\end{equation}
\end{lemma}

\proof 
By \equ{barocco},
\beq{perpieta}
(\AA^T{\rm U}) \act=\act_1 k + \hAA^T\hat \act-  {\textstyle \frac{(\hAA k)\cdot \hat \act}{|k|^2}k}=
\act_1 k + \hAA^T\hat \act- {\textstyle \frac{  \hAA^T \hat \act \cdot k}{|k|^2}k}
=\act_1 k + \proiezione_k^\perp \hAA^T \hat \act\,.
\eeq	
Recalling the definition of $\htk$ in \equ{tikitaka}, we have
\beq{pallinadicacca}\ts
  \partial^2_\act\big(\act_1^2+ \htk(\hat \act)\big)
= 
  \partial^2_\act\big( \act_1^2 + \frac{ \normadue \proiezione^\perp_k \hAA^T \hat\act\normadue ^2}{|k|^2}
 \big)
 \eqby{perpieta}
  \frac{\partial^2_\act\normadue \AA^T{\rm U} I \normadue^2}{|k|^2}  
  =    \frac{  2(\AA^T{\rm U})^T  \AA^T{\rm U}}{|k|^2}
\eeq
and
\beq{puritani}
  | \partial^2_{\hat\act}\htk|
\leq
 2| k|^{-2}  |\AA|^2|{\rm U}|^2  
 \stackrel{\eqref{atlantide},\eqref{UU}}\leq
 2 n^5\,.
 \eeq
 Using that $|k|\le \K/6$,
 by \equ{athlone}, \equ{cerbiatta},  \equ{tikitaka}
and  Cauchy estimates
we get, for a suitable $c'=c'(n)$,
\beq{impepata}
\sup_{\hat D_{\ro}}
| \partial^2_{\hat\act}(\hzk-\htk)|
\leq
\frac{c'}{\K^{14n+2}}\,.
\eeq
By  \equ{puritani} and \equ{impepata},
taking $\K$ large enough (depending only on $n$),
we get 
 the first estimate in \eqref{lovebuzz}.

\nl
Let us prove the second estimate in \equ{lovebuzz}.
Observe that 
\beqno
\ts
2 \det  \partial^2_{\hat\act} \htk 
=
\det  \partial^2_\act\big(\act_1^2+ \htk   \big)
\stackrel{\equ{pallinadicacca}} =
\frac{2^n}{| k|^{2n}} \det \big((\AA^T{\rm U})^T  \AA^T{\rm U}	\big)
\stackrel{(\ref{atlantide},\ref{Fiu})}= \frac{2^n}{| k|^{2n}}\ge \frac{4}{| k|^{2n}}
= 4 \dalla\,,
\eeqno
and that
\beq{oloturia}
|(\partial^2_{\hat\act}\htk)^{-1}|\leq
  |(\partial^2_\act( \act_1^2 +\htk))^{-1}|
 \stackrel{\eqref{pallinadicacca}} \leq 
  {\ts \frac{|k|^{2}}2}  |\AA^{-1}|^2|{\rm U}^{-1}|^2
  \stackrel{\eqref{atlantide},\eqref{UU}}\leq
  {\ts \frac{1}2}
  n^5(n-1)^{n-1} |k|^{2n}  
\,.
\eeq
Then, by \equ{oloturia}, \equ{impepata},  using $|k|\le \K/6$, we get, for a suitable constant $c''=c''(n)$, 
\beq{follia}
|(\partial^2_{\hat\act}\htk)^{-1}|\cdot
 | \partial^2_{\hat\act}(\hzk-\htk)|\leq 
\frac{ c''}{\K^{12n+2}}\,,
\eeq

\nl
We now need an  elementary result on perturbation of positive--definite matrices,  whose proof is given  in   
\Appendix:

\begin{lemma}\label{tonino}
 Let $P,Q$ be $d\times d$ positive--definite  matrices   and assume that 
$\l:=|P^{-1}| |Q|$ is strictly smaller than 1.  Then $\det (P+Q)\geq (1-\l)^d \det P$.
 In particular, if $\l\leq (2d)^{-1}$, then $\det (P+Q)\geq  (\det P)/2$.
\end{lemma}
Then, observe that since $\partial^2_\act\big(\act_1^2+ \htk   \big)$ is   positive--definite  (by \equ{pallinadicacca}), so
is  $\partial^2_{\hat\act} \htk$. Therefore, 
since $ \partial^2_{\hat \act} \hzk=\partial^2_{\hat\act} \htk +\partial^2_{\hat\act}(\hzk-\htk) $,
in view of \equ{follia}, taking $\K\ge \itct$ for a suitable $\itct=\itct(n)>1$, Lemma~\ref{tonino} implies also the second estimate in \equ{lovebuzz} and the proof of Lemma~\ref{laromascaja} is complete. \qed

\nl
{\bf Step 2} Here we define suitable {\sl  coverings of  the sets $\Buik$} defined in \equ{arista},
 \equ{cerbiatta}. Such coverings are
made up of  sets  corresponding  to  zones  close to the separatrices and zones away from them.

\nl
Recall the definitions  given  in
\equ{arista}, \equ{sinistro}
and \equ{wagner}. 
For any
$\logeo\in(0,1/\bfcast)$,  define the following  subsets of\footnote{The constant $\bfcast$ satisfies
\equ{sgallettata}. Recall that for $i$ odd $\acci^i_0\equiv0$ (see \equ{accidenti}).
The number $\logeo$   will be fixed in Proposition~\ref{enigmista} below.}   
$\Buik$:
\allowdisplaybreaks
\begin{align}
&\left\{ 
\begin{array}{l}\Buc^i(\logeo)
:=
\{ \act: 
\bacci^i_{\logeo}(\hat\act)
< \act_1<
\bacci^i(\hat\act)\,,
 \hat \act\in\hat D  \}\,,
 \nonumber
 \\ 
\Bua^i(\logeo)
:= \vallo^i\times \hat D\,,\quad 
\vallo^i:=(0,\bar\bacci^i_{\logeo/2})\,,
\end{array}
\right.
& i \ \ {\rm odd}\,;
\\
\label{tetto}
&\left\{ 
\begin{array}{l}
\Buc^i(\logeo):=
\{ \act: 
\acci^i(\hat\act)
< \act_1<
\acci^i_{\logeo}(\hat\act)\,,
\hat \act\in\hat D
\}
\\
\phantom{\Buc^i(\logeo):=\ }
\cup\{\act: \bacci^i_{\logeo}(\hat\act)
< \act_1<
\bacci^i(\hat\act)\,,
 \hat \act\in\hat D  \}\,,
 \\ 
\Bua^i(\logeo)
:= \vallo^i\times \hat D\,,\quad 
\vallo^i:=(\bar\acci^i_{\logeo/2},\bar\bacci^i_{\logeo/2})\,,
\end{array}
\right.
&i \ \ {\rm even}\,, \ \ i\neq 0, 2N_k\,;
\\
\nonumber
&\left\{ 
\begin{array}{l}\Buc^i(\logeo)
:=
\{ \act: 
\acci^i(\hat\act)
< \act_1<
\acci^i_{\logeo}(\hat\act)
\,,
 \hat \act\in\hat D  \}\,,
 \\ 
\Bua^i(\logeo)
:= \vallo^i\times \hat D\,,\quad 
\vallo^i:=(\bar\acci^i_{\logeo/2},\bacci^i(\hat\act))\,,
\end{array}
\right.
&  i= 0, 2N_k\,.
\end{align}
Then,  one has:
\begin{lemma}\label{mediatore}
Let $0\le i\le 2N$   and assume that\footnote{$\bfcast$ appears in Theorem~\ref{garrincha}, while
$\bfc\le \bfcast$ appears in Theorem~\ref{glicemiak}. Recall   \equ{sgallettata}.} 
\begin{equation}\label{BWV1013}
\logeo< 1/\bfcast\,,\qquad\quad  
\lalla\leq \logeo^2/2^8\bfc^4\,.
\end{equation}        
Then,
$\Buik=\Buc^i(\logeo)\cup\Bua^i(\logeo)$.
\end{lemma}
\proof
We give a detailed proof  in the case  \equ{tetto} ($i$ odd), as there is no extra difficulty in extending the proof to the other cases.
For ease of notation  in this proof we omit the suffix $i$.\\
Since the functions $E\to \bar\act_1(E)$
and  $E\to \act_1(E,\hat\act)$
are positive and strictly increasing (see \eqref{vana}),
the functions
$\loge\to \bar\bacci_\loge$ and 
$\loge\to\bacci_\loge(\hat\act)$
are positive and strictly decreasing.
We claim that
\begin{equation}\label{activia}
\bacci_{\logeo}(\hat \act)
<
\bar\bacci_{\logeo/2}
<
\bacci_{\logeo/4}(\hat \act)
\,,
\qquad
\forall\, \hat\act\in\hat D
\,.
\end{equation}
From such relations the claim follows:  the fact that $\Bua$ is a subset of $\Bu$
follows from the second inequality in 
\eqref{activia}, and then  the equality
$\Bu=\Buc\cup\Bua$
follows from the first inequality in 
\eqref{activia}.\\
Let us prove in detail the first inequality in \eqref{activia} 
(the second one being  analogous).
By \eqref{tetto} and \eqref{sinistro} we 
have\footnote{Recall  that $\bacci_{\logemax}(\hat \act)
=I_1(E_-(\hat \act),\hat\act)=\acci_0(\hat\act)=0$.}
$$
\bacci_{\logeo}(\hat \act)
=\suca \int_{\logeo}^{\logemax}
\partial_E I_1(E_+(\hat \act)-\suca\z ,\hat \act)\, d\z
=
\suca \int_{\logeo+\logesharp}^{\logemax+\logesharp}
\partial_E I_1(\bar E_+-\suca\z ,\hat \act)\, d\z
\,,
$$
where
$\logesharp:=(\bar E_+ -E_+(\hat \act))/\suca$.
Analogously,
$$
\bar\bacci_{\logeo/2}
=\suca \int_{\logeo/2}^{\blogemax}
\partial_E \bar I_1(\bar E_+-\suca\z)\, d\z\,,
$$
where $\blogemax$ was defined in \eqref{sinistro}. 
Note that, by \eqref{october} and \equ{sgallettata}, $|\logesharp|\le 
3\upkappa^3\lalla\leq \bfc\lalla$, and that 
by  \equ{sinistro}, \equ{latooscuro}    we have that
 $\logeo\leq1/\bfcast\le\min\{\logemax/4,\blogemax/8\}$.
Then again by \eqref{latooscuro}, \eqref{sinistro},
\eqref{destro} and \eqref{BWV1013} 
  we get that, for every $\hat\act\in\hat D$,
$$
\ts{
\frac{\logeo}{2},\ \logeo+\logesharp,\ 
\logemax+\logesharp,\ \blogemax
\ \in\ 
\left(\frac{\logeo}8, \blogemax+\frac{1}{\bfc}\right).}
$$
We write
\begin{eqnarray*}
\frac{
\bar\bacci_{\logeo/2}
-
\bacci_{\logeo}(\hat \act)}{\suca}
&=& 
 \int_{\logeo/2}^{\logeo+\logesharp}
\partial_E \bar\act_1(\bar E_+-\suca\z)\, d\z
+
\int_{\logemax+\logesharp}^{\blogemax}
\partial_E \bar\act_1(\bar E_+-\suca\z)\, d\z
\\
&&\!\!\!\!\!\!\!
+
\int_{\logeo+\logesharp}^{\logemax+\logesharp}
\big(\partial_E \bar\act_1(\bar E_+-\suca\z)-
\partial_E \act_1(\bar E_+-\suca\z ,\hat \act)\big)\, d\z\,,
\end{eqnarray*}
observing that, 
for every $\z$ in the three integration intervals
(and for every $\hat\act\in\hat D$),
the quantity $\bar E_+-\suca\z $
belongs to the 
set\footnote{Recall \eqref{autunno2}. Note that \eqref{BWV1013}
implies \equ{bassoraTH} with $\loge=\logeo/8$.} $\mathcal E_{\logeo/8}$.
Then, by \eqref{vana}, \eqref{rosettaTH}
and  \eqref{BWV1013} we get, for every $\hat\act\in\hat D$,
\begin{eqnarray*}
\frac{
\bar\bacci_{\logeo/2}
-
\bacci_{\logeo}(\hat \act)}{\suca}
\geq
\frac{\logeo-2|\logesharp|}{2\bfc\sqrt\suca}
-\bfc^2\frac{|\ln\frac{\logeo}{8}|}{\sqrt\suca}
\big(|\blogemax-\logemax|+|\logesharp|\big)
-8\bfc^2\lalla\frac{\logemax}{\logeo\sqrt\suca}
\\
\stackrel{\eqref{destro},\eqref{latooscuro}}\geq
 \frac{1}{2\bfc\sqrt\suca}\left(
\logeo-2\bfc\lalla-
2^4\bfc^4 \lalla |\ln\ts\frac{\logeo}{8}|
-2^6\bfc^3\lalla/\logeo
\right)
\geq \frac{\logeo}{4\bfc\sqrt\suca}>0
\,. \qedeq	
\end{eqnarray*}

\nl
{\bf Step 3${}^*$} {\sl Non--degeneracy of the twist function in neighborhoods of separatrices}

\nl
Here we show that (a suitable regularization of) the twist determinant $ \det\partial^2_{\act}\mathtt{h}^i$ in \equ{parata} is a non--degenerate function in the sense of Definition~\ref{troll} in {\sl suitable 
 neighborhoods of separatrices}. \\
 Actually, it will be convenient to study the twist directly as a function of the energy, for values    $E=E_\mp^{i}(\hat \act)\pm \suca \z$
 close to critical separatrix values $E_\mp^{i}$. We therefore define:
\beq{mattoni0}
\ddi_{\!{}_\mp}(\z,\hat \act):=
\det\big[ 
\partial^2_{\act} {\mathtt E}
   \big( \act_1(E_\mp^{i}(\hat \act)\pm \suca \z, \hat \act),\hat \act\big)
   +\partial^2_{\act} \hzk (\hat \act)\big] \,.
\eeq

\nl
The study of the twist determinant \equ{mattoni0} will be based on the analytic properties described in Theorem~\ref{glicemiak}. In particular, the properties that we shall use are the same in the plus and the minus case. Hence, we shall consider only the plus  case  and {\sl consider, henceforth,  $\ddi:=\ddi_{\!{}_+}$}.

\nl
The precise statement on the non--degeneracy of $\z\to \ddi(\z,\hat \act)$ (see Proposition~\ref{enigmista} below)  needs  some preparation.

\nl
First of all, we introduce a suitable
`regularization' function $\vulva=\vulva(\z,\hat\act)$
\beq{vulva}
\vulva(\z,\hat\act):=\z \cdot \big(\sqrt\suca\partial_E \act_1 (E_+(\hat \act)-\suca\z ,\hat \act)\big)^3\,,
\eeq
and define  {\sl the regularized twist determinant} $\ddb$ by setting\footnote{Recall \equ{lovebuzz}.}
\beq{mattoni}
\begin{array}{l}\ddb 
=\dt/
\det \partial^2_{\hat \act} \hzk\,,\quad {\rm with}\quad 
       \dt(\z ,\hat \act):= 
\vulva(\z,\hat\act)^n\cdot 
\ddi(\z ,\hat \act)
\,. \phantom{\dst \int}
\end{array}
\eeq
The functions appearing in Theorem~\ref{glicemiak}, as well as the functions in \equ{vulva} and \equ{mattoni} belong to the following ring 
of functions $\cF$.

\begin{definition} 
We denote by $\cF$ the set of functions of the form 
\beq{brendl}
f(\z,\hat\act)=\z^h \sum_{j=0}^{\ell} u_j(\z,\hat\act) \log^j \z\,,
\eeq
where $h,\ell\in\integer$ with $\ell\ge 0$ and the $u_j$ are real analytic functions on a (complex) neighborhood of\footnote{Recall that the domain $\hat D$ is defined in \equ{cerbiatta}, but, essentially plays no r\^ole.} 
$\{z=0\}\times \hat D\subset \complex^n$. \\
We shall also use the following notation: given two functions  $f_i\in \cF$ we say that $f=f_1\oplus f_2$ if there exists two functions $u_i$ real analytic  on a neighborhood of $\{z=0\}\times \hat D$ such that\footnote{E.g., $f$ in \equ{brendl} can be written as $ z^h \big(\bigoplus_{j=0}^{\ell}  \log^j \z\big)$.} $f=u_1f_1+u_2 g_2$.

\nl
We say that  $f(\z,\hat\act)=\cOr(h,\ell)$
if $f\in \cF$ as in 
\equ{brendl} and there exists $\varrho>0$ such that 
$$
\vvvert f\vvvert_\varrho:=\sup_{0\le j\le \ell} \sup_{\sopra{\{z\in \complex: |\z|<\varrho\}}{\hat \act\in \hat D}} |u_j|<+\io\,.
$$
\end{definition} 

\rem (i) The functions 
$(\z,\hat\act)\to f(\z,\hat\act)= \act_1^{i}\big(E_\mp^{i}(\hat \act)\pm \suca \z, \,\hat \act\big)$ in \equ{LEGOk} of  Theorem~\ref{glicemiak} belongs to $\cF$ and, by \equ{pappagallok}, 
$$
\vvvert f\vvvert_{1/\bfc} \le \bfc\sqrt\suca\,;
$$
furthermore, the `algebraic structure' of such function $f$ is given by
$$
f=\sqrt\suca(1\oplus \z\log\z) \,.
$$
(ii) The following elementary properties (which, in particular, show that $\cF$ is a ring) will be often used:
\beqno
\left\{
\begin{array}{l}
 \cOr(h,p)\cdot \cOr(k,q)=\cOr(h+k,p+q)\,,\\
 (\cOr(h,p))^j=\cOr(jh,jp)\,,\\
  \cOr(h,p)+\cOr(k,q)= \cOr(\min\{h,k\},\max\{p,q\})  \,.
\end{array}
\right.  
\eeqno
\erem 

\nl
Finally, define the 
 following linear differential operators:
\beqno
\mathcal L:=L^{3\bn} (\partial_\z \cdot L^{3\bn})^{\bn}\,,\quad {\rm where:}\quad L:=\z \partial_\z\,,\quad \bn:= n-1\,.
\eeqno
Notice that  $\mathcal L$ is a linear differential operator of order $\mmu:=3\bn^2+4\bn=3n^2-2n-1
\geq 7$ and 
there exist  
suitable polynomials $a_j(\z)$
such that\footnote{Actually, $\mathcal L =
\sum_{j=\bn+1}^{\mmu}
a_j z^{j-\bn} 
\partial_\z^j 
$, with $a_j\in \natural$.
For example, if $n=2$,  $\mmu=7$ and $\mathcal L$ is given by:\\
$\phantom{aa...}\mathcal L=\z^6 \partial_\z^7
 + 18 \z^5 
 \partial_\z^6
  + 98 \z^4 
  \partial_\z^5
+ 
184 \z^3
\partial_\z^4
 + 100 \z^2 
\partial_\z^3
 + 8 \z 
 \partial_\z^2
 \,.$}
\begin{equation}\label{ciaone}
\mathcal L 
=
\sum_{j=1}^{\mmu}
a_j(\z)
\partial_\z^j \,.
\end{equation}

\begin{proposition}\label{enigmista} 
There exists\footnote{The constant $\itct$ has been introduced in Lemma~\ref{laromascaja}.} $\bfcd=\bfcd(n,\kappa)>\itct$ such that if $\K\ge \bfcd$, then the following holds.

\nl
{\rm (i)} One has  
 \begin{equation}\label{caciottone}
\mathcal L [\ddb]
=
\bn!^{3\bn+1}(3\bn)!\ 
\balla^{3\bn} 
+\cOr(1,3\bn+1)\,,
\end{equation}
where\footnote{Recall \equ{ciofecak}.}
$$\balla(\hat\act):=-\suca^{-1/2}\psi_+(0,\hat\act) \,,\qquad 
\varrho:=1/\bfc\,,  
$$
and
\begin{equation}\label{basilisco}
1/\bfc \le \inf_{\hat D}|\balla|\leq
\sup_{\hat D}|\balla|\le \bfc\,.
\end{equation}
{\rm (ii)}  There exist suitable positive constants 
$\xish=\xish(n,\upkappa)< 1$ and\footnote{The constant $\bfcast$ has been introduced in Theorem~\ref{garrincha}.}
$\logeo=\logeo(n,\upkappa)<1/\bfcast$,
 such that,  for  $\hat \act\in\hat D$, the function  $\z\to \ddb(\z,\hat \act)$ defined in \equ{mattoni} is $\xish$--non--degenerate 
at order $\mmu=3n^2-2n-1$
on the interval $(0,\logeo)$.
\end{proposition}

\nl
To prove this proposition we need a couple of preparatory lemmata.

\begin{notation}
In  the rest of this section, it is understood  that in an expansion $f=z^h \big(\bigoplus_{j=0}^{\ell}  \log^j \z\big)$, one has  $\vvvert f\vvvert_\varrho\le c$ for a suitable constant $c=c(n,\kappa)$; furthermore, $\cO$ stands for $\cOr$ with $\varrho=1/\bfc$.
\end{notation}

\nl
We shall  consider in detail only  the {\sl inner odd case $0<i<2N$},   since the other cases do not present any new difficulties; for ease of notation, we do not indicate explicitly the labels $k$ and $i$.

\begin{lemma}\label{3-0} If $\vulva$ is as in \equ{vulva},
$\tilde\act=\tilde\act(\z,\hat\act)
:=(\bacci_\z(\hat\act),\hat\act)$
and 
$\lella$   is as in \eqref{pappagallok}, one has ($2\leq i,j\leq n$)
\begin{eqnarray}
\vulva
&=&
z\big(\balla\log \z+(1\oplus\z\log\z) \big)^3
\nonumber\\
&=&
\balla^3  \z\log^3 \z+\cO(1,2) +\cO(2,3)=\cO(1,3)
\,,
\nonumber
\\
\vulva\cdot\partial^2_{\act_1} {\mathtt E}  
|_{I=\tilde\act}
&=&
\balla+\z(1\oplus\log\z) 
=\balla+\cO(1,1)
\,,\label{fausti}
\\
\vulva\cdot\partial^2_{\act_1 \hat \act_i} {\mathtt E}  
|_{I=\tilde\act}
&=&
\lella (1\oplus \z \log \z \oplus \z \log^2 \z )=\lella \,\cO(0,2)
\,,
\nonumber
\\
\vulva\cdot\partial^2_{\hat \act_i \hat \act_j} {\mathtt E}  
|_{I=\tilde\act}
&=&
\lella (1\oplus \z \log \z \oplus \z \log^2 \z 
\oplus \z^2\log^3 \z )=\lella \,\cO(0,3)
\,.\nonumber
\end{eqnarray}
\end{lemma}
\proof
By the chain rule, one has (writing $I_1$ in place of $I_1^i$)
\begin{eqnarray}
&&\partial_{\act_1} {\mathtt E}^i 
=
\frac{1}{\partial_E \act_1 }\,,\qquad
\partial_{\hat \act} {\mathtt E}^i 
=
-\frac{\partial_{\hat \act} \act_1 }{\partial_E \act_1 }\,,
\qquad
\partial^2_{\act_1} {\mathtt E}^i 
=
-\frac{\partial^2_{E} \act_1 }{(\partial_E \act_1 )^3}\,,
\nonumber
\\
&&\partial^2_{\act_1 \hat \act} {\mathtt E}^i 
=
\frac{\partial^2_{E} \act_1  \partial_{\hat \act}\act_1 }{(\partial_E \act_1 )^3}-
\frac{\partial^2_{E\hat \act} \act_1  }{(\partial_E \act_1 )^2}
\,,
\label{daitarn3}
\\
&&\partial^2_{\hat \act} {\mathtt E}^i 
=
-\frac{\partial^2_{\hat \act} \act_1 }{\partial_E \act_1 }
+\frac{\partial_{\hat \act}^T\act_1\ \partial_{\hat \act}(\partial_{E} \act_1)
+ \partial_{\hat \act}^T(\partial_{E} \act_1)\ \partial_{\hat \act}\act_1 }
{(\partial_E \act_1 )^2}
-\frac{\partial^2_{E} \act_1 \  \partial_{\hat \act}^T \act_1\  \partial_{\hat \act}\act_1}{(\partial_E \act_1 )^3}
\,,
\nonumber
\end{eqnarray}
where the derivatives of ${\mathtt E}^i$ and  $\act_1=\act_1^i$ are evaluated in  $\big(\act_1^i (E,\hat \act),\hat \act\big)$
and $(E,\hat \act),$ respectively.
Now,
by \equ{daitarn3} and  \eqref{LEGOk}, we have
\begin{eqnarray}\label{gorgo}
&&\sqrt\suca\partial_E \act_1
=\balla\log \z+(1\oplus\z\log\z) 
=1\oplus \log\z \,,\quad
\partial_{\hat \act_i}\act_1=\lella (1\oplus \z\log\z )\,,
\nonumber
\\
&&\suca^{3/2}\partial^2_E \act_1=
-\balla\z^{-1}+
 (\log\z \oplus  \z^{-1})=\log\z \oplus  \z^{-1}\,,
\\
&&\suca\partial^2_{E\hat \act_i} \act_1=\lella (1\oplus \log\z )\,,\quad
\partial^2_{\hat \act_i\hat \act_j}\act_1=\lella 
\ro^{-1}(1\oplus \z\log\z )
\stackrel{\eqref{cimabue}}=
\lella 
\suca^{-1/2}(1\oplus \z\log\z )\,.
\nonumber
\end{eqnarray}
Finally, by \eqref{gorgo}, \eqref{daitarn3}, \eqref{LEGOk}
and \eqref{pappagallok}
we get 
\begin{eqnarray*}
\sqrt\suca\partial_E \act_1   
&=&
\balla\log \z+(1\oplus\z\log\z) 
=
1\oplus  \log \z 
\nonumber 
\\
\vulva\cdot\partial^2_{\act_1} {\mathtt E}  
&=&
-\suca^{3/2}\z \partial^2_E \act_1 
=\balla+\z(1\oplus\log\z) 
=
1\oplus \z \log \z 
\,,
\nonumber
\\
\vulva\cdot\partial^2_{\act_1 \hat \act_i} {\mathtt E}  
&=&
\suca^{3/2}\z (\partial^2_E \act_1  \partial_{\hat \act_i}\act_1 -
 \partial^2_{E\hat \act_i} \act_1  \partial_E \act_1)
=
\lella (1\oplus \z \log \z \oplus \z \log^2 \z )
\,,
\nonumber
\\
\vulva\cdot\partial^2_{\hat \act_i \hat \act_j} {\mathtt E}  
&=&
\suca^{3/2}\z\big(
- (\partial_E \act_1 )^2 \partial^2_{\hat \act_i\hat \act_j} \act_1 
+2  \partial_E \act_1 \partial_{\hat \act_i}  \act_1\ \partial^2_{\hat \act_j E} \act_1 
- \partial^2_E \act_1 \  \partial_{\hat \act_i} \act_1\  \partial_{\hat \act_j}\act_1\big)
\nonumber
\\
&=&
\lella (1\oplus \z \log \z \oplus \z \log^2 \z 
\oplus \z^2\log^3 \z )
\,. \qedeq
\end{eqnarray*}

\begin{lemma}\label{robin} 
One has
\beq{takeiteasy}
\ddb=\ddb(\z ,\hat \act)=
\balla^{3\bn}  \z^\bn\log^{3\bn} \z
+\cO(\bn+1,3\bn+1)
+\cO(0,3\bn-1)\,,
\quad\bn:=n-1\,.
\eeq
Furthermore, there exists $\bfcd=\bfcd(n,\kappa)>\itct$ such that
if $\K\ge \bfcd$ and 
$\logeo\le 1/\bfcd$,  one has:
\begin{equation}\label{terencehill}
|\ddi(\z ,\hat \act)|\ge \dalla\ |\ddb(\z ,\hat \act)|
\,,
\qquad
\forall\, 0<\z\leq \logeo\,,\  \hat \act\in\hat D\,.
\end{equation}
\end{lemma}
\proof
Recalling \equ{parata} we split $\dt$ in \equ{mattoni}  in two terms.
The first term is
$$
\vulva^n (\partial^2_{\act_1}{\mathtt E})
   \ \det (\partial^2_{\hat \act} \hzk +  \partial^2_{\hat \act}{\mathtt E})
   \stackrel{\eqref{fausti}}=
   \big(\balla+\cO(1,1)\big)\,
   \vulva^{\bn}\det (\partial^2_{\hat \act} \hzk +  \partial^2_{\hat \act}{\mathtt E})\,,
$$
and, by  \eqref{fausti}, we have that
\begin{eqnarray*}
&&\vulva^{\bn} 
   \ \det (\partial^2_{\hat \act} \hzk +  \partial^2_{\hat \act}{\mathtt E})
   \\
   =&&
   \vulva^{\bn} 
   \ \det \partial^2_{\hat \act} \hzk + 
   \sum_{j=1}^{\bn} \vulva^{\bn-j}	
   \lella^j \big(1\oplus \z \log \z \oplus \z \log^2 \z 
\oplus \z^2\log^3 \z\big)^j
\\
=&&
\big(\balla^{3\bn}  \z^\bn\log^{3\bn} \z
+\cO(\bn,3\bn-1)+\cO(\bn+1,3\bn)
\big)
\det \partial^2_{\hat \act} \hzk
\\
&&
+\lella\big(\cO(0,3\bn-3)+\cO(\bn,3\bn-1)+\cO(\bn+1,3\bn) \big)
\\
=&&
\big(\balla^{3\bn}  \z^\bn\log^{3\bn} \z
+\cO(\bn+1,3\bn)
+\cO(0,3\bn-1)
\big)
\det \partial^2_{\hat \act} \hzk\,,
\end{eqnarray*}
where in the last line we used (recall 
\eqref{pappagallok}, \eqref{cerbiatta}, 
 \eqref{lovebuzz}) 
 
\begin{equation}\label{dalla}
\lella\leq |k|^{-2n}=\dalla\leq \inf_{\hat D} \det \partial^2_{\hat \act} \hzk\,.
\end{equation}

\nl
The second term is,
\begin{eqnarray*}
&&\vulva^n
    \det\left(
  \begin{array}{cc}
    0
  &
    \quad\partial_{\hat \act}^T(\partial_{\act_1}{\mathtt E}  ) 
    \\ 
 \partial_{\hat \act}(\partial_{\act_1}{\mathtt E}  )    
    &
 \partial^2_{\hat \act}{\mathtt E} +\partial^2_{\hat \act} \hzk    
  \\ 
  \end{array}
  \right)
  =\lella
  \det\left(
  \begin{array}{cc}
    0
  &
    w^T
    \\ 
w 
    &
{\rm  N}    
  \\ 
  \end{array}
  \right)
\\
  &=&\lella ^2 (1\oplus \z  \log \z  \oplus \z  \log^2 \z )^2
(1\oplus \z  \log \z  \oplus \z  \log^2 \z \oplus 
\z\log^3 \z )^{n-2}
\\
&=&  
\lella ^2 \cO(0,3n-4)
+\cO(n,3n-2)
\,,
\end{eqnarray*}
where $w$ is a  $\bn$--dimensional vector and 
${\rm N}$ is an $(\bn\times\bn)$ matrix satisfying  by \eqref{fausti}
 $$
w_i= \cO(0,0)+\cO(1,2) \,,\ \ 2\leq i\leq n\,,\qquad
{\rm N}_{ij}= \cO(0,0)+\cO(1,3) \,,\ \ 2\leq i,j\leq n \,.
$$
Thus, the second term has the form
$
\lella^2\big(\cO(\bn+1,3\bn+1)
+\cO(0,3\bn-1) \big)$.
Summing up the two terms and using
\eqref{dalla} we get \eqref{takeiteasy}. 

\nl
By the first line in \eqref{fausti}, we see that, taking $\bfcd$ big enough, one has
\beqno
|\vulva(\z,\hat\act)|\leq 1\,,   
\qquad
\forall\, 0<\z\leq \logeo\,,\  \hat \act\in\hat D\,.
\eeqno
Thus, by the definitions in \equ{mattoni0}, \equ{mattoni} and by \equ{lovebuzz}, one obtains \equ{terencehill}.
\qed

\nl
Before giving the proof of Proposition~\ref{enigmista}, we need one more lemma. Define
$$
\mathcal L_{m,k}:= L^{k} (\partial_\z \cdot L^{k})^{m}\,.
$$

\begin{lemma}
Let $0\le \ell<k\le \mmu$, $0\le m, q\le \mmu$, and $f_1=\cOr(0,\ell)$, $f_2=\cOr(m+1,q)$. Then
\beq{haus}
\mathcal L_{m,k}[\z^m\ln^k\z+f_1+f_2]
=(m!)^{k+1}k!+f_3\,,
\eeq
where, for a suitable constant $c$, which depends only on $n$, one has 
$$
f_3=\cO_{\!\varrho/2}(1,\max\{k-1,q\})\,,\quad {\rm and}\quad 
\vvvert f_3\vvvert_{\varrho/2}\le c \max\{\vvvert f_1\vvvert_\varrho\,, \vvvert f_2\vvvert_\varrho\}\,.
$$
\end{lemma}
\proof 
Observing that $L\z^m=m\z^m$, $L\log^{\ell+1}\z=(\ell+1) \log^\ell\z$, one easily checks that, for any $0\le m,\ell\le \mmu$, one has
\beqno
\begin{array}{ll}
L [\cOr(m,\ell)]=\cO_{\!{\frac34\varrho}}(m,\ell)\,, \phantom{\dst \int} \\
 L\cOr(0,\ell+1)=\cO_{\!{\frac34\varrho}}(1,\ell+1)+\cO_{\!{\frac34\varrho}}(0,\ell)\,,\phantom{\dst \int}\\
L^{\ell+1} [\cOr(0,\ell)]=\cO_{\!{\frac34\varrho}}(1,\ell)\,,\phantom{\dst \int}
\end{array}
\eeqno
where the norm $\vvvert\cdot\vvvert_{\!{\frac34\varrho}}$ of the functions in the right hand sides are bounded by $c'=c'(n)$ times the norm
$\vvvert\cdot\vvvert_{\varrho}$ of the functions in the left hand sides\footnote{The algebraic relations are just calculus, while  the estimates follow easily by (iterated) use of Cauchy estimates. 
}. Analogously, from the above relations, it follows that, for any $0\le \ell<k\le \mmu$ and $0\le q\le \mmu$, one has
\begin{eqnarray*}
&&\mathcal L_{m,k}  [\z^m
\log^{k}  \z]=(m!)^{k+1}(k)!+\cO_{\!{\frac12\varrho}}(1,k-1)\,, \phantom{\dst \int} 
\\
&&\mathcal L_{m,k} [\cO_{\!{\frac34\varrho}}(0,\ell)]=\cO_{\!{\frac12\varrho}}(1,\ell)\,, \phantom{\dst \int} \\
&&
\mathcal L_{m,k}[\cO_{\!{\frac34\varrho}}(m+1,q)]=\cO_{\!{\frac12\varrho}}(1,q)\,,\phantom{\dst \int} 
\end{eqnarray*}
where, the norm $\vvvert\cdot\vvvert_{\!{\frac12\varrho}}$ of the functions in the right hand sides are bounded by $c=c(n)>c'$ times the norm
$\vvvert\cdot\vvvert_{\frac34\varrho}$ of the functions in square brackets in the left hand side. From such relations, the claim of the lemma follows easily.
\qed

\proof {\bf of Proposition~\ref{enigmista}} 
The estimates in \eqref{basilisco} follow trivially from \eqref{ciofecak} and  \eqref{pappagallok}.

\nl
To check \equ{caciottone}, observe that   $\mathcal L=\mathcal L_{\bn,3\bn}$, and use   \eqref{takeiteasy}  in Lemma~\ref{robin} and \equ{haus} (with 
$m=\bn,$ $k=3\bn$, $\ell=3\bn-1$, and $q=3\bn+1$).

\nl
It remains to prove claim (ii).
By \eqref{caciottone}, 
and \eqref{basilisco}, we  see that for $\logeo<1/\bfcast$ small enough one  has:
\beqno
\frac1{\bfc^{3\bn}}\le 
\inf_{0<\z\leq \logeo}
\inf_{\hat \act\in\hat D}\big|\mathcal L [\ddb]\big| \stackrel{\equ{ciaone}}{\le} c'' \max_{1\le j\le \mmu}|\partial_\z^j\ddb|
\,,
\eeqno
where $c''=c''(n)$. 
Thus, for  $\hat \act\in\hat D$,  $\z\to \ddb(\z,\hat \act)$ is $\xish$--non--degenerate 
at order $\mmu=3n^2-2n-1$
on the interval $(0,\logeo)$ with $\xish= (c'' \bfc^{3\bn})^{-1}$.
\qed

\nl
{\bf Step 4:} {\sl The Twist Theorem  in  neighborhoods of  separatrices} 

\nl
We can now state and prove the Twist Theorem in neighborhoods of separatrices.

\begin{proposition}\label{zelenskyC} Let $k\in\genKO$, $0\le i\le 2N$, $\fico>0$, and $\logeo$ as in Proposition~\ref{enigmista}--(ii). 
Then,  there exist a  positive constant
${\bfcu}={\bfcu}(n,\upkappa)\ge \bfc$   such that,  if $\K\geq {\bfcu}$, then   
\beq{gibbosobis}
\meas\big(\big\{ \act\in \Buc^i(\logeo):\ |\det \partial^2_{\act}\mathtt{h}^i (\act) |\le \fico\big\}\big)
\le \bfcu \, (|k|^{2 n}\fico)^{1/9 n^4}\, \meas \Buik 
\,.
\eeq
\end{proposition}

\nl
Before the proof, which  will be based on two lemmata, we introduce the following

\begin{notation}\label{weiss}
Given two non negative  functions $f$ and $g$ we say that $f\lessdot g$ if there exists a constant $c=c(n,\kappa)\ge 1$, depending only on $n$ and $\kappa$, such that $f\le c g$. Similarly, given a function $f$ and a non negative function $g$,
we say that $f=O(g)$ if there exists a constant $c=c(n,\kappa)\ge 1$, such that $|f|\le c g$.
\end{notation}

\begin{lemma} There exists a constant $\bfct=\bfct(n,\kappa)>1$ such that, for every $\hat \act\in \hat D$ and 
 $\fico>0$, one has\footnote{$\mmu=3n^2-2n-1$
is defined in Proposition~\ref{enigmista}.}
\begin{equation}\label{prepuzio}
\meas \{ \z\in(0,\logeo]\  {\rm s.t.}\ 
|\ddb(\z,\hat \act)|\leq \fico  \}
\le \bfct \fico^\gu\,,
\qquad 
\ts\gu:=\frac{1}{\mmu(\mmu+3)}
\,.
\end{equation}
\end{lemma}
\proof  
If $z_0\leq 2\fico^\gu$ estimate \equ{prepuzio} is obvious.
Consider the case $z_0> 2\fico^\gu$.
Let 
$$
\logeuno:= \fico^\gu<z_0/2\,.
$$
By  \eqref{takeiteasy}, \eqref{lovebuzz} and \eqref{basilisco} we have that\footnote{Denoting,
as usual, the $\logeuno/2$--complex--neighborhood of the real interval
$[\logeuno,\logeo]$ by $[\logeuno,\logeo]_{\logeuno/2}$.}
$$
\sup_{[\logeuno,\logeo]_{\logeuno/2}}\sup_{ \hat D} |\ddb(\z,\hat \act)|\lessdot
1+|\log^{3n-4} \logeuno|\lessdot 1/\logeuno\,.
$$
By Cauchy estimates
\begin{equation}\label{loosing2}
\sup_{\logeuno\leq \z\leq \logeo}
\sup_{\hat \act\in\hat D}
\max_{ 1\leq j\leq \mmu+1}|\partial_\z^j \ddb(\z,\hat \act)|
\leq c_\flat /\logeuno^{\mmu+2}=:M\,,
\end{equation}
for a suitable $c_\flat\geq1$, depending only on $n,\upkappa$.
Now  we  want to apply
Lemma \ref{nessuno}
  with
$$
f= \ddb\,,\quad
 m= \mmu\,,\quad
  a = \logeuno\,,\quad
   b= \logeo\,,\quad
 \xi= \xish\,,
 \quad
 M \ {\rm as \ in\ } \eqref{loosing2}\,.
$$
 Then,  we get
\begin{equation}\label{asgard}
{\rm meas}  \{ \z\in(\logeuno,\logeo):
|\ddb(\z,\hat \act)|\leq \fico  \} 
\lessdot   \fico^\gu\,,\quad 
\forall\ \hat \act\in \hat D\,.
\end{equation}
Since the interval $(0,\logeuno)$, has length
$\logeuno=\fico^\gu$, from \equ{asgard} we obtain
 the measure estimate  \eqref{prepuzio}.
 \eproof
Now, recalling that $\dalla=|k|^{-2n}$ (see  \equ{balla}),  we have:
\begin{lemma}\label{paradise}
There exists $\bfcq=\bfcq(n,\upkappa)\ge \max\{\bfcd,\bfct\}$  such 
that
for $k\in\genKO$, $i$ odd and $\fico>0$, 
$$
\meas\big(
\big\{ \act\in \Buc^i:\ 
|\det \partial^2_{\act}
\mathtt{h}^i (\act) |
\le \fico
\big\}
\big)
\le \bfcq
\, 
\sqrt\suca\,
(\fico/\dalla)^{{\frac{1}{9 n^4}}}\,
\meas \hat D\,.
$$
\end{lemma}
\proof
Let $\mathcal Z_\fico (\hat \act):= \big( \{ \z\in(0,\logeo]:\  |\ddi(\z ,\hat \act)|\leq \fico  \}$. 
By
\eqref{terencehill}
and 
 \eqref{prepuzio} we get
\begin{equation}\label{missionbell0}
{\rm m}_\fico ={\rm m}_\fico (\hat \act):=\meas(\mathcal Z_\fico (\hat \act))
\le \bfct  (\fico/\dalla)^\gu\,,
 \qquad \forall\  \ \hat \act\in \hat D\,.
 \end{equation}
Note that  since $\logeo\le 1/2$ (see \eqref{BWV1013}), by definition
\begin{equation}\label{primosale}
{\rm m}_\fico\leq \logeo \leq1/2\,.
\end{equation}
Recalling \equ{arista},
we define,  for $\hat \act\in\hat D$ and $\fico>0$,
\beqno
\mathcal \act_\fico(\hat \act):=
\big\{ \act_1\in
\big[\bacci_{\logeo}(\hat\act),
\bacci(\hat \act)\big):\ \ 
\big|\det\big[\partial^2_\act \big(\hzk (\hat \act)+ 
{\mathtt E} (\act)\big)
  \big] \big|\leq \fico  \big\}\,.
\eeqno
We have that
\beq{animals}
 \mathcal \act_\fico(\hat \act) 
 =
  \bacci_{\mathcal Z_\fico(\hat \act)}(\hat \act)
  :=\big\{\act_1 =\bacci_z(\hat \act) :
  \z\in \mathcal Z_\fico(\hat \act)
  \big\}
   \,,
\eeq
since by definition of $\mathcal Z_\fico$, \eqref{mattoni0} and \eqref{tetto}
$\ddi(\z ,\hat \act)=
\det\big[\partial^2_\act \hzk (\hat \act)+ 
\partial^2_\act{\mathtt E}
   \big(\bacci_\z(\hat\act),\hat \act\big) \big]$.
For every $\hat \act\in\hat D$ and $\fico>0$,  making the change
of variable $I_1=\bacci_\z(\hat\act)$, and noticing that 
$\partial_\z \bacci_z(\hat\act)=-\suca \partial_E \act_1(E_+(\hat \act)-\suca\z ,\hat \act)$, 
we get
\begin{eqnarray*}
{\rm meas}(\mathcal \act_\fico(\hat \act))&=&\int_{\mathcal \act_\fico(\hat \act)}d\act_1
\stackrel{\equ{animals}}=
 \int_{\bacci_{\mathcal Z_\fico(\hat \act)}}d\act_1
=\int_{\mathcal Z_\fico(\hat \act)} |\partial_\z 
\bacci_\z(\hat\act)|d\z
\\
&\stackrel{\eqref{fausti}}\lessdot&
\sqrt\suca
\int_{\mathcal Z_\fico (\hat \act)}  |\log\z|d\z\,.
\end{eqnarray*}
Moreover, recalling \eqref{primosale}, 
$$
\int_{\mathcal Z_\fico (\hat \act)}  |\log\z|d\z
\le
\int_0^{{\rm m}_\fico}
 |\log\z|d\z
+\int_{\mathcal Z_\fico (\hat \act)\cap({\rm m}_\fico ,\logeo]} |\log\z|d\z
\le 2 {\rm m}_\fico |\log{\rm m}_\fico |
\,.
$$
Thus, by  \eqref{primosale}, using ${\frac{1}{9 n^4}}<\gu$
(see \equ{prepuzio} and recall the definition of $\mmu=3n^2-2n-1$  in Proposition~\ref{enigmista}), we get
$$
{\rm meas}(\mathcal \act_\fico(\hat \act))
\lessdot
\sqrt\suca\, {\rm m}_\fico |\log{\rm m}_\fico |
\lessdot
\sqrt\suca\, {\rm m}_\fico^{1/(9n^4\, \gu)}
\,.
$$
By \eqref{missionbell0},
${\rm meas}(\mathcal \act_\fico(\hat \act))
\lessdot
\sqrt\suca\,
(\fico/\dalla)^{1/9n^4}$, 
for every $\hat \act\in\hat D$ and $\fico>0$.
By Fubini's Theorem,   the claim follows. 
\eproof

\proof {\bf of Proposition~\ref{zelenskyC}}
By  \eqref{arista} we get
\begin{eqnarray}\label{natalino}
\meas \Buik&=&\int_{\hat D}\bacci(\hat\act)d\hat\act
=\int_{\hat D}
\act_1\big(E_+(\hat \act),\hat \act\big)
d\hat\act
=\int_{\hat D}d\hat\act
\int_{E_-(\hat \act)}^{E_+(\hat \act)}
\partial_E\act_1\big(E_+(\hat \act),\hat \act\big)dE
\nonumber
\\
&\stackrel{\eqref{vana}}\geq&
\frac{E_+(\hat \act)-E_-(\hat \act)}{\bfc\sqrt\suca}\meas \hat D
\stackrel{\eqref{sinistro}}\geq 
\frac{\sqrt\suca}{2\upkappa\bfc}\meas \hat D\,.
\end{eqnarray}
Lemma~\ref{paradise} and 
\eqref{natalino} imply at once \eqref{gibbosobis} if one takes $\bfcu\ge \bfcq$ big enough. The proof of Proposition~\ref{zelenskyC} is complete for $i$ odd. The changes for the inner case with $i$ even  are straightforward .

\nl
Let us  indicate the changes one needs to do in order to prove  the outer case  $i=0,2N$.
Recalling \equ{nestore}, \equ{cerbiatta}, \equ{athlone},  by 
 Cauchy estimates, we have (recall Notation~\ref{weiss})
\begin{equation}\label{appannato}
|\partial_{\hat\act}\gitre|_{\hat D,3 \ro}\lessdot
1/\K^{14n+2}\,,
\qquad
|\partial^2_{\hat\act}\gitre|_{\hat D,3 \ro}\lessdot
1/(\sqrt\e\, \K^{\frac{37}2 n+2})\,.
\end{equation}
Note that the term $\partial^2_{\hat\act}\gitre$ in \eqref{appannato} has a `big' estimate, containing
 a $\sqrt\e$ at the denominator. However
 this does not cause any problem, since,   by the first lines in \eqref{fausti}, \equ{daitarn3} and \equ{gorgo}  one has\footnote{The regularizing term $\vulva$
 defined in Lemma \ref{3-0}.}
 $$
 \vulva\, \partial_{\act_1} {\mathtt E}\,\partial^2_{\hat\act_i\hat\act_j}\gitre=
 \K^{-\frac{37}2 n-2}\, \z (1\oplus \log \z)^2\,,
 $$
 where, the function in brackets belongs to $\cF$ and  has norm $\vvvert\cdot\vvvert_{\!\varrho}$ bounded by a constant depending only on $n$ and $\kappa$.\\
At this point, mimicking the proof for the inner case, one gets easily  \equ{gibbosobis} also in the outer case $i=0,2N$, if one chooses $\bfcu$ big enough.
The proof of of Proposition~\ref{zelenskyC} is  complete.
 \qed

\nl
{\bf Step 5:} {\sl The Twist Theorem far from separatrices in the inner case}

 \begin{proposition}\label{zelenskyA} Let $0<i<2N$. 
 
 \nl
 {\rm (i)}
There exists a      constant
$\bfcs=\bfcs(n,\upkappa)>1$
such that if 
$\K\geq \bfcs$,
$\Nf\leq\noruno{k}\leq\KO$, then,  on\footnote{Recall the definition of $\Bua=\Bua^i$ in \equ{tetto}.}  $\Bua^i$,
$|\det \partial^2_{\act}\mathtt{h} | \geq \dalla/{2^{5}}$.

\nl
{\rm (ii)}
There exists a  suitable    constant
 $\hcteot=\hcteot(n,\upkappa,\upxi, \mtt )\ge \bfcs$,  
such that if 
$\K\geq \hcteot$ and $\fico< \dalla/{2^{5}}$, then
\beq{gibbosoA}
\meas\big(
\big\{ \act\in \Bua^i:\ 
|\det \partial^2_{\act}
\mathtt{h} (\act) |
\le \fico
\big\}
\big)\le \hcteot (|k|^{2n}\fico)^{\frac{1}{\mtt}}\, \meas \Buik
\,.
\eeq
\end{proposition}
 
\rem
Notice that by point (i), the set $\big\{ \act\in \Bua:\  |\det \partial^2_{\act} \mathtt{h} (\act) | \le \fico \big\}$ is, for $\fico<\dalla/{2^{5}}$ empty. Therefore, in proving point (ii) one needs to consider only $\noruno{k}< \Nf$.
\erem
For definiteness, 
{\sl in the proof of Proposition~\ref{zelenskyA}, we consider only  
$i$ odd}, as the case $i$ even can be  treated in a completely analogous way.

 \nl
First, we prove some perturbative estimates on the derivatives of the energy.

\nl
Recalling the definition of $\vallo=\vallo^i$ in \eqref{tetto} with $\logeo$ as in Proposition~\ref{enigmista}--(ii), and notice that  
$\vallo$ only depends on $n,\upkappa$. Then, the following estimates hold.

\begin{lemma}\label{casetta}
 There exists 
 $\bfcc=\bfcc(n,\upkappa)>1$   such that, defining $\rstar:=  \sqrt\suca/\bfcc$, one has,
 for  $\act\in \vallo_{2\rstar}\times\hat D$, 
\begin{equation}\label{ofena2}
\left|\partial_{\act_1} \mathtt E \right|
\le \bfcc 
\sqrt\suca\,,
\qquad
\left|\partial^2_{\act_1} \mathtt E \right|
\le \bfcc 
\,,
\qquad
\left|\partial^2_{\act_1 \hat \act} \mathtt E \right|
\le \bfcc 
\lella
\,,
\qquad
\left|\partial^2_{\hat \act} \mathtt E \right|
\le \bfcc 
\lella
\end{equation}
and
\begin{equation}\label{soclose}
|\partial_{\act_1}{\mathtt E} -
\partial_{\act_1}\bar{\mathtt E}|
\le \bfcc  \sqrt\suca\lalla
\,.
\end{equation}
\end{lemma}
\proof
Recall the definitions in  \equ{arista}  and  \equ{wagner}.
Then:
\begin{eqnarray*}
\ts\act_1(E_+(\hat\act)-\suca\frac{\logeo}{8},\hat\act)-
\bar\bacci_{\logeo}
&\geq&\ts\phantom{-}
\bar\act_1(\bar E_+-\suca\frac{\logeo}{8})
\\
&& \ts-\bar\act_1(\bar E_+-\suca\frac\logeo{2})
-|\act_1(E_+(\hat\act)-\suca\frac{\logeo}{8},\hat\act)-\bar\act_1(\bar E_+-\suca\frac{\logeo}{8})|
\\
&\stackrel{\equ{vana},\equ{LEGOk}}\geq&\ts
\frac{1}{\bfc\sqrt\suca}
\left(\frac{\suca\logeo}{2}-\frac{\suca\logeo}{8}\right)
-
| \phi_+(\frac{\logeo}{8},\hat\act) -  \bar\phi_+(\frac{\logeo}{8})|\\
&&\ts
\phantom{\frac{1}{\bfc\sqrt\suca}\left(\frac{\suca\logeo}{2}-\frac{\suca\logeo}{8}\right)} 
- | \psi_+(\frac{\logeo}{8},\hat\act) -  \bar\psi_+(\frac{\logeo}{8})|\cdot \frac{\logeo}{8}|\ln\frac{\logeo}{8}|
\\
&
\stackrel{\equ{trippa}}\geq&\ts
\frac{3\sqrt\suca\logeo}{8\bfc}
-2\bfc \sqrt\suca\lalla
\stackrel{\equ{BWV1013}}\geq \frac{\sqrt\suca\logeo}{4\bfc}>0\,,
\end{eqnarray*}
which imply 
$\vallo\times\hat D\subseteq \Buik(\logeo/8)$.\\
Now, we can take $\bfcc>1$ big enough 
so that (recall \equ{blueeyes})
$\dst
 \vallo_{2\rstar}
\times\hat D
\subseteq
\big(\Buik (\logeo/8)\big)_{\rho_{\logeo/8}}\,.
$
Thus, by \eqref{ofena}  we get\footnote{Obviously the first two estimates holds 
also for $\bar{\mathtt E}=\mathtt E|_{\lalla=0}$. }  \eqref{ofena2}.
\\
Now, observe that
by the definitions in 
\eqref{tetto} and  \eqref{wagner} 
we get 
$
\bar{\mathtt E}(\vallo)=(\bar E_-, \bar E_+-\suca\logeo/2)
$.
Then, recalling  \eqref{autunno2}, by the first estimate in 	\eqref{ofena2},
we get
\begin{equation}\label{formichina}
\bar{\mathtt E}(\vallo_{2\rstar})
\subseteq
\mathcal E_{\logeo/8}\,,
\qquad
{\mathtt E}(\vallo_{2\rstar},\hat\act)
\subseteq
\mathcal E_{\logeo/8}\,,
\quad
\forall \, \hat\act\in\hat D\,,
\end{equation}
taking $\bfcc$ big enough.
\\
Let us, now,  prove  \eqref{soclose}.
Observe that
$$
\partial_{\act_1}{\mathtt E} (\act)-
\partial_{\act_1}\bar{\mathtt E} (\act_1)=
\big(
\partial_E \bar \act_1(\bar{\mathtt E} (\act_1))
-
\partial_E\act_1({\mathtt E} (\act_1),\hat\act)
\big)
\partial_{\act_1}{\mathtt E} (\act)\cdot 
\partial_{\act_1}\bar{\mathtt E} (\act_1)\,,
$$
so that
\begin{eqnarray*}
&&\sup_{
\vallo_{2\rstar}
\times \hat D
}
|\partial_{\act_1}{\mathtt E} (\act)-
\partial_{\act_1}\bar{\mathtt E} (\act_1)|
\\
&&\qquad \stackrel{\eqref{formichina}}\leq
\sup_{\mathcal E_{\logeo/8}\times \hat D}
\Big|
\partial_E \bar \act_1(E)
-
\partial_E\act_1(E,\hat\act)
\Big|
\cdot \sup_{
\vallo_{2\rstar}
\times \hat D
}
\Big|
\partial_{\act_1}{\mathtt E} (\act)\cdot 
\partial_{\act_1}\bar{\mathtt E} (\act_1)
\Big|
\stackrel{\eqref{rosettaTH},\eqref{ofena2}}
\lessdot \sqrt\suca\lalla\,. \qedeq
\end{eqnarray*}
Next, we provide perturbative estimates on the twist. 

\nl
By Cauchy estimates, from \eqref{soclose}, there follows
\begin{equation}\label{polpettedibollito}
\sup_{\vallo_{\rstar}
\times \hat D}|\partial^2_{\act_1}{\mathtt E}-
\partial^2_{\act_1}\bar{\mathtt E} |
\lessdot\lalla\,.
\end{equation} 
Hence,  by  
\eqref{parata} and \eqref{ofena2}, on $ \vallo_{\rstar}\times \hat D$, we get\footnote{Recall Notation~\ref{weiss}.}
\begin{eqnarray}\label{iut} 
\det \partial^2_\act  \mathtt{h}   
  =
   (\partial^2_{\act_1}{\mathtt E})
   \cdot \det \partial^2_{\hat \act} \hzk 
+
O(\lella )
=
  (\partial^2_{\act_1}\bar{\mathtt E}+O(\lalla))
   \cdot \det  \partial^2_{\hat \act} \hzk 
+
O(\lella )\,.
\end{eqnarray}
Now, by \equ{lovebuzz}, \equ{pappagallok},  one has that  $\dalla^{-1}\le \K^{2n}$ and $\lella/\dalla=O(\K^{-3n})$. Finally, since, by \equ{cerbiatta}, $\lalla=1/\K^{5n}$, from \equ{iut} one gets at once the following
\begin{lemma} Let $\rstar$ be as in Lemma~\ref{casetta}, $0\le i\le 2N$, 
and $\vallo=\vallo^i$ as in \equ{tetto}. Then,	
 \beq{sugo}
\left|\det \partial^2_\act   \mathtt{h}(\act)  \right|
\geq
\dalla
|g(\act)|
\,, \qquad \forall\, I\in  \vallo_{\rstar}
\times\hat D\,,
\eeq
with
\begin{equation}\label{sugotonno}
g(\act)=
\partial^2_{\act_1}\bar{\mathtt E} (\act_1)
+O(\K^{-3n})
 \,,\quad
\forall\, I\in 
 \vallo_{\rstar}
\times\hat D\,.
\end{equation}
\end{lemma}

\proof {\bf of Proposition~\ref{zelenskyA}}
(i) Since $\noruno{k}\ge \Nf$, $\GO= {\textstyle \frac{2\e}{|k|^2}}\, \fproj f$ in \equ{paranoia}  is close to a cosine, as proved in  Lemma~\ref{pollaio}  in \Appendix. Hence,  \equ{cosinelike} in Proposition~\ref{kalevala}  holds, so that by \equ{sugo} and \equ{sugotonno}, taking $\bfcs$ large enough and $\K\ge \bfcs$, the claimed estimate $|\det \partial^2_{\act}\mathtt{h} | \geq \dalla/{2^{5}}$ follows.

\nl
(ii)
Recall \equ{wagner}.
Since $\loge\to\bar\bacci_\loge$ is a decreasing function,
 we get
$\bar\bacci_{\logeo/2}\leq \bar\bacci=\bar\bacci_0$.  
Rescaling we get
$\tilde \vallo:=(0,\bar\bacci_{\logeo/2}/\bar\bacci)\subseteq (0,1)$ so that $\bar\bacci \tilde \vallo= \vallo$.
Recalling \eqref{arista}, 
by \eqref{LEGOk}$\div$\eqref{ciofecak}
we have that 
$ \bar\bacci\lessdot
\sqrt\suca$. Then, choosing $0<\tilde{\mathtt r}\leq 1$ small enough, we have that
\begin{equation}\label{sputin}
\bar\bacci \tilde \vallo_{2\tilde{\mathtt r}}
\subseteq
 \vallo_{\rstar}\,.
\end{equation}
By \eqref{nutella} we get
\begin{equation}\label{barbarossa}
\partial^2_{\act_1}\bar{\mathtt E}  (\act_1)=
{\mathtt F}_\GO (\act_1/\bar\bacci)\,.
\end{equation}
 By \eqref{sugotonno} and \eqref{sputin}
 we get
\begin{equation}\label{sugotonno2}
g(\bar\bacci x,\hat\act)=
{\mathtt F}_\GO(x)
+O(\lalla)
 + O(\lella /\dalla)\quad
 \mbox{uniformly for}\quad 
 (x,\hat\act)\in
\tilde \vallo_{2\tilde{\mathtt r}}
\times\hat D\,.
\end{equation}
By \eqref{ofena2}, \eqref{sputin} and \eqref{barbarossa} we get
\begin{equation}\label{ofena3}
\sup_{\tilde \vallo_{2\tilde{\mathtt r}}}|{\mathtt F}_\GO|\lessdot 1\,.
\end{equation}
For 
$\hat \act\in\hat D$, 
set\footnote{Recall that $ \vallo:=(0,\bar\bacci_{\logeo/2})$
in \eqref{tetto}.}
\begin{eqnarray}\label{ippona}
\mathcal \act_\fico'(\hat \act)
:=
\big\{ \act_1\in
\vallo:\,
\big|\det \partial^2_\act \mathtt{h}^i(\act)
 \big|\leq \fico  \big\}\,,
\quad
\tilde{\mathcal \act}_\fico'(\hat \act)
:=
\{x\in\tilde \vallo:\,  
|g(\bar\bacci x,\hat\act)|\leq \fico/\dalla
\}\,.
\end{eqnarray}
By \eqref{ippona} and \eqref{sugo}
we have that, for every $\hat\act\in \hat D$,
\begin{equation}\label{rosmarino}
\meas \mathcal \act_\fico'(\hat \act)
\leq \bar\bacci\meas \tilde{\mathcal \act}_\fico'(\hat \act)
\lessdot
\sqrt\suca\, \meas \tilde{\mathcal \act}_\fico'(\hat \act)\,.
\end{equation}
using $ \bar\bacci\lessdot
\sqrt\suca$.
Before estimating $\meas \tilde{\mathcal \act}_\fico'(\hat \act)$
we need the following bound:
\begin{equation}\label{pelino}
1/2\leq{\bar\bacci_{\logeo/2}}/{\bar\bacci}\,.  
\end{equation}
Recalling \eqref{arista} we have
\begin{equation}\label{stefania}
\bar\bacci=\bar\act_1(\bar E_+)
=\int_{\bar E_-}^{\bar E_+}\partial_E\bar\act_1
\stackrel{\eqref{vana}}\geq
\frac{\bar E_+-\bar E_-}{\bfc\sqrt\suca}
\stackrel{\eqref{latooscuro}}\geq \frac{\sqrt\suca}{\bfc\upkappa}\,.
\end{equation}
Recalling \eqref{tetto} 
we have by\footnote{Note that condition \eqref{bassoraTH}
reduces here to $\loge\leq 1/\bfc$ since we are considering
$\bar\act_1$, namely the case $\lalla=0.$
In any case one can prove the estimate also directly
by \eqref{LEGOk} and \eqref{pappagallok}.} 
\eqref{rosettaTH} that, for $0<\z\leq \logeo/2,$
$
\partial_E \bar\act_1(\bar E_+-\suca\z)
\leq {\bfc^2}{}|\log \z|/\sqrt\suca$.
Therefore
\begin{eqnarray*}
1-\frac{\bar\bacci_{\logeo/2}}{\bar\bacci}
&=&
\frac{\bar\act_1(\bar E_+)-\bar\act_1(\bar E_+-\suca\logeo/2)}{\bar\bacci}
=
\frac{\suca}{\bar\bacci}\int_0^{\logeo/2}
\partial_E \bar\act_1(\bar E_+-\suca\z)d\z
\\
&\stackrel{\eqref{stefania}}\leq& \bfc^3\upkappa
\int_0^{\logeo/2}|\log\z|d\z
\leq \bfc^4
 \logeo\log|\logeo| \stackrel{\eqref{BWV1013}}\leq \frac12
\,,
\end{eqnarray*}
proving \eqref{pelino}.
\\
Let us come back to the estimate of
$\meas \tilde{\mathcal \act}_\fico'(\hat \act)$.
Recalling Definition \ref{dracula} we have that
${\mathtt F}$ is
$\upxi$--non--degenerate at order $\mtt $.
By \eqref{sugotonno2}, \eqref{ofena3} and Cauchy estimates,
taking 
 $\lalla$ and\footnote{$\dalla$ is defined in \equ{balla}.} $\lella/\dalla$ small enough (i.e.,  $\K\geq \hcteot$ for a suitable  $\hcteot$ arge enough) depending only on $\upkappa$, $n$, $\upxi$ and $\mtt $
 we have that the function $x\mapsto g(\bar\bacci x,\hat\act)$
 is $(\upxi/2)$--non--degenerate at order $\mtt $.
Now  we  want to apply  Lemma~\ref{nessuno} with $\fico$ replaced by  
$\fico/\dalla$, and with the following choices:
$$
f(x)=g(\bar\bacci x,\hat\act) \,,\quad
 m= \mtt\,,\quad
  a = 0\,,\quad
  1/2\stackrel{\eqref{pelino}}\leq
   b= \bar\bacci_{\logeo/2}/\bar\bacci<1\,,\quad
 \xi=\upxi/2\,;
$$
the constant $M$ controlling the derivatives of $f$, by \eqref{sugotonno2}, \eqref{ofena3}
and Cauchy estimates can be bounded by  
$1\leq M\leq c_{n,\upkappa}/\tilde{\mathtt r}^{\mtt+1}$
for\footnote{Recall that $\tilde{\mathtt r}=\tilde{\mathtt r}(n,\upkappa)$ was chosen in \eqref{sputin}
small enough.} for a suitably large constant $c_{n,\upkappa}$ depending only on 
$n$ and $\upkappa$. 
In conclusion, by Lemma~\ref{nessuno},
we get
$$
\meas \tilde{\mathcal \act}_\fico'(\hat \act)
\leq c_\mtt 
\Big(\frac{2c_{n,\upkappa}}{\upxi\tilde{\mathtt r}^{\mtt+1}}+1\Big)
\Big(\frac{\fico}{\dalla\upxi}\Big)^{\frac{1}{\mtt}}\,.
$$
Then \eqref{gibbosoA} follows by \equ{lovebuzz}, \eqref{rosmarino}
Fubini's theorem and \eqref{natalino}. The proof of Proposition~\ref{zelenskyA} is complete.
\qed

\nl
{\bf Step 6:} {\sl Uniform twist  in outer regions far from separatrices}

\nl
Recall the definition of the twist $\dds$ in the outer regions in  \equ{enumaelish}, and that $\dalla=|k|^{-2n}$ (see \equ{balla}).

\begin{proposition}\label{zelenskyAL}
Let $i=0,2N$. Then, 
there exists a  suitable    constant
$\bfcst=\bfcst(n,\upkappa)>1$     
such that if 
$\K\geq \bfcst$,  then on $\Bua^i$
 \beqno
\ts |\dds|\geq  \dalla/2\,.
\eeqno
\end{proposition}
\proof
 Taking $\loge=\logeo/2$ 
 in \eqref{ofena}, on $\Bua^i$, we have that\footnote{The quantities ${\rm \hat v}$ and  ${\rm \hat M}$ are defined in \equ{daitarn4}.}
\begin{equation}\label{cajato}
|{\rm \hat v}|,\ |{\rm \hat M}|\ \lessdot {1}/{\K^{\frac{19}{2}n-1}}\,.
\end{equation}
By  \eqref{cerbiatta} and \eqref{pappagallok}, we also have
$$
\frac{\sqrt\suca}{\ro}\lessdot 
\frac{1}{\K^{\frac{9}{2}n}}
\,,\qquad
\lella\lessdot\frac{1}{\K^{\frac{19}{2}n+1}}\,.
$$
Then, recalling Definition~\ref{piggy} and 
\eqref{cerbiatta},
we get
$
|{\mathtt E}|\leq 2\Ro^2\leq
2\e\K^{9n+4}$.    
\\
Now, it is a general fact that, in the outer case, the unperturbed energy function is strictly concave, as it follows from the following simple consequence of Jensen's inequality\footnote{ $\bar E_i$ is defined in \equ{cesarini}.}.

\begin{lemma}\label{jensen} 
Let  $i=0,2N$. Then, for every 
$E> \bar E_0=\bar E_{2N}$, $\partial^2_{\act_1}\bar{\mathtt E}^{i}(\bar \act_1^{i}(E))
\geq 2$ . 
\end{lemma}
The proof is given in \Appendix.

\nl
Now, since  estimate \eqref{polpettedibollito} still holds in the present case
$i=2N$ we get by Lemma~\ref{jensen},
$
\partial^2_{\act_1} {\mathtt E}\geq \frac12 \partial^2_{\act_1} \bar{\mathtt E}
\geq 1$, so that  the claim
 follows by \eqref{enumaelish}, \eqref{cajato}
and \eqref{lovebuzz}. \qed

\nl
{\bf Step 7:} {\sl Conclusion}

\nl
\proof {\bf of the Twist Theorem~\ref{zelensky}}
Let $\logeo=\logeo(n,\kappa)$ be as in Proposition~\ref{enigmista}--(ii), and 
let $\bfcot=\bfcot(n,\kappa)\ge 1$ be such that the second estimate in \equ{BWV1013} holds if\footnote{Recall \equ{cerbiatta}.} $\K=1/\lalla^{5n}\ge \bfcot$.
Then, by Lemma~\ref{mediatore}, 
\beq{repetita}
\Buik=\Buc^i(\logeo)\cup\Bua^i(\logeo)\,,\quad\quad \forall\ 0\le i\le 2N\,.
\eeq
Define\footnote{Recall that $\bfcu\ge \bfcq\ge \max\{\bfcd,\bfct\}\ge \itct$. }
\beq{torvaianica}
\cteot:=  2\max\{ \bfcu,  \bfcs, \bfcc, \bfcst, \bfcot, \hcteot\}\,.
\eeq
Let us consider first the outer case $i=0,2N$. Recall the definition of $\lippo$ in \equ{gibboso}.
By Lemma~\ref{depechemode}, Proposition~\ref{zelenskyAL}, Proposition~\ref{zelenskyC}, and by \equ{torvaianica}, we find
\beqano
\meas\big(\big\{ \, \act\in \Bik:\  \big| \det \partial^2_{\act}\mathtt{h}^i_k(\act) \big|\leq \fico\big\}\big)&=&\meas\big(\big\{ \, \act\in \Buik:\  \big| \dds^i(\act) \big|\leq \fico \big\}\big)\\
&=&\meas\big(\big\{ \act\in \Buc^i(\logeo):\ |\det \partial^2_{\act}\mathtt{h}^i (\act) |\le \fico\big\}\big)\\
&\le&\bfcu \, (|k|^{2 n}\fico)^{1/9 n^4}\, \meas \Buik \\
&\stackrel{\equ{pizza}}= &\bfcu \, (|k|^{2 n}\fico)^{1/9 n^4}\, \meas B^i_k\\
&\stackrel{\equ{torvaianica}}< &\cteot \, (|k|^{2 n}\fico)^{\lippo}\, \meas B^i_k\,,
\eeqano
proving Theorem~\ref{zelensky} in the outer case $i=0,2N$. 

\nl
In the inner case $0<i<2N$, $B^i_k=\Buik$ (compare \equ{tordobis}) and, since $\K\ge \cteot$, \equ{gibboso} follows by \equ{repetita}, \equ{gibbosobis} in Proposition~\ref{zelenskyC}, and \equ{gibbosoA} in Proposition~\equ{zelenskyA}. \qed


\section{Maximal KAM tori and proof of the main results}

In this final section we show that  primary and secondary maximal KAM tori of $\ham$ span the   complementary of $\Rd\times\T^n$ apart from an exponentially small (in $1/\K$)  set and prove the results in Section~\ref{omnia}.

\nl
To construct such tori we shall use the following `KAM theorem'.
 
 \begin{theorem}[\cite{BCKAM}]
 \label{KAM}
Fix $n\ge 2$ and let 
$\ttD$ be any non--empty, bounded subset of $\real^n$. 
Let
\beqno
\rm H( p, q):= h(p)+f(p,q)
\eeqno
be  real analytic on $\ttD_\rkam\times\T^n_\skam,$
 for some $\mathfrak r>0$ and $0<\skam\le 1$, and having finite norms
\beq{bellini}
\rm \ttM:=  \modulo\partial^2_p h{\modulo}_\rkam \,,\qquad\qquad
{\modulo}f{\modulo}_{\rkam,\skam} \,.
\eeq
Assume that the frequency map $\rm p\in \ttD \to\o=\partial_p h$ is a local diffeomorphism, namely, assume:
\begin{equation}\label{dupa}
\rm d:=  \inf_\ttD|\det \partial^2_p h| >0\,,
\end{equation}
and let
$\rm {\rm d}_*:={d}/{\ttM^{\it n}}$ and  $r_*:= \, {\rm d}_*^2 \rkam$.
Then, there exists  $\ckam=\ckam(n)>1$ such that, if   
\beq{enza} 
\rm \epsilon:= \frac{{\modulo}f{\modulo}_{\rkam,\skam} }{\ttM \rkam^2}\le \frac{{\rm d}_*^8\ \skam^{4(n+1)}}{\ckam}\ , 
\eeq
there exists a  set $\rm{\cal T}\subseteq (\ttD_{{r}_*}\cap\R^n)\times \torus^n$ formed by primary KAM tori such that\footnote{Here `meas' denotes the outer Lebesgue measure.} 
\beq{gusuppo}
\meas \big((\ttD\times\torus^n)\bks{\cal T}\big) \le  \rm C\, \sqrt\epsilon\,,\qquad \rm C:=
\big(\max\big\{ {\rm d}_*^2 \rkam\, ,\diam \ttD\big\}\big)^n \cdot \frac{\ckam}{{\rm d}_*^{{\it n}+5}\ \skam^{3(n+1)}}\,.
\eeq
\end{theorem}

\nl
This statement is an immediate corollary of 
Theorem~1 in\footnote{In Theorem 1 of \cite{BCKAM}  take $\t=n$ and substitute 
$\l$ with its maximal value  $2\cdot n! \, {\rm d}_*^{-1}$ (see  (14) of \cite{BCKAM}).
}
\cite{BCKAM}. 

\rem
(i) Note that in the formulation of Theorem~\ref{KAM} the action domain $\ttD$ is a completely arbitrary bounded set and that the smallness quantitative condition \equ{enza} depends on $\ttD$ only through its  diameter, which in our application depends on $k$.
For a similar statement, which takes into account the geometry of $\ttD$, see \cite{CK}.
\\
(ii)
We point out that  the smallness condition \eqref{enza} can be rewritten as  
    \beq{enzabis} 
\rm {\modulo}f{\modulo}_{\ttD,\rkam,\skam}
\le  \frac{ \rkam^2 d^8\ \skam^{4{\it n}+4}}{\ckam \,\ttM^{8{\it n}-1}}\,.
\eeq
(iii) Finally, observe that, since\footnote{Indeed the absolute value of any eigenvalues
 of the symmetric matrix $\rm \partial^2_p h$ is bounded by $\ttM$,
 which implies
$\rm d\leq \sup_\ttD|\det \partial^2_p h|\leq \ttM^{\it n}$. }  ${\rm d}_*\leq 1$,
estimate \eqref{gusuppo} implies 
\beq{gusuppobis}
\rm
\meas \big((\ttD\times\torus^{\it n})\bks{\cal T}\big) 
\le\rm
\big(\max\big\{ \rkam\, ,\diam \ttD\big\}\big)^{\it n} \cdot \frac{\ckam\, \ttM^{{\it n}^2+5{\it n}-1/2}}{d^{{\it n}+5}\ \skam^{3{\it n}+3} \rkam}
\, \sqrt{{\modulo}{\rm f}{\modulo}_{\ttD,\rkam,\skam} }
\,.
\eeq
\erem

\subsection*{KAM tori in the non--resonant region}

\begin{proposition}\label{pecco} Let the assumptions of Theorem~\ref{normalform} hold.
There exists a constant $\ttcnr=\ttcnr(n,s)\geq \bfco$ such that, if
$\KO\ge  \ttcnr$, then
 there exists a family of primary
maximal KAM tori ${\cal T}^0$ invariant for the Hamiltonian 
$\ham$ in \equ{ham}, satisfying
\beq{lazie}
 \meas \big((\Rz\times\torus^{\it n})\bks {\cal T}^0\big)  \le \ttcnr\, \sqrt\e\, e^{-\KO s/6}\,.  
\eeq
\end{proposition}

\rem 
The above result is essentially classical, and, in fact,  no genericity assumptions on the potentials are needed. However, there is one delicate point related to the KAM tori near the boundary. Indeed, primary tori oscillates, in general, by a quantity of order $\sqrt\e$, and naive applications of classical KAM theorems would leave out regions near the boundary of the phase space of measure $\sim \sqrt\e$. Such a problem is overcome by using the second covering  in \equ{codino} in Theorem~\ref{normalform}, which is introduced so that  \equ{surge} holds; compare, also, Remark~\ref{munro}--(ii).
\erem

\proof {\bf of Theorem~\ref{KAM}} 
We apply the KAM Theorem~\ref{KAM} to the nearly--integrable Hamiltonian $\hamo$ in Theorem~\ref{normalform}--(ii). More precisely, we let\footnote{Recall Theorem~\ref{normalform},  
the definitions \equ{athlone}, \eqref{dublino} and \eqref{neva}.}
\begin{eqnarray*}
{\rm h}({\rm p})=\frac{|{\rm p}|^2}2+\e g^{\rm o}({\rm p})\,,
\ \ 
{\rm f}=\e
 f^{\rm o}\,,
 \ \ 
 \ttD=\Rzt\,,
 \ \
 \rkam=\frac{r_{\rm o}'}{2}=\frac{\sqrt\e \K^{\frac{9}{2}n+2}}{16\KO}\,,
 \ \
 \skam=\min\{\ts\frac{s}{2}\,,1\}\,.
\end{eqnarray*}
By \eqref{552} and  Cauchy estimates we get
$$
\ttM\leq 2\,,\qquad
{\modulo}{\rm f}{\modulo}_{\rkam,\skam}
\leq  \e e^{-\KO s/3}\,,\qquad
\rm d\geq 1/2\,.
$$
If $\KO$ is taken large enough (larger than a constant despending on $n$ and $s$)
the KAM smallness condition \equ{enzabis} is satisfied, and the KAM Theorem~\ref{KAM} yields the existence of 
 a set $\widetilde{\cal T}^{0}$ of invariant tori for the Hamiltonian 
$\hamo$ in Theorem~\ref{normalform}--(ii), which, by\footnote{Notice that the hypothesis $\K<\e^{-1/(9n+4)}$ implies that $\rkam<1$, so that
$\max\big\{ \rm d^2 \ttM^{-2{\it n}} \rkam\, ,\diam \ttD\big\}= 2$.}  
\eqref{gusuppobis}, satisfy 
\beq{dignita}
\meas \big((\Rzt\times\torus^{\it n})\bks \widetilde{\cal T}^{0}\big)   
\leq
\ttcnr \sqrt\e\,e^{-\KO s/6}\,,   
\eeq
for a suitable constant  $\ttcnr=\ttcnr(n,s)$ large enough (so that also the condition on $\KO$ is met).
Since the  map $\Psi_{\! \rm o}$  in \eqref{surge} is symplectic,  
the family of tori 
$
{\cal T}^0:= \Psi_{\! \rm o}(\widetilde{\cal T}^{0})
$
 is formed by KAM invariant for $\ham$ in \equ{ham}. The first relation in
\equ{surge} and the bound \equ{dignita} imply \equ{lazie}. \qed

\subsection*{KAM tori near simple resonances}

Now, we turn to the construction, in all neighbourhoods of simple resonances,  of families of  primary tori for the nearly--integrable Hamiltonians  $\cH_k^i$ of Theorem~\ref{garrincha}, for all $k\in\genKO$ and $0\le i\le 2N_k$. Note that such tori correspond, in the inner case $0<i<2N_k$, to  {\sl secondary tori for the  Hamiltonian $\ham$}.

\nl
Let us introduce  zones $\Bik(\loge,\fico) \subseteq \Bik$, which are  $\loge$--away in energy from separatrices and   where the twist is bounded away from zero by a quantity $\fico>0$, namely (recall \equ{tales}, \equ{tordobis}), let us define:
 \beq{topgun}
 \Bik(\loge,\fico):=\{ \act\in \Bik(\loge)  \ {\rm s.t.}\  
|\det \partial^2_{\act}\mathtt{h}^i_k(\act)|> \fico\}\subset \Bik\,.
\eeq

\begin{proposition}\label{noAC} {\bf (KAM tori for $\cH^i_k$)} 
Let the assumptions of Theorem~\ref{garrincha} hold.  There exist positive constants $\ttcq=\ttcq(n,s,\b)>1$ and 
${\ttCu}=\ttCu(n,s,\b,\d) \ge \bfcast$  
such that the following holds.
Let $k\in\genKO$, $0\le i \le 2N_k$; 
$0<\loge\leq 1/\bfcast$ and $0<\fico<1/2$. Then, if
\begin{equation}\label{svitolina}
\K\geq {\ttCu}\ln \frac{1}{\loge \fico}\,,
\end{equation} 
there exists a set 
${\cal T}^i_k$ of  maximal KAM  tori for the Hamiltonian $\cH_k^i$ 
 in \eqref{kant} such that 
\beq{gusuppotto}
\meas \big( (\Bik(\loge,\fico)\times\torus^{\it n})\bks {\cal T}^i_k\big)
\le
\ttcq e^{-\K s/7}
\,.
\eeq
\end{proposition}

\proof
We apply the KAM Theorem~\ref{KAM} 
to the Hamiltonian 
$\cH^i_k$ of Theorem~\ref{garrincha} with (recall
\equ{kant} and \eqref{tales}):
\begin{eqnarray}
&&{\rm h}=\hik=
{\ts \frac{|k|^2}{2}}\mathtt{h}^i_k
\,,\qquad
{\rm f}=\e \fik
\,,
 \qquad
 \ttD=\Bik(\loge,\fico)\,,
\nonumber
\\
 &&
 \rkam=\rhs=\frac{\sqrt\suca} {\bfcast\KO^{n}}\, \loge |\log \loge| 
\,,
 \qquad\quad
 \skam=\shs=\frac{1}{\bfcast\KO^{n}|\log \loge|}
\,.
\label{astenuti}
\end{eqnarray}
Note that, by \equ{sgallettata} and \equ{kappa},   $0<\loge\leq 1/\bfcast\le1/8\itcd$, which implies easily   $\rkam\le \ro$ and
 $\skam\le 1$. Also, since $\bfcast\ge \hbfc$ (see Theorem~\ref{garrincha}) and $\KO^n\ge 2^n\ge n$, one has $\rhs\leq \rhol/n$.
 
 \nl
 In the following arguments  we denote by $c(\cdot)$
possibly different constants depending only on the quantities inside brackets.
 \\
 We first have to estimate $\ttM$ in
 \equ{bellini}, namely,  $\partial^2_\act \mathtt{h}^i_k$.
 By \eqref{ofena}, \equ{pappagallok}
 and \eqref{caviale2} we get
 \beq{medie}
\sup_{(\Buik(\loge))_{\rhol}}\big|\partial^2_{\act} \mathtt E^i \big|
\leq \frac{n \hbfc}{\loge} \,.
\eeq 
In the case
$0<i<2N_k$, by  \equ{tales}, we have $\Buik(\loge) =\Bik(\loge)$. Therefore, 
recalling \equ{kant}, 
we can bound  $|\partial^2_\act\mathtt{h}^i_k|$ by $c(n,s,\b)/\loge$.\\
The estimate on $|\partial^2_\act\mathtt{h}^i_k|$ in the case $i=0,2N_k$ needs some extra attention.  
In particular  fix
 $i=2N_k$ (the case $i=0$ being analogous). 
Recalling the definition of $\Istar$ in 
\equ{tordissimo}, \equ{tordobis} we have that 
$\partial^2_{\act} \Istar$ depends only on $\hat\act$ and not
on $\act_1$.
Moreover
by \equ{athlone}, \equ{atlantide},  \eqref{tess},
\equ{cerbiatta}
and Cauchy estimates  
we get
 \beq{voloillatte} 
  \sup_{\hat \act\in \hat D_{3\ro}}|\partial_{\act } \Istar|\leq c(n)
  \,,\qquad
  \sup_{\hat \act\in \hat D_{3\ro}}|\partial^2_{\act } \Istar|\leq 
  \frac{c(n)|k|^2}{\sqrt\e \K^\nu}\,.
  \eeq
 Recalling Definition~\ref{piggy}, \equ{ofena}
  and  \eqref{cerbiatta},
   we have that 
  $$
  \sup_{(\Bu^{2N_k}_k(\loge))_{\rhol}}
  |\mathtt E^{2N_k}|\leq 4\Ro^2=\frac{4\e \K^{2\nu}}{|k|^4}
  \,.$$
 Then, by \eqref{ofena}
   and \eqref{cerbiatta} we get 
   $$
    \sup_{(\Bu^{2N_k}_k(\loge))_{\rhol}}
   |\partial_{\act_1}\mathtt E^{2N_k}|\leq
   \hbfc\sqrt{8\ttcs\e+4\e \K^{2\nu}|k|^{-4}}
   \leq 
  4\hbfc\sqrt{\e} \K^{\nu}|k|^{-2}$$
    (taking $\K\geq \ttcs$  
    defined in \equ{cerbiatta}).
 Finally, recalling also
\eqref{tales}, \eqref{cioccolata}, \equ{medie},
\eqref{voloillatte} , we  get by the chain rule
 $$
\sup_{(B_k^{2N_k}(\loge))_\rkam}\big|\partial^2_{\act} \big(\mathtt E^{2N_k}\circ \act_*\big) \big|
\leq \frac{c(n,s,\b)}{\loge} \,.
$$
By \eqref{bellini}, \equ{astenuti}, \eqref{kant}, \eqref{lovebuzz}  
(and  that $\rkam\leq  \chr$),
we finally get
\begin{equation}\label{liceo}
\ttM\leq |k|^2\frac{c(n,s,\b)}{\loge}\,,\qquad \forall\ 0\le i\le 2N_k\,.
\end{equation}
Next, by \eqref{astenuti} and \eqref{schiaffino}, 
\begin{equation}\label{apocalisse2}
{\modulo}{\rm f}{\modulo}_{\rkam,\skam}
\leq
\e\, e^{-\K s/3}\,.
\end{equation}
By \eqref{dupa}, \eqref{astenuti} and \equ{topgun},
we get
 \begin{equation}\label{apocalisse3}
 {\rm d}\geq 2^{-n}|k|^{2n} \fico\qquad \mbox{and}\qquad
 \frac{\ttM^n}{{\rm d}}\leq \frac{c(n,s,\b)}{\loge^n \fico}\ .
\end{equation} 
By \eqref{cerbiatta},   \eqref{P1+}
and using $\KO\leq 6\K $,  
we have:
\beq{apocalisse4}
\frac{\e}{\suca} 
  \leq\frac{\KO^{n+2}}{8 \ttcs\d}\ e^{\KO s}
  {\leq}\frac{\K^{n+2}}{6^{n+3} \ttcs\d}\ e^{\K s/6}\,.
\eeq
It is now  easy to check, by \equ{liceo},
 \equ{apocalisse2}, \equ{apocalisse3} and \equ{apocalisse4}, that the KAM smallness condition \eqref{enzabis} is satisfied
taking $\K$ as in \equ{svitolina} with $\ttCu$ large enough.
By the KAM Theorem~\ref{KAM} we, then, obtain a set ${\cal T}^i_k$ of invariant tori for the Hamiltonian 
 in \eqref{kant}, which, in view of  \eqref{gusuppobis} and by \equ{liceo},
 \equ{apocalisse2}, \equ{apocalisse3} and \equ{apocalisse4}, satisfies
 \equ{gusuppotto} with a suitable constant $\ttcq=\ttcq(n,s,\b)$; 
 in particular, note that, by  \equ{diametro} and \equ{astenuti}, the maximum in \eqref{gusuppobis}
 is estimated by $c(n)\KO^{n^2}$.
 \qed

\nl
Putting together these KAM statements and the Twist Theorem~\ref{zelensky}, 
 the proof of the  results stated in Section~\ref{omnia} follow easily.

\subsection*{Proof of Theorem \ref{prometeo} and its corollaries} 

By Lemma~\ref{telaviv}, since $f\in \Gns$, there exist $\d,\b>0$ such that 
\equ{P1+} and \equ{P2+} hold with $\Nf$ as in \equ{enne}. Let 
\beqno
\KO:=\K/6\,,
\eeqno
with $\K\ge 12$ and let $\a$ be as in   \equ{athlone}. Then, Assumptions~\ref{assunta} hold, and we may let the Definitions~\ref{assunto} hold.
Let $\itcd=\itcd(n)$ and $\bfco=\bfco(n,s,\d)$ be as in Theorem~\ref{sivori}, and
assume that\footnote{Eq. \ref{K1} implies that $\K>12$.}
\beq{K1}
\K\ge 6 \max\{\itcd,\bfco\}\,.
\eeq
Then, Theorem~\ref{sivori} holds and we may define the parameters $\upxi>0$ and  $\mtt \geq 1$ as in Definition~\ref{dracula} with respect to standard Hamiltonians $\Hsharp$ (with $\noruno{k}\le \KO$) of Theorem~\ref{sivori}--(ii).

\nl
We now let $\lippo<1$ as in \equ{gibboso}, $\ttCu=\ttCu(n,s,\b,\d)$ be as in Proposition~\ref{noAC},
and define 
\begin{equation}\label{nadal}
\fico:=e ^{- \frac{\K}{{\ttCu}(1+\lippo)}}\,,\qquad \loge:=\fico^\lippo\,.
\end{equation}
Notice that, with such definitions, it is  
\beq{pamuk}
\K= {\ttCu}\ln \frac{1}{\loge \fico}\,,
\eeq
(compare \equ{svitolina}).

\nl
With these premises, let us turn to the proof of the claims of Theorem~\ref{prometeo}.

\nl
Claim (ii) has already been  proven in  Lemma~\ref{coverto} above.

\nl
Next, we define the set of maximal KAM tori $\cT$ for $\ham$ as it appears in item (iv) of the theorem.

\nl
Let $\ttcnr=\ttcnr(n,s)$ as in Proposition~\ref{pecco}.There exists a constant 
$$\cgoth=\cgoth(n,s,\b,\d,\mtt)\ge \max\{\ttcnr, 2\ttCu \bfcast/\lippo \}\,,
$$
such that, if $\K\ge \cgoth$, then
\beqno
\quad \K^{2n}\fico \stackrel{\equ{nadal}}= \K^{2n} e ^{- \frac{\K}{{\ttCu}(1+\lippo)}} \le 1\,.
\eeqno
Assume that \beq{K2}
\K\ge \cgoth\,.
\eeq
Then, $\loge=\fico^\lippo$ in \equ{nadal} is smaller than $1/\bfcast$ and (recall \equ{balla})
\beq{falla}
\fico\le \frac\dalla{2^5}<\frac12\,.
\eeq
Thus, in view of \equ{pamuk}, by \equ{K2} the assumptions of Propositions~\ref{pecco} and \ref{noAC} are satisfied, and
 we can define the following families of  tori\footnote{${\cal T}^i_k $ is defined    in Proposition~\ref{noAC}, 
 $\cT^0$ in Proposition~\ref{pecco} and $\upphi_k^i$  in \equ{tordo}.
 }:
\beq{torio}
\left\{
\begin{array}{l}
 {\cal T}^{1,k}_i:=
\upphi_k^i({\cal T}^i_k )
\,,\quad
\dst {\cal T}^{1,k}:=\bigcup_{0\leq i\leq 2N_k} {\cal T}^{1,k}_i\,,\quad
 {\cal T}^1:=\bigcup_{k\in\genKO}{\cal T}^{1,k}\,,
\\ 
{\cal T}:={\cal T}^0\cup{\cal T}^1\,.
\end{array}\right.
\eeq
Observe that ${\cal T}^i_k$ are invariant tori for  $\cH^i_k$ in \equ{kant},  while
${\cal T}^{1,k}_i$, $\cT^1$ and $\cT^0$ are invariant for the original Hamiltonian $\ham$. 
\\
Thus,  ${\cal T}$ is a  family of maximal KAM tori for $\ham$ as in   item (iv) of Theorem~\ref{prometeo}.

\nl
Claim (i)  follows, now,  immediately by \equ{sonnosonnosonno}, setting
\beq{arno}
\cA:=\big((\Rz\cup \Ru)\times\torus^n\big) \bks {\cal T}\,.
\eeq
It remains to prove claim (iii), namely, the exponential measure estimate on $\cA$.

\nl
Observe that by \equ{arno} and \equ{torio}
\beq{insiemi}
\cA\subseteq \big((\Rz\times\torus^n) \bks {\cal T}^0\big) \cup \big((\Ru\times\torus^n) \bks {\cal T}^1\big)
\subseteq
\big((\Rz\times\torus^n) \bks {\cal T}^0\big) \cup
\bigcup_{k\in\genKO} 
(\Ruk\times\torus^n) \bks\cT^{1,k}
\,.
\eeq
We now need the following elementary result, whose proof is given  in \Appendix.
\begin{lemma}\label{geremia} If $f\in \Bns$ satisfies \equ{P2+}, then,
 for any
  $k\in \gen$,   the number  $2N_k$ of critical points  of 
  $\fproj f$ is bounded by ${\bar{\mathtt c}}:=\max\{4,\pi\sqrt{8/\b}\}$.
\end{lemma}
Obviously, the hypothesis of this lemma are met by our fixed potential in $\Gns$, and the following measure estimate  holds.

\begin{lemma}
Let $\loge$ as above in \equ{nadal} and ${\bar{\mathtt c}}$  as in Lemma \ref{geremia}.
Then,
for any $k$ in $\genKO$, one has 
\beq{periferia1}
\begin{array}{l}
\phantom{aa}
\meas\big((\Ruk\times \torus^n)\bks \cT^{1,k}\big)
\\
\le \dst
\bfcast
\meas\big(\Rukt\times\T^n\big)
\loge |\log \loge|
 +{\bar{\mathtt c}}\, \max_{0\le i\le 2N_k} 
\meas\Big(\big(\Bi_k(\loge)\times\T^n  \big)
\setminus {\cal T}^i_k \Big)\,.
\end{array}
\eeq

\end{lemma}
\proof  
Since $\upphi_k^i$   in Theorem~\ref{garrincha} is a diffeomorphism, one has
\beqano
&&
(\Ruk\times\torus^n) \bks\cT^{1,k}
\stackrel{\equ{torio}}=
(\Ruk\times\torus^n) \bks
\Big(\bigcup_{0\leq i\leq 2N_k}
\upphi_k^i({\cal T}^i_k )\Big)
\\
&\subseteq&
\Big(
\big(\Ruk\times\T^n\big)\ \setminus\ 
\bigcup_{0\leq i\leq 2N_k}
\upphi_k^i
\big(\Bik(\loge)\times \T^n\big)
\Big)
\cup
\bigcup_{0\leq i\leq 2N_k}
\upphi_k^i\Big(
\big(\Bik(\loge)\times \T^n\big)
\setminus 
{\cal T}^i_k
\Big)\,,
\eeqano
then, passing to measures,  using \equ{inthecourt2}, the fact that $\upphi_k^i$ is symplectic  and 
Lemma \ref{geremia}, we get \equ{periferia1}. \qed
Now, assume that, together with \equ{K1} and \equ{K2}, it is also $\K\ge \cteot$. Then, recalling  \equ{falla}, 
Theorem~\ref{zelensky} holds.
Thus, recalling \equ{topgun},  observing that 
$$
 \Bik(\loge)\ =\  
 \big\{ \act\in \Bik  \ {\rm s.t.}\  
|\det \partial^2_{\act}\mathtt{h}^i_k(\act)|\leq \fico\big\}
\cup \Bik(\loge,\fico) \,,
$$
by \eqref{gibboso} and \eqref{gusuppotto}
we get
\begin{equation}\label{pesantezza}
\meas\big( (\Bik(\loge)\times\T^n)
\setminus {\cal T}^i_k \big)
\ \leq \ \cteot  (|k|^{2n} \fico)^\lippo
\meas \Bi_k\, +\, \ttcq e^{-\K s/7}\,.
\end{equation}
Now, by \equ{insiemi}, \equ{lazie}, \equ{periferia1},  \equ{pesantezza}, \equ{diametro}, \equ{nadal}  and since $\noruno{k}\le \KO=\K/6$
we get, for a suitable constant\footnote{To get \equ{verace}, use the following: $\e\le\K^{-\gamma}<1$ (compare \equ{torino}); $\meas\big(\Rukt\times\T^n\big)\le c(n)$, as $\Rukt\subset \{y:|y|\le 2\}$; $|\log\loge|=\frac{\lippo}{{\ttCu}(1+\lippo)}\K$; $\meas \Bik \le \itcast$ by \equ{diametro}; \#$\genKO<(2\KO+1)^n$.}
 $\cgotu=\cgotu(n,s,\d,\b, \upxi,\mtt)$,
\beq{verace}
\ts \meas(\cA)\le \cgotu \, \K^{2n} e^{- \K/\cteos}\,,\qquad \cteos:= \max\big\{ 36/s, 2\ttCu/\lippo\big\}\,.
\eeq
Finally, let   
\beq{finally}
\cteo=\cteo(n,s,\d,\b, \upxi,\mtt)\ge 1+\cteos
\eeq
be such that, if 
$\K\geq \cteo$, then $\cgotu \K^{2n} e^{- \K/\cteos}\le e^{-\K/(1+\cteos)}$. Then, if $\K\ge \cteo$, claim (iii) follows, and the 
proof of Theorem \ref{prometeo} is complete.
\qed

\rem \label{tortore}
Notice that $\cT^0$ is a family of maximal primary tori for $\ham$, and so are the families $\cT^{1,k}_i$ for all $k\in\genKO$ and $i=0,2N_k$. On the other hand,   $\cT^{1,k}_i$ for all $k\in\genKO$ and $0<i<2N_k$
are families of {\sl maximal secondary tori} for $\ham$. In particular these families do not bifurcate from integrable tori.
\erem

\subsubsection*{Proofs of Corollaries \ref{prometeobis} and \ref{prometeotris}}

\proof {\bf of Corollary \ref{prometeobis}}
As already pointed out in \S~\ref{omnia}, Corollary~\ref{prometeobis}  follows trivially from
 Theorem~\ref{prometeo} and the measure estimate \equ{lapalisse},
by taking 
$\K:=\cteo|\ln\e|$,
$\bcteo:=1+ (2\pi)^n \cdr \cteo^\gamma$, and $\e_{\rm o}$ so that $\e\K^\g<1$  
 for $\e< \e_{\rm o}$.
\qed

\proof {\bf of Corollary \ref{prometeotris}}
Let $n=2$.
We claim that $\Rd$ in
\equ{sonnosonnosonno} satisfies
\beq{alesia}
\Rd\subseteq \{y\in\real^2:|y|<\e^{a/2}\}\,.
\eeq
Fix $y\in \Rd$. Then $|y|<1$ and there exists
$k\in \genKO$ such that $|y\cdot k|\leq \a/4$ (since $y\notin \Rz$).
Moreover, since $y\notin \Ruk$,  there exists 
$\ell\in \genK\bks\Z k$ such that
\beq{gergovia}
|\pk y\cdot \ell|\leq \frac{6 \a \K}{|k|}\,.
\eeq
Then, $|\proiezione_k y|< \a/(4|k|)\leq\a/4$.
Moreover, since  $\ell\notin \Z k$,
$$
|\proiezione_k^\perp \ell|=\frac{|k_1 \ell_2-k_2 \ell_1|}{|k|}\geq \frac{1}{|k|}
\,,\qquad 
|\proiezione_k^\perp y \cdot \ell |=
|\proiezione_k^\perp y|\, |\proiezione_k^\perp \ell| 
\geq \frac{|\proiezione_k^\perp y|}{|k|}\,,
$$
which implies, by \eqref{gergovia}, $|\proiezione_k^\perp y|\leq 6\a\K$.
In conclusion
\beq{ventilatore}
|y|=|\proiezione_k y + \proiezione_k^\perp y|
< 7\a\K\stackrel{\eqref{athlone}}=7\sqrt\e \K^{12}\,.
\eeq
Now, let $\hat a:=(1-a)/24$ and  
$\K:=1/(\sqrt[12]7 \e^{\hat a})$.
Then, \equ{alesia} follows by \equ{ventilatore}.
\\
Finally, let $\e_{\rm o}<1$ be so small that $\e\K^\g<1$
is satisfied for any $\e<\e_{\rm o}$.
Then, 
by the estimate in Theorem~\ref{prometeo}--(iii), we get
$$
\meas(\cA)\le 
\meas\big(\big(\{\e^{a/2}<|y|<1\}\times\T^n\big)
\,\setminus\,{\cal T}
\big)\ \leq \  
e^{-\frac{1}{7^{1/12} \cteo\e^{\hat a}}}<e^{-\frac{1}{{2\cteo \e^{\hat a}}}}
\,. \qquad  \qedeq. 
$$

\appendix


\small

\section{Proofs of elementary lemmata}\label{appendicite}

\subsection*{Proof of Lemma~\ref{telaviv}}

Assume $f\in \Gns$ for some $s>0$ and 
let $0<\d_0\le 1$ be smaller than $$
\varliminf_{\sopra{\noruno{k}\to+\io}{k\in\gen}} 
|f_k| e^{\noruno{k}s} \noruno{k}^n\stackrel{\equ{P1}}>0\,.$$
Then, there exists $N_0$ such that $|f_k|>\d_0 \noruno{k}^{-n} e^{-\noruno{k}s}$, for any $\noruno{k}\ge N_0$, $k\in\genKO$.
Since $\lim_{\d\to 0} \Nf=+\infty$, there exists $0<\d<\d_0$ such that $\Nf>N_0$. Hence, if $\noruno{k}\ge \Nf$ and $k\in\genKO$, \equ{P1+} holds.
\\
Since $\fproj f$ is, for any $\noruno{k}\le \Nf$, a Morse function with distinct critical values one can, obviously, find a $\b>0$ for which
\equ{P2+} holds. 

\nl
To prove the `if part', we need two lemmata.  The first lemma can also be found in \cite{BCuni} (compare Proposition~1.1 there); for completeness, we reproduce the simple but instructive proof also here.

\begin{lemma}\label{pollaio} 
Let $f\in\Bns$ such that
 \eqref{P1+} holds. Then, 
 for any  $k\in \gen $ with  $ \noruno{k}\geq \Nf$,   
there exists $\sa_k\in[0,2\pi)$ 
so that 
\begin{equation}\label{struzzo}
\fproj f(\sa)=2 |f_k |\big(\cos(\sa+ \sa_k)+F^k_\star(\sa)\big)\,,\quad F^k_{\! \varstar}(\sa):=\frac{1}{2|f_k|}\sum_{|j|\geq 2}f_{jk}e^{\ii j \sa}\,,
\end{equation}
with $F^k_{\! \varstar}\in\hol_1^1$ and
$\modulo F^k_{\! \varstar} \modulo_1\leq 2^{-40}$.
\end{lemma}
\proof 
We write $\fproj f$ as 
\beqno
\fproj f(\sa):= \sum_{j\in\integer\bks\{0\}} f_{jk} e^{\ii j \sa}
= \sum_{|j|=1} f_{jk}e^{\ii j\sa} + \sum_{|j| \ge 2} f_{jk}e^{\ii j\sa} \,,
\eeqno
and, defining $\sa_k\in[0,2\pi)$ so that $e^{\ii \sa_k}= f_k/|f_k|$, one has
$$
\frac{1}{2|f_k|}\sum_{|j|=1}f_{jk}e^{\ii j \sa}=\Re \Big( \frac{f_k}{|f_k|} e^{\ii \sa}\Big)=\Re e^{\ii (\sa+\sa_k)}
=\cos (\sa+\sa_k)\,,
$$
which is equivalent to \equ{struzzo}. Now, since
 $\| f\|_s\leq 1$ (so that $|f_k|\le e^{-\noruno{k}s}$), one finds, for $\noruno{k}\ge \Nf$, 
\beqano
&&\modulo F^k_{\! \varstar} \modulo_1\stackrel{\equ{struzzo}}\leq
\frac{1}{2|f_k|}\sum_{|j|\geq 2}|f_{jk}|e^{|j|}
\stackrel{\equ{P1+}}\leq
\frac{\noruno{k}^n e^{\noruno{k}s}}{2\d}\sum_{|j|\geq 2}|f_{jk}|e^{|j|}
\le
\frac{\noruno{k}^n e^{\noruno{k}s}}{2\d}\sum_{|j|\geq 2}e^{-|j|(\noruno{k}s-1)}
\\
&&\le
\frac{2 e^2 \noruno{k}^n}{\d} \ e^{-\noruno{k}s}
=\frac{2^{n+1} e^2 }{s^n\d} e^{-\frac{\noruno{k}s}2}\ \ \Big(\frac{\noruno{k}s }2\Big)^n e^{-\frac{\noruno{k}s}2}
\le 
\Big(\frac{2n}{e s}\Big)^n\, \frac{2e^2}{\d} \,  e^{-\frac{\noruno{k}s}2}
\le 2^{-40}\,,
\eeqano
where  last inequality follows 
since $\noruno{k}\ge \Nf$ (see \equ{enne}).
\qed

\noindent
Now, assume that \equ{P1+} and \equ{P2+} hold for some $\d\in (0,1]$ and $\b>0$. Then, from \equ{P1+} follows immediately \equ{P1}. 
It remains to prove that \equ{P2} holds for any $k\in \gen$ with $\noruno{k}> \Nf$. 
In view of Lemma~\ref{pollaio},    the thesis follows from the following elementary

\begin{lemma}\label{pennarello}
Let $F\in C^2(\T,\R)$, $\bar\sa$  and $0<\cgot<\frac12$ are such that\footnote{
 $\|F\|_{C^2}:=\max_{0\leq k\leq 2}\sup|F^{(k)}|$. Note that, by Cauchy estimates,  $\|F\|_{C^2}\le 2|F|_1$.}
  $\|F-\cos (\sa+\bar \sa)\|_{C^2}\le \cgot$.
 Then,  $F$ has only two critical points and it is $(1-2 \cgot)$--Morse.
 \end{lemma}
\proof By considering the translated function $\sa\to F(\sa-\bar\sa)$, one can reduce oneself to the case $\sa=0$ (note that  $F$ is $\b$--Morse, if and only if $\sa\to F(\sa-\bar\sa)$ is $\b$--Morse).\\
Thus, set $\bar \sa=0$, and note that, by assumption $|F'|=|F'+\sin\sa-\sin\sa|\ge|\sin \sa|-\cgot $, and, analogously, $|F''|\ge |\cos\sa|-\cgot $. Hence, $|F'|+|F''|\ge|\sin\sa|+|\cos\sa|-2\cgot \ge 1-2\cgot $. Next, let us show that $F$ has a unique strict maximum $\sa_0\in I:=(-\pi/6,\pi/6)$ (mod $2\pi$). Writing $F=\cos\sa+g$, with $g:=F-\cos \sa$, one has that $F'(-\pi/6)=1/2+g'(\pi/6)\ge 1/2 - \cgot >0$, and, similarly $F'(\pi/6)\le -1/2 +\cgot $, thus $F$ has a critical point in $I$, and, since $-F''=\cos\sa -g''\ge \cos\sa-\cgot \ge \sqrt3/2-\cgot >0$, $F$ is strictly concave in $I$, showing that such critical point is unique and it is a strict local minimum. In fact, similarly one shows that $F$ has a second critical point $\sa_1\in (\pi-\pi/6,\pi+\pi/6)$ where $F$ is strictly convex, so that $\sa_1$ is a strict local minimum; but, since
 in the complementary of these intervals $F$ is strictly monotone (as it is easy to check),  it follows that $F$ has a unique global strict maximum and a unique global strict minimum. 
Finally, $F(\sa_0)-F(\sa_1)\ge \sqrt3-2\cgot >1-2\cgot $ and the claim follows. \qed

\subsection*{Proof of  Lemma \ref{container}}
First note that by \eqref{pasqua}, \eqref{cimabue} and 
\eqref{alce}
\beq{gatti}
(1-\lalla)p_1^2- (1+\lalla)2^{-16}\ro^2\leq
\Hpend(p,q_1)
\leq (1+\lalla)(p_1^2+ 2^{-16}\ro^2)\,.
\eeq
By the first inequality in \eqref{gatti}
we have that if $(p,q_1)\in\cM(\hat p)$
then $p_1^2\le 
\frac{\Eb+ 2^{-16}\ro^2(1+\lalla)}{1-\lalla}$,
which is indeed smaller than $(\Ro+\ro/2)^2$
by \eqref{minore} and \eqref{alce}.
This proves the second inclusion in 
\eqref{gonatak}.
\\
By the second inequality in \eqref{gatti}
we have that if $|p_1|\leq \Ro+\ro/3$
then
$\Hpend(p,q_1)
\leq (1+\lalla)\big((\Ro+\ro/3)^2+ 2^{-16}\ro^2\big)$,
which is smaller than $\Eb$
again by \eqref{minore} and \eqref{alce}.
This proves the first inclusion in 
\eqref{gonatak}.
\qed

\subsection*{Proof of  Lemma \ref{schiena}}

Let $(J,\psi)\in D_{\zeta_{} r_{}}\times \T_{\zeta_{} s_{}}^n$. Then there exists $J_0\in D$
such that $|J-J_0|<\zeta_{} r_{}$.
Set 
$$
\ts w(z):=\big(J_0+\frac\zeta{z}(J-J_0),\Re \psi +\frac\zeta{z}\Im\psi \big)\,,
$$
with $\Re \psi:=(\Re \psi_1,\ldots,\Re\psi_n)$, and analogously for $\Im\psi$.
Note that  $w(\zeta_{})=(J,\psi)$, and 
that $w(z)\in D_{r_{}}\times \T_{s_{}}^n$ for every $|z|<1$.
Consider the holomorphic function $G(z):={g}(w(z))$ defined for $|z|<1$.
Then,  $|\Im G|\leq \xi$ for $|z|<1$.
Let $u$ and $v$ be real harmonic functions such that $G(z)=u(x,y)+\ii v(x,y)$, where
$z=x+\ii y$.
Since by hypothesis 
$\sup_{|z|<1}|v|\leq \xi$, by interior estimate of derivatives of for harmonic functions\footnote{See Theorem 2.10 in \cite{GT}.} 
we have that 
$\sup_{|z|\leq 1/2}|v_x|\leq 4\xi$
and
analogously for $v_y$.
By Cauchy--Riemann equations, the same estimate holds for
$u_x$.
Therefore 
$\sup_{|z|\leq 1/2}|G'|=\sup_{|z|\leq 1/2}|u_x+\ii v_x|
\leq 8\xi$.
Since $w(0)=(J_0,\Re\psi)\in D\times T^n$
and ${g}$ is real analytic, we have that $G(0)={g}(J_0,\Re\psi)\in\R.$
Then, for any $0<\zeta\le 1/2$, by the the mean value theorem, we have that 
$$
|\Im {g}(J,\psi)|=|\Im G(\zeta_{})|=|\Im G(\zeta_{})-\Im G(0)|
\leq |G(\zeta_{})-G(0)|\leq 
8\zeta_{}\xi\,. \qedeq
$$

\subsection*{Proof of Lemma \ref{tonino}}
Let us consider first the case $P={\rm I_d}.$ Consider the unitary matrix $U$ diagonalizing $Q$, namely
$U^{-1}QU=\Lambda={\rm diag}_{1\leq j\leq d}\lambda_j$. Note that $|Q|=|\Lambda|=\max_{1\leq j\leq d} |\lambda_j|=\l.$
Then $U^{-1}({\rm I_d}+Q)U=I+\Lambda$ and $\det({\rm I_d}+Q)=\det(I+\Lambda)\geq (1-\l)^d,$ proving the case
$P={\rm I_d}.$
\\
Consider now the general case. Write $P+Q=P^{1/2}({\rm I_d}+P^{-1/2}QP^{-1/2})P^{1/2}.$
Note that, since $P^{-1/2}$ is symmetric, then $P^{-1/2}QP^{-1/2}$ is symmetric too.
Since $|P^{-1/2}QP^{-1/2}|\leq |P^{-1/2}|^2|Q|=|P^{-1}| |Q|$ and 
$\det P^{1/2}=(\det P )^{1/2}$, from the previous case the general case follows. 
\\
The final claim in Lemma~\ref{tonino} follows, as  $(1-\l)^d\geq 1-d\l$.
\eproof

\subsection*{Proof of Lemma~\ref{jensen}}
First observe that the cases $i=0$ and 
$i=2N$ are identical since
$$\bar\act_1^{0}(E)=\bar\act_1^{2N}(E)\,,\qquad
\bar{\mathtt E}^{0}(\act_1)=\bar{\mathtt E}^{2N}(\act_1)\,.
$$
Let us then consider the case $i=2N$. By definition,
\begin{equation}\label{bolsena}
 \bar \act_{1}^{2N}(E)
=
\frac{1}{2\pi}\int_0^{2\pi}\sqrt{E-\GO(x)}dx\,,
 \end{equation}
so that, by Jensen's inequality,
 $$
 (2\partial_E\bar \act_1^{i}(E))^3=
 \Big(\frac{1}{2\pi}\int_0^{2\pi}\frac1{\sqrt{E-\GO(x)}}dx\Big)^3
 \leq 
 \frac{1}{2\pi}\int_0^{2\pi}\frac1{(E-\GO(x))^{3/2}}dx
 =-4 \partial^2_{E}\bar \act_1^{i}(E)\,,
 $$
 and the claim follows by \eqref{sonnifera}.
 \qed

\subsection*{Proof of Lemma \ref{geremia}}
Consider first the case 
 $\noruno{k}\ge \Nf$.
 By Lemma~\ref{pollaio}, $F^k_{\! \varstar}:=\fproj f/2|f_k|$
  satisfies 
  $$|F^k_{\! \varstar}-\cos(\sa+\sa_k)|_1\leq 2^{-40}\,.$$
Thus,  by Cauchy estimates we get
  $\|F^k_{\! \varstar}-\cos(\sa+\sa_k)\|_{C^2}\leq 2^{-39}$, so that by
  Lemma~\ref{pennarello} it follows that $2N_k=4$.

\nl
For  the case $\noruno{k}\le \Nf$ we need the following elementary observation:
 
\begin{lemma}\label{solero}
If $G$ is $\b$--Morse, then the number $2N$ of its critical points  is bounded by
$\pi\sqrt{2\max_\R |G''|/\b}$.
\end{lemma}
\proof
If $\sa_i$ and $\sa_j$ are different critical points of $G$, then, by Taylor expansion at order two and by \equ{ladispoli}  one has
$\b\le|G(\sa_i)-G(\sa_j)|\le \frac12 (\max_\R|G''|) |\sa_i-\sa_j|^2$, which implies that the minimal distance between two critical points is at least $\sqrt{2\b/\max_\R |G''|}$, from which the claim follows.
\eproof
Now, by \equ{P2+} we know that $\fproj f$ is $\b$--Morse, and  since $\|f\|_s\leq 1$
  we have
 $\dst
 \sup_\R|(\fproj f)''|\leq
 \sum_{j\neq 0} |f_{jk}| j^2\le 
\sum_{j\neq 0} e^{-|j|}j^2<4$.
Then, by Lemma \ref{solero}, the  claim follows also in this case.   \qed



\begin{thebibliography}{99}     


\footnotesize

\bibitem{A63}
V. I. Arnold, \textit{Proof of A. N. Kolmogorov's theorem on the conservation of conditionally periodic motions with a small variation in the Hamiltonian}, Russian Math. Surv, 1963, vol. 18, no 5, pp. 9--36.

\bibitem{A63b}
V. I. Arnold, \textit{Small denominators and problems of stability of motion in classical and celestial mechanics}. Usp. Mat. Nauk 18, No.6, 91--192 (1963). Engl. transl.: Russ. Math. Surv. 18, No.6, 85--191 (1963)

\bibitem{A}
V. I. Arnol'd, {\sl Instability of dynamical systems with many degrees of freedom}, Dokl. Akad. Nauk SSSR, 156 (1964), pp. 9--12

\bibitem{AKN1} V.~I. Arnol'd, V.~V.  Kozlov, and A.~I. Neishtadt.
Dynamical systems. III.
Translated from the 1985 Russian Edition by A. Iacob. Encyclopaedia of Mathematical Sciences, 3. Springer--Verlag, Berlin, 1988. xiv+291 pp.

\bibitem{AKN}
V.~I. Arnol'd, V.~V.  Kozlov, and A.~I. Neishtadt.
{\sl Mathematical aspects of classical and celestial mechanics},
volume~3 of Encyclopaedia of Mathematical Sciences.
Springer-Verlag, Berlin, third edition, 2006.
[Dynamical systems. III], Translated from the 2002 Russian original by E. Khukhro

\bibitem{BFS} D. Bambusi, A. Fus\`e, M. Sansottera. 
{\sl Exponential stability in the perturbed central motion}.
Regular and Chaotic Dynamics, 23, 821--841 (2018)

\bibitem{BKZ}
P. Bernard, V. Kaloshin, and K. Zhang, {\sl Arnol'd diffusion in arbitrary degrees of freedom and 3--dimensional normally hyperbolic invariant cylinders}, Acta Math., 217:1 (2016), pp. 1--79

\bibitem{BBB}
M. Berti, L. Biasco, and  P. Bolle, {\sl Drift in phase space: a new variational mechanism with optimal diffusion time}. J. Math. Pures Appl. (9) 82 (2003), no. 6, 613--664

\bibitem{BClin} L. Biasco, and L. Chierchia. 
{\sl On the measure of Lagrangian invariant  tori in nearly--integrable mechanical systems}. 
Atti Accad. Naz. Lincei Rend. Lincei Mat. Appl. {\bf 26} (2015), 1--10 

\bibitem{BCKAM} L. Biasco, and L. Chierchia. 
{\sl Explicit estimates on the measure of  primary KAM tori}.
Ann. Mat. Pura Appl. (4), 2017, 261--281

\bibitem{BCnonlin} L. Biasco, and L. Chierchia. 
{\sl On the topology of nearly-integrable Hamiltonians at simple resonances}.
Nonlinearity 33 (2020) 3526--3567

\bibitem{BC-DCDS} L. Biasco, and L. Chierchia. 
{\sl On the measure of KAM tori in two degrees of freedom}.
Discrete and Continuous Dynamical Systems, Volume 40, Number 12 (2020) 
6635-6648

\bibitem{BClin2} L. Biasco, and L. Chierchia. 
{\sl Quasi--periodic motions in generic nearly--integrable mechanical systems}, 
Atti Accad. Naz. Lincei Cl. Sci. Fis. Mat. Natur. 33 (2022), no. 3, pp. 575--580,
\underline{\href{https://arxiv.org/pdf/2206.01055.pdf}{\sl arXiv:2206.01055}}


\bibitem{BCaa23}  L. Biasco, and L. Chierchia. 
{\sl Complex Arnol'd--Liouville maps}, Regular and Chaotic Dynamics, 2023. \\
\underline{\href{https://link.springer.com/article/10.1134/S1560354723520064}{\sl doi 10.1134/S1560354723520064}}. 
\underline{\href{https://arxiv.org/abs/2306.00875v1}{\sl arXiv:2306.00875v1}}

\bibitem{BCuni} L. Biasco, and L. Chierchia. 
{\sl Global  properties of  generic real--analytic nearly--integrable Hamiltonian systems},
\underline{\href{https://arxiv.org/abs/2306.13527}{\sl 	arXiv:2306.13527}}



\bibitem{CC} A. Celletti, and L. Chierchia
{\sl A Constructive Theory of Lagrangian Tori and Computer-Assisted Applications}, in Dynamics Reported, Dynam. Report. Expositions Dynam. Systems (N. S.), vol. 4, Berlin: Springer, 1995, pp. 60--129

\bibitem{CdL22}
Q. Chen; R. de la Llave, 
{\sl Analytic genericity of diffusing orbits in a priori unstable Hamiltonian systems}. 
Nonlinearity 35 (2022), no. 4, 1986--2019

\bibitem{CG} L. Chierchia, and G. Gallavotti.
{\sl Drift and diffusion in phase space}
Ann. Inst. Henri Poincar\'e, Phys. Th\'eor., 60, 1--144 (1994).
See also: Erratum, Ann. Inst. Henri Poincar\'e, Phys. Th\'eor., 68, no. 1, 135 (1998)

\bibitem{CK} L. Chierchia, and C.E. Koudjinan.
{\sl V.I. Arnol'd's ``Global'' KAM Theorem and geometric measure estimates}.
Regular and Chaotic Dynamics, 2021, Vol. 26, no. 1, 61--88

\bibitem{DH}
A. Delshams and G. Huguet. {\sl Geography of resonances and Arnol'd diffusion in a priori unstable Hamiltonian
systems}. Nonlinearity, 22(8):1997--2077, 2009


\bibitem{DLS}
A. Delshams, R. de la Llave, and T.M. Seara, 
{\sl Instability of high dimensional Hamiltonian systems: multiple resonances do not impede diffusion}.
Adv. Math. 294 (2016), 689--755

\bibitem{KAMstory} H. S. Dumas, {\sl The KAM Story}, 
World Scientific, 2014

 
\bibitem{E} H. Eliasson, {\sl Perturbation of linear quasi-periodic systems. } Dynamical Systems and Small Divisors (Cetraro, Italy, 1998), 1-60, Lect. Notes Math.
1784, Springer, 2002

\bibitem{GT}
D. Gilbarg, and N. S. Trudinger 
{\sl Elliptic partial differential equations of second order}. Springer--Verlag,
Berlin--New York, 1977. Grundlehren der Mathematischen Wissenschaften, Vol. 224

\bibitem{GLS}
M. Gidea,  R. de la Llave, and T. M. Seara, 
{\sl A general mechanism of diffusion in Hamiltonian systems: qualitative results} 
Comm. Pure Appl. Math. 73 (2020), no. 1, 150--209

\bibitem{KZ}
V. Kaloshin and K. Zhang. {\sl Arnol'd diffusion for smooth systems of two and a half degrees of freedom}, volume
208 of Annals of Mathematics Studies. Princeton University Press, Princeton, NJ, 2020

\bibitem{Ko54} A. N. Kolmogorov, \textit{On the conservation of conditionally periodic motions under small perturbation of the Hamiltonian}, Dokl. Akad. Nauk. SSR 98 (1954), 527-530. English translation by Helen Dahlby in: Stochastic behavior in classical and quantum Hamiltonian systems, Volta Memorial conference, Como, 1977, Lecture Notes in Physics, 93, Springer, 1979.

\bibitem{Laz} V.F. Lazutkin, {\sl Concerning a theorem of Moser on invariant curves. Probl. Dyn. Theory Propag. Seism. Waves 14, 109--120. Leningrad: Nauka. (1974) (Russian)}

\bibitem{LMS} L. Lazzarini, J-P. Marco, and D. Sauzin, {\sl 
Measure and capacity of wandering domains in Gevrey near-integrable 
exact symplectic systems}, Memoirs of the AMS, {\bf 257}, (2019),  1-110



\bibitem{Ma}
J. N. Mather,
{\sl Arnol'd diffusion by variational methods}. Essays in mathematics and its applications, 271--285, Springer, Heidelberg, 2012

\bibitem{MNT} A.G. Medvedev, A.I.  Neishtadt, D.V. Treschev, 
{\sl Lagrangian tori near resonances of near--integrable Hamiltonian systems}, 
Nonlinearity, {\bf 28}:7 (2015), 2105--2130

\bibitem{Mo62}  J. K. Moser, \textit{On invariant curves of area-preserving mappings of an annulus}, Nach. Akad. Wiss. G\"ottingen, Math. Phys. Kl. II {\textbf 1} (1962), 1-20

\bibitem{Mo68}
J. Moser: \textit{Lectures on Hamiltonian Systems}. Mem. Am. Math. Soc. 81. American Mathematical Society, Providence, R I. (1968), 60 pp.

\bibitem{Nei}
A. I. Neishtadt, 
{\sl Estimates in the Kolmogorov theorem on conservation of conditionally periodic motions}.
J. Appl. Math. Mech. 45 (1981), no. 6, 766--772

\bibitem{P82}
J. P\"oschel,
{\sl Integrability of Hamiltonian systems on Cantor sets}.
Comm. Pure Appl. Math. 35 (1982), no. 5, 653--696

\bibitem{P93}
J. P\"oschel,  
{\sl Nekhoroshev estimates for quasi--convex Hamiltonian systems}. 
Math. Z. {\bf 213}, pag. 187 (1993)

\bibitem{pyartli}
 A. S.  Pyartli,{\sl Diophantine approximations on submanifolds of Euclidean space}. Functional
Analysis Appl. 3, 303--306 (1970) (Translation from Russian, Funkts. Analys. Prilozh. 3, 59--62,  1969)

\bibitem{Sva} A.A. Svanidze, {\sl Small perturbations of an integrable dynamic system with integral invariant}. Proc. Steklov Inst. Math. 2, 127--151 (1981)

\bibitem{T}
D. Treschev. {\sl Arnol'd diffusion far from strong resonances in multidimensional a priori unstable Hamiltonian
systems}. Nonlinearity, 25(9):2717--2757, 2012

\bibitem{Z}
Zhang, Ke 
{\sl Speed of Arnol'd diffusion for analytic Hamiltonian systems}
Invent. Math. 186 (2011), no. 2, 255--290.







\end{thebibliography}
\end{document}